\documentclass[sn-mathphys-num]{sn-jnl}

\usepackage{graphicx}%
\usepackage{multirow}%
\usepackage{amsmath,amssymb,amsfonts}%
\usepackage{amsthm}%
\usepackage{mathrsfs}%
\usepackage[title]{appendix}%
\usepackage{xcolor}%
\usepackage{textcomp}%
\usepackage{manyfoot}%
\usepackage{booktabs}%
\usepackage{algorithm}%
\usepackage{algorithmicx}%
\usepackage{algpseudocode}%
\usepackage{listings}%
\usepackage{tikz}
\usetikzlibrary{positioning}
\usetikzlibrary{cd}
\usetikzlibrary{arrows}
\usetikzlibrary{decorations.pathmorphing, fit}
\definecolor{canaryyellow}{rgb}{1.0, 0.94, 0.0}
\definecolor{brightgreen}{rgb}{0.4, 1.0, 0.0}
\definecolor{jazzberryjam}{rgb}{0.65, 0.04, 0.37}

\theoremstyle{thmstyleone}%
\newtheorem{theorem}{Theorem}
\newtheorem{proposition}{Proposition}

\theoremstyle{thmstyletwo}%

\theoremstyle{thmstylethree}%
\newtheorem{definition}{Definition}%
\newtheorem{obs}{Observation}

\begin{document}

\title[Theoretical results for Perfect Location signed Roman domination problem]{Theoretical results for Perfect Location signed Roman domination problem}


\author*[1]{\fnm{Bojan} \sur{Nikoli\'c}}\email{bojan.nikolic@pmf.unibl.org}

\author[1]{\fnm{Milana} \sur{Grbi\'c}}\email{milana.grbic@pmf.unibl.org}
\equalcont{These authors contributed equally to this work.}

\author[1]{\fnm{Dragan} \sur{Mati\'c}}\email{dragan.matic@pmf.unibl.org}
\equalcont{These authors contributed equally to this work.}

\affil*[1]{\orgdiv{Faculty of Natural Sciences and Mathematics}, \orgname{University of Banja Luka}, \orgaddress{\street{Mladena Stojanovi\'ca 2}, \city{Banja Luka}, \postcode{78000},  \country{Bosnia and Herzegovina}}}


\abstract{The study of Roman domination has evolved to encompass a variety of challenging extensions, each contributing to the broader understanding of domination problems in graph theory. This paper explores the Perfect Location Signed Roman Domination (PLSRD) problem, a novel combination of the Perfect Roman, Locating Roman, and Signed Roman Domination paradigms. In PLSRD, each weak vertex, assigned  the label  -1, must be protected by exactly one strong vertex, with additional limitation that two weak vertices cannot share the same strong vertex, while the total sum of labels in the closed neighborhood of each vertex must remain positive.
	
	This paper provides exact values for the PLSRD number in several well-known graph classes, including complete graphs, complete bipartite graphs, wheels, paths, cycles, ladders, prism graphs, and $3\times n$ grids. Additionally, we establish a lower bound for a general 3 regular graph, as well as the upper bounds for flower snarks graphs, highlighting the intricate interplay between the PLSRD constraints and the structural properties of these graph families.}




\maketitle

\section{Introduction}
In contemporary graph theory, domination problems play a crucial role in understanding and optimizing network structures. 
One significant aspect of this field is the Roman domination problem (RDP), which was inspired by the military strategies employed during ancient Rome when Emperor Constantine the Great developed military tactics for defending Roman provinces. The basic Roman domination problem involves labeling each vertex with values 0, 1, or 2, ensuring that every ``weak'' vertex, i.e. vertex labeled with 0, has a neighboring vertex labeled with 2, all while minimizing the total sum of the labels. Since the early works concerning Roman domination \cite{stewart1999defend,revelle2000defendens,cockayne2004roman}, dozens of variations of this problem have emerged in the literature. Of particular interest in this paper are three variants of the basic problem. The signed Roman domination \cite{ahangar2014signed} focuses on the labeling of vertices with -1, 1, or 2, with the constraints that each ``weak'' vertex, which is in this case a vertex labeled with -1 must have a corresponding ``strong'' neighbor labeled with 2, and the sum of labels in the closed neighborhood of each vertex must be positive. The perfect variant of RDP, called Perfect Roman domination \cite{HENNING2018235,banerjee2019perfect}, imposes an additional constraint that every ``weak'' vertex must have exactly one ``strong'' neighboring vertex. Lastly, an additional condition considered in the problem addressed in this paper is the so-called locating property \cite{slater1987domination,slater1988dominating}, which, in the context of protecting vertices, requires that two distinct weak vertices cannot share the same set of strong neighbors.


While this approach has yet to be tested in practical applications, the perfect locating signed Roman domination problem offers promising avenues for resource allocation and strategic placement in various fields. For instance, in emergency management, this model could optimize the placement of rescue teams or medical units in urban areas. By treating weak vertices as regions that require coverage and strong vertices as the locations of emergency responders, this method would ensure that each vulnerable area is serviced by exactly one responder, minimizing response time and maximizing efficiency. Additionally, in sensor networks, this framework could be employed to position sensors in a way that each monitored zone is supported by a dedicated sensor while preventing overlap in coverage. Since each weak vertex is associated with a unique strong vertex, this configuration suggests that a graph can effectively defend itself against a simultaneous attack on all vertices, ensuring that every weak vertex is adequately protected.
\subsection{Basic notation}

In this paper, we assume that all graphs are finite and contain neither loops nor parallel edges. Let $G=(V,E)$ represent a graph, where $V$ is set of vertices and  $E$ is set of edges. By  $N(v) = \{u \in V \mid vu \in E\}$ we denote the open neighborhood of the vertex $v$. The closed neighborhood is denoted by $N[v] = N(v)\cup\{v\}$. The smallest degree of any vertex in the graph is denoted by $\delta (G)$.  Let $T$ be a tree. A \textit{leaf} vertex of a tree $T$ is a vertex that has a degree of 1, meaning it is connected to exactly one other vertex from $T$. A \textit{support} vertex is a vertex in a tree that is adjacent to one or more leaf vertices.

\subsection{Literature review}

The \textit{Perfect Locating Signed Roman Domination Problem (PLSRDP)} is introduced in~\cite{akwu2026perfect} as a problem that addresses multiple aspects of vertex protection. As already mentioned, it integrates the concepts of perfect, locating, and signed Roman domination. In the following paragraphs, we will provide an overview of the relevant literature pertaining to these three problems.

The perfect domination set of a graph $G = (V, E)$ is a subset $S$ of $V$ such that every vertex in $V$ is either in $S$ or is adjacent to exactly one vertex in $S$~\cite{weichsel1994dominating}. In the field of Roman domination problems, the concept of perfect domination ensures that each ``weak'' vertex has exactly one ``strong'' protector. The perfect Roman domination problem (PRDP) was formally introduced in~\cite{HENNING2018235}, where the authors showed that the upper bound for the perfect Roman domination number (PRDN) for a tree with $n$ vertices does not exceed $\frac{4n}{5}$. This result was further improved in~\cite{darkooti2019perfect}, by proving that for every tree $T$ with $n$ vertices, with $l$ leaves and $s$ support vertices, PRDN is less or equal to  $\frac{4n}{5} - \frac{s}{n}$.  PRDP is proven to be NP-complete for chordal, planar, and bipartite graphs \cite{banerjee2019perfect}. 
In recent years, the problem of determining PRDN for different classes of graphs has been extensively studied in the literature. The authors in~\cite{almulhim2022perfect} showed that PRDN of the Cartesian product of a path and a path, a path and a cycle, and a cycle and a cycle does not exceed two-thirds of the total number of vertices. Additionally, the problem has been investigated for regular graphs~\cite{henning2018perfect}, where upper bounds have been established for cubic graphs and \( k \)-regular graphs with $k \geqslant 4$.

The locating problem was first introduced in~\cite{slater1988dominating}. Consider a dominating set $S$ in a graph $G$. For each vertex $u \in V \setminus S$, let $S(u)$ represent the set of vertices in $S$ that are adjacent to $u$. The set $S$ is referred to as a locating dominating set (LD-set) if, for any two distinct vertices $u$ and $w$ in $V \setminus S$, it holds that $S(u) \neq S(w)$. The smallest possible size of an LD-set in $G$ is called the location domination number of $G$. Building on the concept of locating domination, the authors in~\cite{rad2019locating} extended the idea to Roman dominating functions by introducing an additional condition. They considered Roman dominating functions $f = (V_0, V_1, V_2)$ with the requirement that for each vertex $x \in V_0$, the set $N(x) \cap V_2$ must be unique. In other words, any two distinct vertices $x, y \in V_0$ are distinguished by the existence of a vertex $v \in V_2$ such that $|N(v) \cap \{x, y\}| = 1$. Such a function is called a locating Roman dominating function (LRDF) if $N(u) \cap V_2 \neq N(v) \cap V_2$ for any two distinct vertices $u, v \in V_0$. The minimum weight of an LRDF in a graph $G$ is referred to as the locating Roman domination number (LRDN). In simpler terms, for any two ``weak'' vertices, their sets of guardians must be different. In~\cite{jafari2017bounds}, the authors investigate LRDN in trees, deriving bounds in terms of the tree's order, number of leaves, and support vertices, and characterizing trees that achieve equality for these bounds. In \cite{guo2020locating}, the authors integrated the concepts of  total domination and locating domination, deriving bounds for the locating-total domination number specifically for grid graphs. Combining the conditions related to perfect Roman domination and locating Roman domination ensures that the sets of guards for two ``weak'' vertices are distinct and singleton. In other words, each ``weak'' vertex has its unique ``strong'' guard, and each ``strong'' guard protects at most one ``weak'' vertex.

In~\cite{ahangar2014signed}, the concept of signed Roman domination in graphs is introduced. Recall that a signed Roman dominating function (SRDF) assigns a value from $\{-1,1,2\}$ to each vertex of a graph, ensuring that the sum of values in the closed neighborhood of each vertex is at least one and that every vertex assigned $-1$ is adjacent to at least one vertex assigned $2$. The goal is to minimize the total weight of such a function across all vertices, which defines the signed Roman domination number (SRDN) of the graph. The same paper establishes upper bounds for SRDN for general graphs and for bipartite graphs, and characterizes the graphs that achieve equality in these bounds. It is also known that determining SRDN is an NP-complete problem, even when the problem is restricted to bipartite and planar graphs~\cite{shao2017signed}. This problem has been extensively studied in recent years, and numerous results have been obtained regarding the SRDN for various types of graphs. Specifically, the signed Roman domination number has been determined for the following graph classes: stars, complete graphs, complete bipartite graphs  $K_{p,p}$, cycles, and paths~\cite{volkmann2016signed}, complete bipartite graphs and wheels~\cite{yan2017signed}, the join of cycles, wheels, fans, and friendship graphs~\cite{behtoei2014signed}, the spider and double star graphs~\cite{hong2020signed}, convex polytopes~\cite{zec2021signed}, and digraphs~\cite{sheikholeslami2015signed}.

By combining the previously introduced problems of perfect, locating, and signed Roman domination, the concept of perfect locating signed Roman domination (PLSRD) arises. This problem can be seen as a natural extension in which each of the earlier concepts contributes to the formulation of the new problem. Specifically, perfect domination ensures that every vertex is uniquely dominated, meaning that each ``weak'' vertex is adjacent to exactly one ``strong'' vertex. Locating domination guarantees that each vertex can be uniquely identified by the set of its adjacent strong vertices. In the context of signed Roman domination, a specific assignment of values from the set $\{-1, 1, 2\}$ is made to the vertices, with the additional requirement that each weak vertex has exactly one strong guard and that the sum of the values in the closed neighborhood of each vertex must be at least one. The PLSRD problem combines all of these conditions, ensuring that each weak vertex has a unique guardian and that the guardians' sets are distinguishable. This leads to a situation where every weak vertex is uniquely dominated, every guardian has exactly one weak vertex to protect, and the signed Roman domination conditions are met with the additional constraints of perfect and locating domination. 

To the best of our knowledge this problem was investigated only in~\cite{akwu2026perfect}, where the upper and lower bounds of the PLSRD function for trees are presented. Additionally, for grid graph $G=(V,E)$, it is shown that $\gamma_{PLSR}(G) \leq \frac{3}{4}|V|$.
However, upon reviewing Theorem 2.3 from the same paper, we would like to highlight a potential clarification. The theorem asserts: ``If \(T\) is a tree of order \(n \geq 3\), then \(\gamma^{P}_{LSR}(T) \leq \frac{3}{4}n\), where \(T\) is not a star of order \(n \geq 5\).'' After careful examination, it seems that the required conditions of the theorem must be stronger. More precisely, the claim might hold if the tree \(T\) does not contain subtrees that are stars of order at least 5. In the current formulation, the authors appear to exclude only trees that are stars of order \(n \geqslant 5\), which is not sufficient for the claim.

\section{Problem definition}
In this section we formally define the perfect locating signed Roman dominating function and the perfect locating signed Roman domination number.

\begin{definition}\label{def:1}
	A \textit{perfect locating signed Roman dominating} (PLSRD) function on a graph $G = (V, E)$ is defined as a function $f: V \rightarrow \{-1, 1, 2\}$ that adheres to the following three conditions:

	\begin{enumerate}
		\item[\textnormal{C1:}] Each vertex $v$ for which $f(v) = -1$ must be adjacent to exactly one vertex $u$ where $f(u) = 2$.
		\item[\textnormal{C2:}] Distinct vertices $v$ and $w$ with $f(v) = f(w) = -1$ cannot have a common neighbor $u$ that satisfies $f(u) = 2$.
		\item[\textnormal{C3:}] For every vertex $v$, the total of the function value at $v$ and the values at its neighboring vertices must be positive, i.e. $\sum_{u \in N[v]} f(u) \geq 1$. 
	\end{enumerate}
\end{definition}

The weight of a PLSRD function $f$ of the graph $G=(V,E)$ is calculated as the sum of its values across all vertices. If we denote the set of vertices labeled with $i$ as $V_i$, $i\in \{-1,1,2\}$, then this can be expressed as $f(V) = |V_1| + 2|V_2| - |V_{-1}|$. Typically, the function $f$ can also be identified as a triple $(V_{-1},V_1,V_2)$.

\begin{definition}
	The perfect locating signed Roman domination number (PLSRDN) of $G$, denoted by $\gamma_{plsR}(G)$, is the minimum weight of all PLSRD functions of $G$.
	
\end{definition}

\section{Main theoretical results}

In this section, we present the primary theoretical contributions of the paper. We begin with two observations.

\begin{obs} \label{obs1} Let $G=(V,E)$ be a graph and $f = (V_{-1},V_1,V_2)$ be a PLSRD function, then
	$|V_{-1}|\leqslant|V_2|$.
\end{obs}
\begin{obs}
	\label{obs:x} Let $G=(V,E)$ be a graph with $n$ vertices and $f = (V_{-1},V_1,V_2)$ be a PLSRD function where $|V_{-1}| = l$. Then $f(V) \geqslant n-l$, while the equality holds when $|V_{-1}| = |V_2|$.
\end{obs}
\begin{proof}
	
	Let	$|V_{-1}| = l$. From Observation~\ref{obs1}, at least $l$ vertices are labeled with 2, while the rest of at most $(n-2l)$ vertices are labeled with at least 1. So
	
	$f(V)  \geqslant -l+2l+(n-2l) = l + (n-2l) = n-l$.
	
	Obviously, the equality holds if $|V_{-1}| = |V_2|$.
\end{proof}

\begin{proposition} If $K_n$ is a complete  graph with $n$  vertices, $n \geqslant 2$, then
	
	$$	\gamma_{plsR}(K_n) = n-1$$
\end{proposition}
\begin{proof}
	Let us denote the vertices with $v_1,v_2,\ldots, v_n$.
	
	\textit{Upper bound.} Let us define the PLSRD function $f$ with  $V_{-1} = \{v_1\}$,  $V_{2} = \{v_2\}$, $V_1 = V\setminus (V_{-1}\cup V_2)$. Obviously, $f(V) = -1+2 + n - 2 = n-1$.

	\textit{Lower bound.} Let $f$ be an arbitrary PLSRD function. Suppose that $|V_{-1}|\ > 1$ and WLOG let $v_1,v_2\in V_{-1}$. Since the graph is complete, both vertices are protected with the same vertex labeled with 2, which is a contradiction with condition~(C2) from problem definition. Therefore, $|V_{-1}|\ \leqslant 1$ and from Observation~\ref{obs:x}, we have that $f(V)\geqslant n-1$, which concludes the proof.
\end{proof}

\begin{proposition} Let $K_{n,p}$ be a complete bipartite graph, $2\leqslant p\leqslant n$. Then 
	
	$$	\gamma_{plsR}(K_{p,n}) = p+n-2$$
\end{proposition}
\begin{proof} 
	Let the set $V$ of vertices of $K_{n,p}$ be partitioned in sets $X$ and $Y$ and let $X=\{u_1,\ldots u_p\}$ and $Y=\{v_1,\ldots v_n\}$.
	
	\textit{Upper bound}.

	We define $V_{-1} = \{u_1,v_1\}$, $V_2 = \{u_2,v_2\}$ and $V_1 = V\setminus (V_{-1}\cup V_2)$.
	Obviously, $f$ is a PLSRD function with two vertices labeled  -1, so $f(V) = p+n-2$.
	
	Therefore, $	\gamma_{plsR}(K_{n,p})\leqslant p+n-2$.

	\textit{Lower bound}.
	
	Let $f$ be an arbitrary PLSRD function. Suppose that $|V_{-1}|>2$, i.e. there are at least three vertices labeled with -1. In that case, at least two such vertices $u,v$ belong to the same partition set. WLOG suppose that both $u,v \in X$. Since $f(u) = -1$, there is a vertex $w\in Y$ labelled with 2, which is adjacent with $u$. However, $w$ is also adjacent with $v$, which is a contradiction with condition (C2) from Definition~\ref{def:1}. Therefore, 
	$|V_{-1}|\leqslant 2$ and because of Observation~\ref{obs:x}, $f(V) \geqslant p+n-|V_{-1}|\geqslant p+n-2$.
	
	Therefore, $	\gamma_{plsR}(K_{n,p})\geqslant p+n-2$, which concludes the proof.
	
\end{proof}

For the sake of completeness, we mention the result for star graphs from \cite{akwu2026perfect}.
\begin{proposition} \cite{akwu2026perfect}  Let $S_n$ be a star graph with $n+1$ vertices, where a single central vertex is adjacent to all other $n$ vertices, and those vertices are not adjacent to each other. Then
	
	$$	\gamma_{plsR}(S_n) = n$$
\end{proposition}

\begin{proposition} Let $W_n$, $n\geqslant 4$ be a wheel graph with $n+1$ vertices, where one central vertex is adjacent to all vertices of an $n$-cycle. Then 
	
	$$	\gamma_{plsR}(W_n) =
	\begin{cases}
		\lceil\frac{n}{2}\rceil + 1, & \text{if } n \equiv 0 \pmod{4}, \\
		\lceil\frac{n+1}{2}\rceil+1, & \text{otherwise} 
	\end{cases}
	$$
\end{proposition}
\begin{proof}
	Let $v$ be the central vertex and let $c_1,c_2,\ldots, c_n$ be the vertices in the cycle.
	
	\textit{Upper bound.}

	Let  $n = 4k+r$, where $r\in\{0,1,2,3\}$. For any such value of $r$, we define the function $f:V\rightarrow \{-1,1,2\}$ as follows.
	
	$f(v) = 1$,
	
	$f(c_{4i+2}) = f(c_{4i+3}) = -1$, $i \in \{0,\ldots, k-1\}$,
	
	$f(c_{4i+1}) = f(c_{4i+4}) = 2$, $i \in \{0,\ldots, k-1\}$.

	In other words, the first $4k$ vertices of the cycle are split into $k$ blocks of 4 consecutive vertices, assigned with 2, -1, -1, 2, respectively. The central vertex is assigned with 1.
	
	For cases $n = 4k+r$, where $r>0$ we assign the rest of vertices as follows:
	
	Case $r=1$: $f(c_{4k+1}) = 1$
	
	Case $r=2$: $f(c_{4k+1})  = f(c_{4k+2}) = 1$
	
	Case $r=3$: $f(c_{4k+1})= 2$, $f(c_{4k+2})= -1$, $f(c_{4k+3})= 1$.
	
	Regardless the value of $r$, it can easily be checked that such defined function $f$ is a PLSRD function. In each case, it holds that $|V_{-1}| = |V_2|$ and we get:
	
	Case $r=0$: $f(V)= 2k+1 = \lceil\frac{n}{2}\rceil + 1$; 
	
	Case $r=1$: $f(V)=(2k+1) + 1 = 2k+2 = \lceil\frac{n+1}{2}\rceil + 1$;
	
	Case $r=2$: $f(V)= (2k+1) + 1 + 1 = 2k+3 = \lceil\frac{n+1}{2}\rceil + 1$ 
	
	Case $r=3$: $f(V)= (2k+1) + 2 + 1 - 1 = 2k+3 = \lceil\frac{n+1}{2}\rceil + 1$ 
	
	which concludes the proof for the upper bound.

	\textit{Lower bound.}
	
	Let $f$ be an arbitrary PLSRD function.
	First, we  analyze the case $f(v)\neq 1$.
	Suppose that $f(v) = -1$. Then there is a vertex $c$ from the cycle, such that $f(c) = 2$. Suppose that there is a vertex $c'$ from the cycle, such that $f(c') = -1$. Since $f$ is a PLSRD function, one of its neighbor vertices is assigned with 2. However, that vertex is also adjacent to $v$, which implies that $v$ has two different neighbors assigned with 2, which is a contradiction with condition (C2) from Definition~\ref{def:1}.
	Suppose now that $f(v) = 2$. In that case, only one vertex from the cycle can be assigned with -1. 	
	In both cases we get that $|V_{-1}|\leqslant 1$, implying that 
	$f(V) \geqslant n$. 
	Since  $n$ is strictly greater than the previously proven upped bound, we conclude that the weight of the function $f$ also exceeds this bound, for $f(v)\neq 1$.

	Therefore, we can now suppose that $f(v) = 1$.

	It should be noted that situations where there are three consecutive vertices in the cycle assigned the value -1, or where two vertices assigned the value -1 share a common neighbor in the cycle, are not possible. Therefore, among any four consecutive vertices in the cycle, at most two can be assigned the value -1. 
	
	Let  $n = 4k+r$, where $r\in\{0,1,2,3\}$.
	
	Case 1. $n = 4k$. Then $|V_{-1}|\leqslant 2k = \frac n 2$. From Observation~\ref{obs:x}, $f(V)\geqslant n +1 -\frac n 2 = \frac n 2 + 1 = 	\lceil\frac{n}{2}\rceil + 1$.
	
	Case 2. $n = 4k+1$. Suppose that $|V_{-1}|\geqslant 2k+1$. Since each vertex assigned with -1 has to have the unique adjacent vertex assigned with 2, we have  $|V_{2}|\geqslant 2k+1$. Now, we have $|V| = |V_{-1}|+|V_1|+|V_2| \geqslant 2k+1+1+2k+1 = 4k+3>4k+2 = n+1 = |V|$, which is a clear contradiction. Therefore,  $|V_{-1}|\leqslant 2k$ and from Observation~\ref{obs:x}, $f(V)\geqslant n + 1-\frac {n-1} 2  = \frac {n+1} 2 + 1 = \lceil\frac{n+1}{2}\rceil+1$.
	
	Case 3. $n = 4k+2$. Suppose that $|V_{-1}|\geqslant 2k+1$. This implies that  $|V_{2}|\geqslant 2k+1$, i.e. $|V_{-1}| = |V_2| = 2k+1$. From this and the fact that $f$ is a PLSRD function, one can easily conclude that the cycle vertices of the wheel are divided into $2k+1$ pairs of vertices assigned with -1 and 2, forming the array with the assigned values
	$-1\ 2\ 2\ -1\ -1\ 2\ldots$
	However, since the number of pairs are odd, the last vertex in this array is assigned with 2. Since the vertices forms the cycle, that vertex is adjacent with the first vertex which is assigned with -1. This leads to a contradiction, because the first vertex of the array (assigned with -1) has two neighbor vertices assigned with 2. Therefore, $|V_{-1}|\leqslant 2k$ and from Observation~\ref{obs:x}, $f(V)\geqslant n + 1-\frac {n-2} 2 = \frac {n+2} 2 + 1 = \lceil\frac{n+1}{2}\rceil+1$.
	
	Case 4. $n = 4k+3$. Similarly as in Case 2, suppose that  $|V_{-1}|\geqslant 2k+2$. This implies that  $|V_{2}|\geqslant 2k+2$ and $|V| = |V_{-1}|+|V_1|+|V_2| \geqslant 2k+2+1+2k+2 = 4k+5>4k+4 = n+1 = |V|$, which is a contradiction. Therefore,  $|V_{-1}|\leqslant 2k+1$ and from Observation~\ref{obs:x}, $f(V)\geqslant n + 1-\frac {n-1} 2  = \frac {n+1} 2 + 1 = \lceil\frac{n+1}{2}\rceil+1$.

\end{proof}

In the following two propositions, we establish a connection between the PLSRD problem and the maximum 2-packing problem. The \textit{maximum 2-packing problem} involves finding the largest set of vertices in a graph such that there are no two vertices in this set that are within a distance of 2 from each other. In other words, a set $S$ is a 2-packing set if for any two different vertices $v$ and $u$ from $S$, there are no path  between $v$ and $u$  within the graph that is of length 2 or less. The cardinality of a largest 2-packing set of $G$ is the 2 - packing number.

\begin{proposition} 
	Let $G=(V,E)$ be such graph that $\delta (G)\geqslant 2$, then $\gamma_{plsR}(G)\leqslant n - |S|$, where $S$ is maximum $2$-packing set of graph $G$.
\end{proposition}
\begin{proof}
	Let $S$ be a maximum $2$-packing set of graph $G$. We define $PLSRD$ function $f=(V_{-1},V_1,V_2)$ as follows. Let $V_{-1}=S$. 
	Since $\delta(G)\geqslant 2$, the graph does not have isolated vertices. For each vertex $v\in V_{-1}$, we choose exactly one vertex $u_v\in N(v)$.  The set of all such chosen vertices $u_v$ forms the set $V_2$, defined as $V_2 = \{u_v : v \in V_{-1}\}$. Finally, we define $V_1 = V \setminus (V_{-1} \cup V_2)$.
	
	To verify that the function $f$ is a PLSRD function, we proceed as follows:
	
	Let $v\in V_{-1}$. From the definition of the set $V_2$, $v$ has exactly one neighbor $u_v$ labeled with 2, thus satisfying condition (C1) from Definition~\ref{def:1} is satisfied. 
	Suppose that $u_v$ has another adjacent vertex $v_1\in V_{-1}$. This would imply that $v$ and $v_1$ share a common neighbor, contradicting the fact that $S$ is the $2$-packing set. Therefore, condition (C2) from Definition~\ref{def:1} is also satisfied.
	
	Condition (C3) from Definition~\ref{def:1} is evidently satisfied as well.  Specifically, from the definition of the set $V_{-1}$, each vertex $v\in V_{-1}$ has no adjacent vertices in $V_{-1}$ and  exactly one adjacent vertex in $V_2$.  Thus, we have,  $f(v) + \sum_{u \in N(v)} f(u) \geqslant -1+2 = 1$. Similarly, each vertex $v\in V_2$, the expression  $f(v) + \sum_{u \in N(v)} f(u) \geqslant 2+(-1) = 1$ holds. Finally, consider any vertex $v\in V_1$.  Since  $V_{-1} = S$, $v$ can have at most one adjacent vertex from $V_{-1}$. Given that $\delta(G)\geqslant 2$, $v$ has at least two neighbors. Consequently, $f(v) + \sum_{u \in N(v)} f(u) \geqslant 1+(-1)+1 = 1$.
	
	From the definition of the function $f$ it is evident that the number of vertices labeled -1 and 2 is equal. Therefore, from Observation~\ref{obs:x}, we have
	$f(V)  = |V| - |S|$, which implies the statement of the proposition.
	
\end{proof}

The following proposition connects the PLSRD number and maximum $2$-packing number in trees. 
\begin{proposition}
	Let $T$ be a tree with $n$ vertices, then $\gamma_{plsR}(T)\leqslant n - |S|$, where $S$ is maximum $2$-packing set of tree $T$.
\end{proposition}

\begin{proof}
	Let $S$ be a maximum 2-packing set of the tree $T$. We construct the function $f$ as follows.
	For each $v\in S$ we establish the following three rules:
	\begin{itemize}
		\item [i.]$v$ is a leaf vertex. Then, $f(v) = -1$ and $f(s)=2$, where $s$ is its support vertex.
		\item[ii.] $v$ is a support vertex. Then it has at least one leaf vertex $u$. We assign $f(v) = 2$ and $f(u) = -1$. All other neighbors of $v$ are assigned the value 1.
		\item[iii.] If $v$ is neither leaf nor support vertex,  then $f(v) = -1$. We assign one of its neighbors the value 2 and all other neighbors the value 1.
	\end{itemize}
	The remaining vertices (if any) are assigned the value 1. An example of the construction of such a function, based on a 2-packing set, is shown in Figure~\ref{fig:tree}.
	\begin{figure}
		\centering
		\begin{tikzpicture}[scale=0.45, line cap=round,line join=round,>=triangle 45,x=1cm,y=1cm]
			\draw [line width=1pt] (-11,0) node[above]{$-1$}-- (-9,0);
			\draw [line width=1pt] (-9,0) node[above]{$2$}-- (-7,0);
			\draw [line width=1pt] (-9,0) -- (-9,-2) node[left]{$1$};
			\draw [line width=1pt] (-7,0) node[above]{$1$}-- (-5,0);
			\draw [line width=1pt] (-5,0) node[above]{$-1$}-- (-3,0)node[above]{$2$};
			\draw [line width=1pt] (-3,0) -- (-2,-2)node[above]{$1$};
			\draw [line width=1pt] (-3,0) -- (-2,2)node[left]{$1$};
			\draw [line width=1pt] (-2,2) -- (0,2)node[above]{$2$};
			\draw [line width=1pt] (0,2) -- (2,2)node[above]{$-1$};
			\draw [line width=1pt] (-2,-2) -- (0,-2)node[above]{$2$};
			\draw [line width=1pt] (0,-2) -- (2,-2)node[above]{$-1$};
			\draw [line width=1pt] (0,-2) -- (2,0)node[above]{$1$};

			\draw [line width=1pt] (-5,0) -- (-4,-2) node[above]{$1$};
			\draw [line width=1pt] (-5,0)-- (-6,-2) node[above]{$1$};
			\draw [line width=1pt] (-4,-2) -- (-5,-4) node[above]{$2$};
			\draw [line width=1pt] (-4,-2) -- (-3,-4) node[above]{$2$};
			\draw [line width=1pt] (-6,-2) -- (-7,-4) node[above]{$2$};
			\draw [line width=1pt] (-7,-4) -- (-7.5,-5.5) node[left]{$-1$};
			\draw [line width=1pt] (-5,-4) -- (-5.5,-5.5) node[left]{$-1$};
			\draw [line width=1pt] (-3,-4) -- (-3.5,-5.5) node[left]{$-1$};
			\draw [line width=1pt] (-3,-4) -- (-2,-5.5) node[above]{$1$};
			\draw [line width=1pt] (-3,-4) -- (-1,-4) node[right]{$1$};
			\draw [fill=black] (-5,0) circle (5.5pt);
			\draw [fill=black] (-7,0) circle (2.5pt);
			\draw [fill=black] (-9,-2) circle (2.5pt);
			\draw [fill=black] (-9,0) circle (2.5pt);
			\draw [fill=black] (-11,0) circle (5.5pt);
			\draw [fill=black] (-4,-2) circle (2.5pt);
			\draw [fill=black] (-6,-2) circle (2.5pt);
			\draw [fill=black] (-3,0) circle (2.5pt);
			\draw [fill=black] (-7,-4) circle (2.5pt);
			\draw [fill=black] (-5,-4) circle (2.5pt);
			\draw [fill=black] (-3,-4) circle (2.5pt);
			\draw [fill=black] (-7.5,-5.5) circle (5.5pt);
			\draw [fill=black] (-5.5,-5.5) circle (5.5pt);
			\draw [fill=black] (-3.5,-5.5) circle (5.5pt);
			\draw [fill=black] (-2,-5.5) circle (2.5pt);
			\draw [fill=black] (-1,-4) circle (2.5pt);
			\draw [fill=black] (-2,-2) circle (2.5pt);
			\draw [fill=black] (-2,2) circle (2.5pt);
			\draw [fill=black] (0,2) circle (5.5pt);
			\draw [fill=black] (2,2) circle (2.5pt);
			\draw [fill=black] (0,-2) circle (5.5pt);
			\draw [fill=black] (2,-2) circle (2.5pt);
			\draw [fill=black] (2,0) circle (2.5pt);
			
		\end{tikzpicture}
		\caption{An illustration of labeling of a tree, based on  a 2-packing set $S$. Enlarged vertices belong to $S$. Four leaf vertices from $S$ are labeled with -1. At the right side of the  graph two support vertices from $S$ are labeled with 2, while one leaf vertex per each is labeled with -1. Finally, the remaining vertex from $S$ in the central part of the graph is labeled with -1, one its neighbor with 2, while others are labeled with 1. The rest of vertices is labeled with 1.}\label{fig:tree}
	\end{figure}
	
	To verify that the function $f$ is well-defined, note that the vertices assigned with -1 (items (i) and (iii)) are either vertices from $S$ or a leaf vertex adjacent  to only one vertex from $S$ (item (ii)), which is a support vertex. Obviously, the distance between any two vertices assigned with -1 is at least 3.
	Therefore, any of the three rules defined above perform assignments exclusively to the vertex $v\in S$ and its neighbors, without affecting any other vertex from $S$ or their respective neighbors.
	
	From this, it can also be concluded that each vertex assigned with -1 has exactly one adjacent neighbor assigned with 2. Additionally, since the distance between any two vertices from $V_{-1}$ is at least three, no pair of such vertices share a common neighbor assigned with 2.
	
	It remains to examine whether the sum over the closed neighborhood of each vertex is positive.
	
	Case 1: $v\in V_{-1}$. 
	
	If $v$ is a leaf, then its only neighbor is assigned the value 2, so $\sum_{u\in N[v]} f(u)= -1+2 = 1$. 
	If  $v$ has at least two neighbors, by the definition of the function $f$, one of them is assigned the value 2 and all others are assigned the value 1. Therefore  $\sum_{u\in N[v]}f(u) \geqslant -1+2 + 1 = 1$.
	
	Case 2: $v\in V_2$. It has already been shown that each vertex in $V_2$ has exactly one adjacent vertex from $V_{-1}$, so $\sum_{u\in N[v]}f(u) \geqslant -1+2=1$.
	
	Case 3: $v\in V_1$. If $v$ is a leaf, its only neighbor is not assigned the value -1, so, $\sum_{u\in N[v]}f(u) \geqslant 2$. If $v$ has at least two neighbors, at most one of them is assigned the value -1, therefore   $\sum_{u\in N[v]}f(u) \geqslant 1+(-1) + 1 = 1$.
	
	Therefore, $f$ is a PLSRD function. Using  $|V_{-1}| = |V_2| = |S|$ and Observation~\ref{obs:x} we have $f(V) = n-|S|$, which completes the proof of the proposition.

\end{proof}

In the following theorem, we propose a lower bound for the PLSRDN of an arbitrary 3-regular graph. To facilitate the proof, we first state and prove a useful observation.

\begin{obs}\label{obs:3reg}
	Let $G$ be a 3-regular graph and let $f = (V_{-1},V_1,V_2)$ be a PLSRD function. Then the inequality $2|V_1|\geqslant |V_{-1}|$ holds.
\end{obs}
\begin{proof}
	Every vertex labeled $-1$ must have at least one neighbor labeled $1$, while  no vertex labeled $1$ can have more than two neighbors labeled $-1$. It follows that the number of vertices in $V_{-1}$ cannot exceed twice the number of vertices in $V_1$.
\end{proof}
\begin{theorem}\label{thm:3reg}
	Let $G$ be a 3-regular graph with $n$ vertices. Then  the inequality
	$$\gamma_{plsr}(G)\geqslant \left\lceil\frac{3}{5}n\right\rceil$$
	holds.
\end{theorem}
\begin{proof}
	Let $f = (V_{-1},V_1,V_2)$ be an arbitrary PLSRD function and let $|V_{-1}| = l$. Then we have $|V_{2}| \geqslant l$ and $|V_{1}| \leqslant n-2l$.
	
	From Observation~\ref{obs:3reg} it follows that $2|V_1|\geqslant |V_{-1}|$, which implies $n-2l\geqslant |V_1|\geqslant \frac12 l$, or equivalently $l \leqslant \frac25 n$. Using Observation~\ref{obs:x}, we deduce
	
	$$f(V)\geqslant n-l \geqslant n-\frac2 5 n = \frac3 5 n.$$
	Since  PLSRDN must be an integer, we conclude that 	$$\gamma_{plsr}(G)\geqslant \left\lceil\frac{3}{5}n\right\rceil.$$
\end{proof}



The following theorem proposes the exact value of PLSRDN for paths and cycles.

\begin{theorem}

	Let $P_n$ be a path of the length $n$ and let $C_n$ be a cycle of the length $n$, $n\geqslant 3$. Then	
	$$	\gamma_{plsR}(P_n) = \left\lfloor\frac{2n}{3}\right\rfloor \textnormal{\quad and\quad} \gamma_{plsR}(C_n) = \left\lceil\frac{2n}{3}\right\rceil.$$
	
\end{theorem}

\begin{proof}
	
	\textit{Upper bound}.	
	Observe that the path $P_3$ can be labeled with $-1\ 2\ 1$. It is easy to see that, by joining $k\geqslant 1$ such paths, the valid solution for the path $P_{3k}$ can be obtained and that total weight of such solution equals $2k$. Let us construct the proper PLSRD functions $f$ and $g$ for a path $P_n$ and a cycle $C_n$, respectively.
	
	Case 1.  $n = 3k$. Both functions  $f$ and $g$ are constructed in the same way, by joining $k$ blocks labeled with $-1\ 2\ 1$, as described above. It holds that $f(V(P_n)) =g(V(C_n))= 2k = \frac {2n}{3}=\left\lfloor\frac{2n}{3}\right\rfloor =\left\lceil\frac{2n}{3}\right\rceil$. 
	
	Case 2.  $n = 3k+1$. In this case,  the functions $f$ and $g$ are constructed by joining $(k-1)$ blocks $-1\ 2\ 1$. 
	For the function $f$, the remaining four vertices are labeled with the pattern  $-1\ 2\ 2\ -1$. For the function $g$, the remaining four vertices are labeled with the pattern  $-1\ 2\ 1\ 1$.
	
	Then, $f(V(P_n))  = 2(k-1)+2 = \frac {2n-2}{3} = \left\lfloor\frac{2n}{3}\right\rfloor$ and  $g(V(C_n)) = 2(k-1)+3 = \frac {2n+1}{3} = \left\lceil\frac{2n}{3}\right\rceil$.
	
	Case 3. $n = 3k+2$. The functions $f$ and $g$ are constructed by joining $k$ blocks $-1\ 2\ 1$, and labeling the rest of 2 vertices with patterns $2\ -1$ and $1\ 1$, respectively. Then $f(V(P_n)) = 2k+1 = \frac {2n-1}{3} = \left\lfloor\frac{2n}{3}\right\rfloor$ and $g(V(C_n)) = 2k+2 = \frac {2n+2}{3} = \left\lceil\frac{2n}{3}\right\rceil$.
	
	\textit{Lower bound}
	
	Let $f$ be an arbitrary PLSRD function for the path $P_n$ and let $n = 3k+r$, where $r\in \{0,1,2\}$.
	
	Any path of the length $n$ can be divided into $k$ blocks of length 3, plus the remaining $r$ vertices. We claim that each of the $k$ blocks can contain at most one vertex labeled with -1. Indeed, if a block of three vertices, say $u-v-w$, contains two vertices labeled with -1,  then  $\sum_{u\in N[v]} f(u)\leqslant 0$ holds, which is not allowed by condition (C3) from Definition~\ref{def:1}.
	
	Case 1. $n = 3k$. Then $|V_{-1}|\leqslant k = \frac n 3$. From Observation~\ref{obs:x}, $f(V)\geqslant n-\frac n 3 = \frac{2n}3=\left\lfloor\frac{2n}{3}\right\rfloor$.
	
	Case 2. $n = 3k+1$. Then $|V_{-1}|\leqslant k+1 = \frac {n+2} 3$. From Observation~\ref{obs:x}, $f(V)\geqslant n-\frac {n+2} 3 = \frac{2n-2}3 = \left\lfloor\frac{2n}{3}\right\rfloor$.
	
	Case 3. $n = 3k+2$. Then $|V_{-1}|\leqslant k+1 = \frac {n+1} 3$. From Observation~\ref{obs:x}, $f(V)\geqslant n-\frac {n+1} 3 = \frac{2n-1}3 = \left\lfloor\frac{2n}{3}\right\rfloor$.
	
	Let now $g$ be an arbitrary PLSRD function for the cycle $C_n$. It can be observed that each vertex $v\in V_{-1}$ generates a subpath of a form $u-v-w$, where $u\in V_1$ and  $w\in V_2$, such that $N(u)\cap V_{-1} = \{v\}$ and   $N(w)\cap V_{-1} = \{v\}$. This means that the observed subpaths, generated by different vertices from $V_{-1}$, are disjoint.  Let $n=3k+r$, $r\in\{0,1,2\}$ Suppose that $|V_{-1}|\geqslant k+1$. Then the total number of vertices in the cycle is at least $3(k+1)$, i.e. $|V(C_n)|\geqslant 3(k+1)>n$, which is a contradiction. Therefore $|V_{-1}|\leqslant k$ and from Observation~\ref{obs:x}, $g(V(C_n))\geqslant n-k=2k+r = \left\lceil\frac{2n}{3}\right\rceil$.
	
\end{proof}

Let $L_n$ be a ladder graph with $V(L_n) = \{a_1,a_2,\ldots,a_n,b_1,b_2,\ldots,b_n\}$ and $E(L_n)=\{a_1a_2,\ldots, a_{n-1}a_n, b_1b_2,\ldots, b_{n-1}b_n, a_1b_1,\ldots ,a_nb_n \}$. The ladder graph is formed as the Cartesian product of path graphs $P_n$ and  $P_2$.

\begin{theorem} \label{thm:ladder} 
	For a ladder graph $L_n$, it holds
	
	$$
	\gamma_{plsR}(L_n) = \left\lceil\frac{6n}{5}\right\rceil$$
\end{theorem}

\begin{proof}
	\textit{Upper bound}
	
	Let us consider the following two labelings of the graph $L_5$.
	
	$X: \begin{tikzcd}[column sep=tiny]
		-1 & 1 & -1 & 2 & 2 \\[-20pt]
		2 & 2 & -1 & 1& -1 
	\end{tikzcd}$ and $Y: \begin{tikzcd}[column sep=tiny]
		2 &2 & -1 & 1 & -1 \\[-20pt]
		-1 & 1 & -1 & 2& 2 
	\end{tikzcd}$
	
	It can be noticed that both labelings satisfies the conditions of a PLSRD function of graph $L_5$, with the total weight equals 6.
	We will use the alternating join of the blocks of $X$ and $Y$ to construct a proper labeling of ladders $L_n$ with $n >5$.

	By adding block $Y$ to the right of block $X$, a valid PLSRD function for the ladder graph of dimension $2 \times 10$ is obtained.
	
	$$X|Y: \begin{tikzcd}[column sep=1pt]
		-1 & 1 & -1 & 2 & 2&\mid& 2 &2 & -1 & 1 & -1 \\[-20pt]
		2 & 2 & -1 & 1& -1&\mid &-1 & 1 & -1 & 2& 2 
	\end{tikzcd}
	$$
	
	
	The only vertices that need to be checked  are the vertices of the last column of block $X$and the first column of block $Y$,  which can be easily verified.
	
	Similarly, the same  holds if the block $X$ is added right to the block $Y$.

	\begin{itemize}[label=\textbullet]
		\item 	Case 1. Let $n \equiv 0 \pmod{5}$, i.e. $n = 5k$
		
		For even $k$, we define the function $f$ by the labeling obtained by joining of $\frac k 2$ blocks $XY$. From the consideration above, such labeling yields to the proper $PLSRD$ function, with the value equals to $6k$.
		
		For odd $k$, the proper labeling is obtained by joining of  $\frac {k-1} 2$ blocks $YX$, with the additional block $Y$ at the end. This labeling also yields to the proper $PLSRD$ function, with the same value  $6k$.
		
		For both cases we get $f(V(L_n)) = 6k =  \left\lceil\frac{6n}{5}\right\rceil$.

		Notice that regardless of the parity of $k$, the labeling ends with block  $Y$, i.e. 
		
		$\dots X|Y|X|Y\ldots X|Y: \begin{tikzcd}[column sep=0pt]
			\ldots&-1 & 1 & -1 & 2 & 2&\mid& 2 &2 & -1 & 1 & -1 \\[-20pt]
			\ldots&2 & 2 & -1 & 1& -1 &\mid&-1 & 1 & -1 & 2& 2 
		\end{tikzcd}
		$

		We will also use this construction of the function $f$ for the rest of the cases, by adding the proper number of additional columns, as follows.

		\item	Case 2. Let $n \equiv 1 \pmod{5}$, i.e. $n = 5k + 1$.

		We construct the function $f$ by joining the column $\begin{tikzcd}[column sep=0pt]
			1 \\[-20pt]
			1 
		\end{tikzcd}
		$ 	 to the solution for the ladder $L_{5k} $. 
		
		We get $f(V(L_{n})) = 6k + 2 =  \left\lceil\frac{6n}{5}\right\rceil$.
		
		\item	Case 3. Let $n \equiv 2 \pmod{5}$, i.e. $n = 5k + 2$.
		
		To the solution for the 	ladder $L_{5k} $ the block   $\begin{tikzcd}[column sep=0pt]
			1&2 \\[-20pt]
			1&-1 
		\end{tikzcd}
		$  is added and
		
		$f(V(L_{n})) = 6k + 3 =  \left\lceil\frac{6n}{5}\right\rceil$. 

		\item	Case 4. Let $n \equiv 3 \pmod{5}$, i.e. $n = 5k + 3$.
		
		To the solution for the 	ladder $L_{5k} $ the block    $\begin{tikzcd}[column sep=0pt]
			1&-1&2 \\[-20pt]
			2&-1 &1
		\end{tikzcd}
		$ 
		is added and
		
		$f(V(L_{n})) = 6k + 4 =  \left\lceil\frac{6n}{5}\right\rceil$. 

		\item	Case 5. Let $n \equiv 4 \pmod{5}$, i.e. $n = 5k + 4$.
		
		To the solution for the 	ladder $L_{5k} $ the block
		$\begin{tikzcd}[column sep=0pt]
			1&-1&2&2 \\[-20pt]
			2&-1 &1&-1
		\end{tikzcd}
		$ 
		is added and
		
		$f(V(L_{n})) = 6k + 5 = \left\lceil\frac{6n}{5}\right\rceil$. 
	\end{itemize}
	In all cases 2--5, the valid $PLSRD$ function is obtained, which can be easily checked.
	
	\textit{Lower bound}
	
	Let $f:V(L_n)\rightarrow \{-1,1,2\}$ be an arbitrary $\mathrm{PLSRD}$ function on the ladder graph $L_n$.
	
	To prove inequality 	$f(V(L_n))\geqslant\left\lceil\frac{6n}{5}\right\rceil$, we employ a discharging method, the well known combinatorial technique originally introduced in the context of graph coloring problems \cite{wernicke1904kartographischen,cranston2017introduction}. This method has also been adapted in numerous studies on domination problems and weight distributions in graphs (see, for example, \cite{zec2024double} and references therein). For the purpose of this theorem, we design a specific column-based discharging procedure ($\mathrm{CBDP}$) on the labeling of the ladder graph $L_n$ with respect to the function 
	$f$. The procedure ensures that, after redistribution, the total weight of each column in the graph is at least $\frac{6}{5}$. To ensure clarity and effectiveness, the procedure is designed to meet the following conditions:

	$(\mathrm{CBDP_1})$ Every column of the ladder graph $L_n$ with a total weight $\leqslant 1$ must be charged
	by its neighboring column(s) in order to ensure that the total weight of this column becomes equal to $\frac{6}{5}$.
	
	$(\mathrm{CBDP_2})$ Every column of the ladder graph $L_n$ with a total weight $\geqslant 2$ can transfer some
	of its weight to the neighboring column(s), but the remaining weight of this column must be greater than or equal to
	$\frac{6}{5}$.
	
	Now, we describe a column-based discharging procedure satisfying the conditions $(\mathrm{CBDP_1})$ i $(\mathrm{CBDP_2})$: 
	
	\begin{itemize}[label=\textbullet]
		\item A column with a total weight of $-2$ must receive $\frac{16}{5}$ of charge from its neighboring column(s).
		
		\item A column with a total weight of $0$ must receive $\frac{6}{5}$ of charge from its neighboring column(s).
		
		\item A column with a total weight of $1$ must receive $\frac{1}{5}$ of charge from its neighboring column(s).
		
		\item A column with a total weight of $2$ can transfer no more than $\frac{4}{5}$ of its charge in total to the neighboring columns. 
		
		\item A column with a total weight of $3$ can transfer no more than $\frac{9}{5}$ of its charge in total to the neighboring columns. 
		
		\item A column with a total weight of $4$ can transfer no more than $\frac{14}{5}$ of its charge in total to the neighboring columns. 
	\end{itemize}
	Since every column with weight $\leqslant 1$ must be charged exclusively by importing weight from its adjacent column(s), we  define CBDP locally for each  column that requires additional charge. This localization involves redistributing weights among the neighboring column(s) of the observed column, and begins  with the immediate neighbors and progressively extending to their neighbors, along with their weight redistribution. The process continues until a column is reached on both sides of the starting column  that holds sufficient weight to potentially transfer it and, if necessary, fully charge the next column in the sequence.
	This approach guarantees that, beyond this stage, the described $\mathrm{CBDP}$ no longer require keeping track of borrowed weights, thereby allowing it to be reset for the remaining portion of the graph.

	Although the structure of a ladder graph is relatively simple, the charge redistribution process across columns involves various cases, each addressing different configurations and behaviors of the columns. To maintain clarity and the flow of the main discussion, these cases are detailed in Appendix~\ref{app:A}. The appendix covers all scenarios where a column in the ladder graph needs to be charged, categorized into boundary columns (type (a)) and interior columns (type (b)).
	Once all relevant cases for each category have been analyzed, we can conclude that the CBDP procedure leaves a charge of at least $\frac{6}{5}$ on each column. This leads to the conclusion that $f(V(L_n)) \geqslant \left\lceil \frac{6n}{5} \right\rceil$, thereby completing the proof.

\end{proof}

\textit{Prism graphs} of \textit{circular ladder graphs}, denoted as $\Pi_n$ is a graph corresponding to the skeleton of an $n$-prism. A prism graph with $2n$ vertices is a 3-regular graph, which can be visualized by taking a ladder graph and adding edges that connect corresponding boundary vertices, thus forming a prism-like structure. Additionally, it is well known that a prism graph $\Pi_n$ is isomorphic to a Generalized Petersen graph of the form $GP(n,1)$.

The following theorem establishes the exact bound for the PLSRDN of a prism graph.

\begin{theorem}
	For $n\geqslant 3$, the PLSRDN of the prism graph  $\Pi_n$  is given by: 
	$$
	\gamma_{plsR}(\Pi_n) =
	\begin{cases}
		\left\lceil\frac{6n}{5}\right\rceil+1, & \text{ if } n\equiv 4(\bmod 10) \text{ or } n\equiv 5(\bmod 10)\\
		\left\lceil\frac{6n}{5}\right\rceil, & \text{ otherwise}.
	\end{cases}
	$$
\end{theorem}

\begin{proof}
	
	\textit{Upper bound}
	
	To define a proper labeling of the prism graph $\Pi_n$, we adopt blocks used in the proof of the upper bound of the ladder graph $L_n$. However, since each column in the prism graph $\Pi_n$ is adjacent to two other columns, we need to carefully arrange these blocks to ensure they yield  a feasible PLSRD function for $\Pi_n$. This means resolving potential feasibility problems related to the adjacency of the ``first'' and ``last'' columns in the patterns for the ladder graph $L_n$. Specifically, we will concatenate an appropriate number of copies of the block $X|Y: \begin{tikzcd}[column sep=1pt]
		-1 & 1 & -1 & 2  & 2 & \mid & 2 & 2 & -1 & 1 & -1\\[-20pt]
		2 & 2 & -1 & 1 & -1 & \mid &  -1 & 1& -1 & 2 & 2
	\end{tikzcd}$, which was introduced earlier in the proof of Theorem \ref{thm:ladder}. Additionally, we append the necessary columns to the right and ensure that, after applying the column-enclosing process, the feasibility of the given labeling remains preserved.
	More precisely, the labeling that define function $f:V(\Pi_n)\rightarrow\{-1,1,2\}$ is constructed as follows:
	
	Case 1.  $n\equiv 0(\bmod 10)$;
	
	We concatenate $\frac{n}{10}$ copies of block $X|Y$ to
	form pattern $\underbrace{\underline{X|Y}\;\underline{X|Y}\dotsc \underline{X|Y}}_{\frac{n}{10}}$. 
	By interconnecting the boundary columns of this pattern, the pattern at the merging point assumes the form 
	$\begin{tikzcd}[column sep=-2pt]
		\dotsc & 1&-1& \mid&-1 & 1 & \dotsc\\[-20pt]
		\dotsc &2& 2& \mid&2 & 2 & \dotsc
	\end{tikzcd},$ 		 
	which can be easily verified as part of a valid labeling. In this case, $f(V(\Pi_n))=12\cdot\frac{n}{10}=\frac{6n}{5}=\lceil\frac{6n}{5}\rceil$.
	
	Case 2.  $n\equiv 1(\bmod 10)$;
	
	Using the block $Y_1:\begin{tikzcd}[column sep=1pt]
		2 & 2 & -1 & 1 & 1 &-1\\[-20pt]
		-1 & 1& -1 & 2 & 1 & 2
	\end{tikzcd}$, we form the pattern 
	
	\noindent$\underbrace{\underline{X|Y}\;\underline{X|Y}\dotsc \underline{X|Y}}_{\lfloor\frac{n}{10}\rfloor-1}X|Y_1$. The merging points are illustrated in the following scheme.

	$X\mid Y_1\mid X:\ \begin{tikzcd}[column sep=-2pt]
		\dots &2&2 & \mid&	2 & 2 & -1 & 1 & 1 &-1& \mid&-1&1&\dots\\[-20pt]
		\dots &1&-1& \mid&	-1 & 1& -1 & 2 & 1 & 2 &\mid&2&2&\dots
	\end{tikzcd}$

	It is straightforward to verify that  $f$ is a PLSRD function. Furthermore,
	
	$f(V(\Pi_n))=12\left(\lfloor\frac{n}{10}\rfloor-1\right)+6+8=12\lfloor\frac{n}{10}\rfloor+2=\lceil\frac{6n}{5}\rceil$.
	
	Case 3.  $n\equiv 2(\bmod 10)$;
	
	Using the block $Y_2:\begin{tikzcd}[column sep=1pt]
		2 & 2 & -1 & 1 & 2 & -1 & 1\\[-20pt]
		-1 & 1& -1 & 2 & 1 & -1 & 2
	\end{tikzcd}$, we form  pattern 							
	$\underbrace{\underline{X|Y}\;\underline{X|Y}\dotsc \underline{X|Y}}_{\lfloor\frac{n}{10}\rfloor-1}X|Y_2$. 
	
	At the merging points, this pattern assumes the form
	
	$X\mid Y_2\mid X:\ \begin{tikzcd}[column sep=-2pt]
		\dots &2&2 & \mid&	2 & 2 & -1 & 1 & 2 &-1& 1& \mid&-1&1&\dots\\[-20pt]
		\dots &1&-1& \mid&	-1 & 1& -1 & 2 & 1 & -1&2 &\mid&2&2&\dots
	\end{tikzcd}$,
	
	which can be easily verified as part of a valid labeling. In this case, $f(V(\Pi_n))=12\left(\lfloor\frac{n}{10}\rfloor-1\right)+6+9=12\lfloor\frac{n}{10}\rfloor+3=\lceil\frac{6n}{5}\rceil$.
	
	Case 4.  $n\equiv 3(\bmod 10)$;
	
	For $n=3$, we use the pattern $\begin{tikzcd}[column sep=1pt]
		1 & -1& 2\\[-20pt]
		2 & -1& 1 
	\end{tikzcd}$. We get $f(V(\Pi_3))=4=\lceil\frac{6\cdot 3}{5}\rceil$.
	
	For $n>3$, using the block $Y_3:\begin{tikzcd}[column sep=1pt]
		2 & 2 & -1 & 1 & 2 &-1 & 1 &-1\\[-20pt]
		-1 & 1& -1 & 2 & 1 & -1 &2 &2
	\end{tikzcd}$, we form pattern $\underbrace{\underline{X|Y}\;\underline{X|Y}\dotsc \underline{X|Y}}_{\lfloor\frac{n}{10}\rfloor-1}X|Y_3$. 
	
	At the merging point, this pattern takes the following form and we get the proper labeling
	
	$X\mid Y_3\mid X:\ \begin{tikzcd}[column sep=-2pt]
		\dots &2&2 & \mid&	2 & 2 & -1 & 1 & 2 &-1& 1&-1& \mid&-1&1&\dots\\[-20pt]
		\dots &1&-1& \mid&	-1 & 1& -1 & 2 & 1 & -1&2 &2&\mid&2&2&\dots
	\end{tikzcd}$.
	
	$f(V(\Pi_n))=12\left(\lfloor\frac{n}{10}\rfloor-1\right)+6+10=12\lfloor\frac{n}{10}\rfloor+4=\lceil\frac{6n}{5}\rceil$.
	
	Case 5.  $n\equiv 4(\bmod 10)$;
	
	Using the block $Y_4:\begin{tikzcd}[column sep=1pt]
		1 & -1 & 2 & 1 \\[-20pt]
		2 & -1& 1 & 1 
	\end{tikzcd}$, we form pattern $\underbrace{\underline{X|Y}\;\underline{X|Y}\dotsc \underline{X|Y}}_{\lfloor\frac{n}{10}\rfloor}Y_4$. At the merging point, this pattern assumes the form
	
	$Y\mid Y_4\mid X:\ \begin{tikzcd}[column sep=-2pt]
		\dots &1&-1 & \mid&	1 & -1 & 2 & 1 \mid&-1&1&\dots\\[-20pt]
		\dots &2&2& \mid&	2 & -1& 1 & 1\mid&2&2&\dots
	\end{tikzcd}$,

	which can be easily verified as part of a valid labeling.  In this case, $f(V(\Pi_n))=12\lfloor\frac{n}{10}\rfloor+6=\lceil\frac{6n}{5}\rceil+1$.
	
	Case 6.  $n\equiv 5(\bmod 10)$;
	
	Using the block $Y_5:\begin{tikzcd}[column sep=1pt]
		-1 & 1 & -1 & 2 & 1 \\[-20pt]
		2 & 2 & -1& 1 & 1 
	\end{tikzcd}$, we form pattern $\underbrace{\underline{X|Y}\;\underline{X|Y}\dotsc \underline{X|Y}}_{\lfloor\frac{n}{10}\rfloor}Y_5$. At the merging point, this pattern takes the following form
	
	$Y\mid Y_5\mid X:\ \begin{tikzcd}[column sep=-2pt]
		\dots &1&-1 & \mid&	-1 & 1&-1& 2 & 1 \mid&-1&1&\dots\\[-20pt]
		\dots &2&2& \mid&	2 & 2&-1& 1 & 1\mid&2&2&\dots
	\end{tikzcd}$
	
	In this case, $f(V(\Pi_n))=12\lfloor\frac{n}{10}\rfloor+7=\lceil\frac{6n}{5}\rceil+1$.
	
	Case 7.  $n\equiv 6(\bmod 10)$;
	
	Using the column $Y_6:\begin{tikzcd}[column sep=1pt]
		1 \\[-20pt]
		1 
	\end{tikzcd}$, we form pattern $\underbrace{\underline{X|Y}\;\underline{X|Y}\dotsc \underline{X|Y}}_{\lfloor\frac{n}{10}\rfloor}XY_6$. At the merging point,  we get the proper labeling
	
	$X\mid Y_6\mid X:\ \begin{tikzcd}[column sep=-2pt]
		\dots &2&2 & \mid&	1  \mid&-1&1&\dots\\[-20pt]
		\dots &1&-1& \mid&	1\mid&2&2&\dots
	\end{tikzcd}$.
	
	In this case, $f(V(\Pi_n))=12\lfloor\frac{n}{10}\rfloor+6+2=\lceil\frac{6n}{5}\rceil$.
	
	Case 8.  $n\equiv 7(\bmod 10)$;
	
	Using the block $Y_7:\begin{tikzcd}[column sep=1pt]
		2 & 1 \\[-20pt]
		-1 & 1 
	\end{tikzcd}$, we form pattern $\underbrace{\underline{X|Y}\;\underline{X|Y}\dotsc \underline{X|Y}}_{\lfloor\frac{n}{10}\rfloor}XY_7$. At the merging point, this pattern takes the  form
	
	$X\mid Y_7\mid X:\ \begin{tikzcd}[column sep=-2pt]
		\dots &2&2 & \mid&2&	1  \mid&-1&1&\dots\\[-20pt]
		\dots &1&-1& \mid&-1&	1\mid&2&2&\dots
	\end{tikzcd}$.
	
	In this case, $f(V(\Pi_n))=12\lfloor\frac{n}{10}\rfloor+6+3=\lceil\frac{6n}{5}\rceil$.
	
	Case 9.  $n\equiv 8(\bmod 10)$;
	
	Using the block $Y_8:\begin{tikzcd}[column sep=1pt]
		2 & -1 & 1 \\[-20pt]
		1 & -1 & 2 
	\end{tikzcd}$, we form pattern $\underbrace{\underline{X|Y}\;\underline{X|Y}\dotsc \underline{X|Y}}_{\lfloor\frac{n}{10}\rfloor}XY_8$. At the merging point, this pattern assumes the form
	
	$X\mid Y_8\mid X:\ \begin{tikzcd}[column sep=-2pt]
		\dots &2&2 & \mid&2&	-1 &1 \mid&-1&1&\dots\\[-20pt]
		\dots &1&-1& \mid&1&	-1&2\mid&2&2&\dots
	\end{tikzcd}$.
	
	In this case, $f(V(\Pi_n))=12\lfloor\frac{n}{10}\rfloor+6+4=\lceil\frac{6n}{5}\rceil$.
	
	Case 10.  $n\equiv 9(\bmod 10)$;
	
	Using the block $Y_9:\begin{tikzcd}[column sep=1pt]
		2 & -1 & 1 & -1\\[-20pt]
		1 & -1 & 2 & 2
	\end{tikzcd}$, we form the pattern $\underbrace{\underline{X|Y}\;\underline{X|Y}\dotsc \underline{X|Y}}_{\lfloor\frac{n}{10}\rfloor}XY_9$. At the merging points, this pattern takes the following form:
	
	$X\mid Y_9\mid X:\ \begin{tikzcd}[column sep=-2pt]
		\dots &2&2 & \mid&2&	-1 &1 &-1\mid&-1&1&\dots\\[-20pt]
		\dots &1&-1& \mid&1&	-1&2&2\mid&2&2&\dots
	\end{tikzcd}$
	
	In this case, $f(V(\Pi_n))=12\lfloor\frac{n}{10}\rfloor+6+5=\lceil\frac{6n}{5}\rceil$.
	
	Considering all the cases examined, we conclude that $f$ is a valid PLSRD function, with 
	$$
	f(V(\Pi_n)) =
	\begin{cases}
		\left\lceil\frac{6n}{5}\right\rceil+1, & \text{ if } n\equiv 4(\bmod 10) \text{ or } n\equiv 5(\bmod 10)\\
		\left\lceil\frac{6n}{5}\right\rceil, & \text{ otherwise}.
	\end{cases}
	$$
	Therefore, 
	$$
	\gamma_{plsR}(\Pi_n) \leqslant
	\begin{cases}
		\left\lceil\frac{6n}{5}\right\rceil+1, & \text{ if } n\equiv 4(\bmod 10) \text{ or } n\equiv 5(\bmod 10)\\
		\left\lceil\frac{6n}{5}\right\rceil, & \text{ otherwise}.
	\end{cases}
	$$
	
	\textit{Lower bound}.
	Since the total number of vertices of $\Pi_n$ is $2n$, it follows directly from  Theorem~\ref{thm:3reg} that $\gamma_{plsR}(\Pi_n)\geqslant \left\lceil\frac{6n}{5}\right\rceil$. To conclude the proof,  it remains to demonstrate that  $\gamma_{plsR}(\Pi_n)>\left\lceil\frac{6n}{5}\right\rceil$ holds for $n\equiv 4(\bmod 10) \text{ or } n\equiv 5(\bmod 10)$. Let us analyze these two cases.
	
	Case 1. $n=10m+4,m\geqslant 0$;
	
	Let us assume, by the sake of contradiction, that $\gamma_{plsR}(\Pi_n)=\left\lceil\frac{6n}{5}\right\rceil$. Thus, there exists a labeling via a PLSRD function $f$ such that $f(V(\Pi_n))=\left\lceil\frac{6n}{5}\right\rceil =\left\lceil\frac{60m+24}{5}\right\rceil = 12m+5$. From Observation~\ref{obs:x} it follows $|V_{-1}|\geqslant 2n-\gamma_{plsR}(\Pi_n)=20m+8-12m-5 = 8m+3$. We will show that this assumption implies that $|V_{-1}| = |V_{2}| = 8m+3$.
	Suppose that $|V_{-1}|\geqslant 8m+4$. Then $|V_{2}|\geqslant 8m+4$ and $|V_{1}|\leqslant 20m+8-16m-8 = 4m<\frac12|V_{-1}|$, which is not possible because of Observation~\ref{obs:3reg}. Therefore  $|V_{-1}|=8m+3$. Further, if $|V_2|\geqslant 8m+4$, we  have  $|V_{1}|\leqslant 20m+8-8m-3-8m-4 = 4m+1<\frac12|V_{-1}|$, which is again in a contradiction with Observation~\ref{obs:3reg}. Therefore,  $|V_2|= 8m+3$ and $|V_1| = 4m+2$.
	
	In other words, $8m+3$ vertices from
	$V_{-1}$ are covered by $4m+2$ vertices from $V_1$. Due to the structure of the prism graph and the constraints of the function $f$, there are only two possible configuration of this type:
	
	$(i)$ One vertex from $V_{-1}$ is exclusively covered by a single vertex from $V_{1}$, while the  remaining $8m+2$ vertices from $V_{-1}$ are covered by $4m+1$ vertices from $V_1$, arranged into $4m+1$ disjoint chains of the form $-1\, 1 -1$. The considered chain $-1\, 1$ can  either be aligned in the same column or in the same row (only non-symmetrical configurations are considered), as illustrated below. On the following schemes the vertices from the considered chain are bolded.
	
	$(i_1)$ Chain $-1\, 1$ is a column: 	
	$\begin{tikzcd}[column sep=tiny]
		\textbf{-1} & 2 & 1 & \boxed{-1}\\[-20pt]
		\textbf{1} & 2 & \boxed{-1} & 1
	\end{tikzcd}$
	
	$(i_2)$  Chain $-1\, 1$ is within a row. We have three possibilities.
	
	$(i_{21})$ 	$\begin{tikzcd}[column sep=tiny]
		\textbf{-1} & \textbf{1} & 1 & \boxed{-1}\\[-20pt]
		2 & 2 & \boxed{-1} & 1
	\end{tikzcd}$
	
	$(i_{22})$ 	$\begin{tikzcd}[column sep=tiny]
		2 &  \textbf{-1} &  \textbf{1} & 2 \\[-15pt]
		2 & -1 & \boxed{1} & -1
	\end{tikzcd}$
	
	$(i_{23})$
	$\begin{tikzcd}[column sep=tiny]
		2 & \textbf{-1} & \textbf{1} & 2 & -1 & \boxed{1} & -1 \\[-20pt]
		1 & -1 & 2 & 2 & -1 & \boxed{1} & -1 
	\end{tikzcd}$
	
	However, neither of the given configurations is valid. In $(i_1)$ and $(i_{21})$, the single boxed labels belong to two non-disjoint chains of the form $-1\,\, 1 -1$, while the single boxed labels in $(i_{22})$ and $(i_{23})$ fail to satisfy the sum condition.
	
	$(ii)$ There is one chain of the form $-1\,\, 1 -1\,\, 1 -1$, while the  remaining $8m$ vertices from $V_{-1}$ are covered by $4m$ vertices from $V_1$, arranged into $4m$ disjoint chains of the form $-1\,\, 1 -1$;
	
	The chain $-1\,\, 1 -1\,\, 1 -1$ must be placed in one row, because any other placement results to a non-feasible labeling. As a result, $m\geqslant 1$ holds in this case. The portion of the graph that involves this chain is
	uniquely depicted by the following block (the considered chain is bolded):
	
	$$B:\begin{tikzcd}[row sep=tiny,column sep=tiny]
		-1&	1 & -1 & 2 & 2 & 2 & -1 & 1 & -1 & \\
		2& 2 &  \textbf{-1} & \textbf{1} & \textbf{-1}& \textbf{1}&  \textbf{-1} & 2 & 2 &
	\end{tikzcd}$$
	
	Moreover, using the fact that in the remaining graph each 1 covers exactly two -1s, the only way block $B$ can continue is if the left side takes the form $\dotsc YXY XY B$, while the right side takes the form $B\, XYXYX \dotsc$ 	
	However, considering that there are $10m-5=5(2m-1)$ columns outside of block $B$, at the point of merging, one of the blocks $X$ or $Y$ must be connected to itself, resulting in a non-feasible labeling. 
	
	The obtained contradiction indicates that, in the case $n\equiv 4(\bmod 10)$, condition $\gamma_{plsR}(\Pi_n)\geqslant\left\lceil\frac{6n}{5}\right\rceil+1$ is satisfied.

	Case 2. $n=10m+5,m\geqslant 0$;
	
	Let us assume, for the sake of contradiction, that $\gamma_{plsR}(\Pi_n)=\left\lceil\frac{6n}{5}\right\rceil$. Thus, there exists a labeling via a PLSRD function $f$ such that $f(V(\Pi_n))=\left\lceil\frac{6n}{5}\right\rceil$. Similarly as in Case 1, it can be shown that this implies  $|V_{-1}| =  |V_{2}|=8m+4$ and $|V_{1}| = 4m+2$. In other words, $8m+4$ vertices from
	$V_{-1}$ are covered by $4m+2$ vertices from $V_1$. The only possible configuration satisfying this fact is that $8m+4$ vertices from $V_{-1}$ are covered by $4m+2$ vertices from $V_1$, arranged into $4m+2$ disjoint chains of the form $-1\,\, 1 -1$.
	These chains can either be $L-$ shaped or aligned within the same row. If at least one of these chains is $L$-shaped,  the following configuration arises (symmetries are ignored, and vertices from a portion of the chain are bolded):  
	$\begin{tikzcd}[column sep=tiny, ]
		\boxed{\textbf{-1}} & 1 \\[-20pt]
		\textbf{1} & \boxed{\textbf{-1}} 
	\end{tikzcd}$.
	
	However, this configuration is not valid, since the boxed labels belong to two non-disjoint chains of the form $-1\,\, 1 -1$. Therefore, we can assume that all $-1\,\, 1 -1$ chains are arranged in rows, as illustrated below (one such chain is bolded).
	
	$\begin{tikzcd}[column sep=tiny]
		2&2&\textbf{-1} & \textbf{1} & \textbf{-1} & -1 & 1 & -1 & 2 & 2\\[-10pt]
		-1&1&-1 & 2 & 2 & 2 & 2 & -1 & 1 & -1
	\end{tikzcd}$.
	Essentially, the given configuration is formed by alternately joining block $X$ with block $Y$. However, considering that there are $10m+5=5(2m+1)$ columns, at the point of merging, one of the blocks $X$ or $Y$ must be connected to itself, resulting in a non-feasible labeling. 
	
	The obtained contradiction indicates that, in the case $n\equiv 5(\bmod 10)$, condition $\gamma_{plsR}(\Pi_n)\geqslant\left\lceil\frac{6n}{5}\right\rceil+1$ is satisfied.
	
	This completes the proof.
	
\end{proof}

In the next theorem, we provide the exact value of PLSRDN for $Grid_{3\times n}$. 
A grid graph $Grid_{3\times n}$ is a graph consisting of $n$ rows and 3 columns, where the vertices are arranged in a rectangular scheme. Two vertices are adjacent if they occupy neighboring cells in the scheme, either horizontally or vertically. 
It should be noted  that the upper bound for a general grid graph $Grid_{m\times n}$ was established in \cite{akwu2026perfect} and is given by  $\frac{3}{4}mn$. However, for the specific case where $m$ is divisible by 3, the authors derived a construction that provides a tighter upper bound, equals $\frac{2}{3}mn$. More precisely, when $m=3$, the construction leads to an upper bound of  $2n$.  In the following theorem, it will be shown that $2n$ also serves as the lower bound for $Grid_{3\times n}$.

\begin{theorem}\label{thm:grid3} Let $Grid_{3,n}$ be a grid graph, where $n\geqslant 1$. Then 
	$$\gamma_{plsR}(Grid_{3\times n}) = 2n$$
\end{theorem}
\begin{proof}
	
	\textit{Upper bound}
	
	It is straightforward to observe that  the following sketch of the configuration, also used in \cite{akwu2026perfect}
	
	$$\begin{tikzcd}[column sep=tiny]
		-1 & 1 & -1 & 1 &\dotsc \\[-20pt] 
		2 & 2 & 2  &2 & \dotsc\\[-20pt] 
		1 & -1 & 1 & -1 & \dotsc\\[-20pt]
	\end{tikzcd}$$
	
	defines a proper PLSRD function, where the sum of the labels of each column equals 2, concluding that $\gamma_{plsR}(Grid_{3\times n}) \leqslant 2n$.
	
	\textit{Lower bound}

	Let $f$ be an arbitrary PLSRD function on the grid graph $Grid_{3\times n}$. The idea of the proof is to show that it is always possible to construct a family of pairwise disjoint subsets of $V_{-1}\cup V_1\cup V_2$, consisting of $|V_{-1}|$ subsets, each of which corresponds to a triplet labeling of the form 
	$\{-1,1,2\}$. Once this is verified, we observe that, since the grid graph $Grid_{3\times n}$ contains $3n$ vertices, it follows that $|V_{-1}|\leqslant n$. Consequently, we deduce that $f(V(Grid_{3\times n}))\geqslant 3n-n=2n$.
	
	From the definition of the PLSRD function $f$ and  the structure of the $Grid_{3\times n}$ graph, it follows that each vertex in $V_{-1}$ is uniquely paired with a neighboring vertex from $V_2$, and is adjacent to at least one vertex from $V_1$. Each vertex $v\in V_{-1}$ belongs to one of the following two types:
	
	Type 1. There is a vertex $u\in N(v)\cap V_1$ such that $N(u)\cap (V_{-1}\setminus\{v\})=\emptyset$;
	
	In this case, there is a one-to-one correspondence between the vertices $u$ and $v$, meaning that labeling of these vertices, along with the labeling of a vertex from  $N(v)\cap V_2$, naturally  forms a triplet $\{-1,1,2\}$ .
	
	Type 2. For each vertex $u\in N(v)\cap V_1$, it holds $N(u)\cap (V_{-1}\setminus\{v\})\ne\emptyset$;
	
	Since any vertex from $V_1$ cannot have more than two adjacent vertices from $V_{-1}$, we have that $N(u)\cap (V_{-1}\setminus\{v\})=\{v'\}$, for some vertex $v'$. The labeling of the vertices $v$, $u$ and $v'$ yields an alternating chain $-1\,\,1-1$. Several examples of such a configuration is shown in Figure~\ref{fig:labellingsGrid}.
	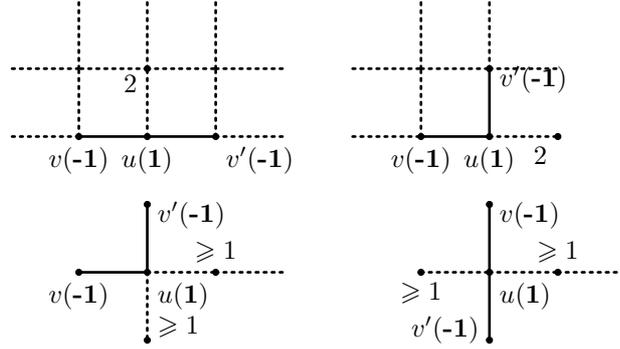
\begin{figure}
		\centering
		\begin{tikzpicture}[scale=0.45, line cap=round,line join=round,>=triangle 45,x=1cm,y=1cm]
			\draw [line width=1pt] (-11,0) node[below]{$v(\textbf{-1})$}-- (-9,0);
			\draw [line width=1pt] (-9,0) node[below]{$u(\textbf{1})$}-- (-7,0);
			\draw [line width=1pt,dotted] (-9,0) 	-- node[above left]{2}(-9,2);
			
			\draw [line width=1pt,dotted] (-7,0) node[below right]{$v'(\textbf{-1})$}-- (-5,0);
			
			\draw [line width=1pt,dotted] (-7,0)-- (-7,2);
			\draw [line width=1pt,dotted] (-11,0)-- (-13,0);
			\draw [line width=1pt,dotted] (-11,2)-- (-11,4);
			\draw [line width=1pt,dotted] (-9,2)-- (-9,4);
			\draw [line width=1pt,dotted] (-7,2)-- (-7,4);
			\draw [line width=1pt,dotted] (-11,2)-- (-13,2);
			\draw [line width=1pt,dotted] (-9,2)-- (-11,2);
			\draw [line width=1pt,dotted] (-9,2)-- (-7,2);
			\draw [line width=1pt,dotted] (-7,2)-- (-5,2);
			\draw [line width=1pt,dotted] (-11,0)-- (-11,2);
			\draw [fill=black] (-7,0) circle (2.5pt);
			\draw [fill=black] (-9,2) circle (2.5pt);
			\draw [fill=black] (-9,0) circle (2.5pt);
			\draw [fill=black] (-11,0) circle (2.5pt);
			
			\draw [line width=1pt] (-1,0) node[below]{$v(\textbf{-1})$}-- (1,0);
			\draw [line width=1pt,dotted] (-1,2) -- (1,2);
			\draw [line width=1pt,dotted] (-3,2) -- (-1,2);
			\draw [line width=1pt,dotted] (1,2) -- (3,2);
			\draw [line width=1pt,dotted] (1,0) node[below]{$u(\textbf{1})$}-- node[below right]{2}(3,0);
			\draw [line width=1pt] (1,0) 	-- node[above right]{$v'(\textbf{-1})$}(1,2);
			\draw [line width=1pt, dotted] (-3,0) -- (-1,0);
			\draw [line width=1pt, dotted] (-1,0) -- (-1,2);
			\draw [line width=1pt,dotted] (-1,2) -- (-1,4);
			\draw [line width=1pt,dotted] (1,2) -- (1,4);
			
			\draw [fill=black] (-1,0) circle (2.5pt);
			\draw [fill=black] (1,2) circle (2.5pt);
			\draw [fill=black] (1,0) circle (2.5pt);
			\draw [fill=black] (3,0) circle (2.5pt);
			
			\draw [line width=1pt] (-11,-4) node[below]{$v(\textbf{-1})$}-- (-9,-4);
			\draw [line width=1pt,dotted] (-9,-4) node[below right]{$u(\textbf{1})$}-- (-7,-4);
			\draw [line width=1pt] (-9,-4) 	-- node[above right]{$v'(\textbf{-1})$}(-9,-2);
			
			\draw [line width=1pt,dotted] (-7,-4) node[above]{$\geqslant 1$}-- (-5,-4);

			\draw [line width=1pt,dotted] (-9,-4)-- node[below right]{$\geqslant 1$}(-9,-6);
			
			\draw [fill=black] (-7,-4) circle (2.5pt);
			\draw [fill=black] (-9,-2) circle (2.5pt);
			\draw [fill=black] (-9,-4) circle (2.5pt);
			\draw [fill=black] (-11,-4) circle (2.5pt);
			\draw [fill=black] (-9,-6) circle (2.5pt);
			
			\draw [line width=1pt,dotted] (-1,-4) node[below]{$\geqslant 1$}-- (1,-4);
			\draw [line width=1pt,dotted] (1,-4) node[below right]{$u(\textbf{1})$}-- (3,-4);
			\draw [line width=1pt] (1,-4) 	--node[above right]{$v(\textbf{-1})$} (1,-2);
			
			\draw [line width=1pt,dotted] (3,-4) node[above]{$\geqslant 1$}-- (5,-4);

			\draw [line width=1pt] (1,-4)-- node[below left]{$v'(\textbf{-1})$}(1,-6);
			
			\draw [fill=black] (3,-4) circle (2.5pt);
			\draw [fill=black] (1,-2) circle (2.5pt);
			\draw [fill=black] (1,-4) circle (2.5pt);
			\draw [fill=black] (-1,-4) circle (2.5pt);
			\draw [fill=black] (1,-6) circle (2.5pt);
			
		\end{tikzpicture}
		\caption{Segments of several different configurations illustrating the mutual positioning of vertices $v$, $u$, and $v'$, labeled as -1, 1, and 1, respectively.}\label{fig:labellingsGrid}
	\end{figure}
	If vertex $v'$ is of Type 1, it has its exclusive vertex labeled with 1 as a neighbor. In that case, we can add a vertex from $V_1$ to the chain, where the vertex correspondences with $v'$, thereby forming two triplets $\{-1,1,2\}$. Otherwise, vertex $u$ is the only vertex from $V_1$ adjacent to the endpoints of the alternating chain $-1\,\,1-1$, so we need to find the vertex from $V_1$ that completes the triplet elsewhere. More generally, for the labels of neighboring vertices, we can form an alternating chain with the pattern $-1\,\,1-1\dotsc-1\,\,1-1$ (starting and ending with $-1$). WLOG, we assume that this alternating chain is of maximal length, meaning it is not part of a larger alternating chain with  $-1$s as endpoints.
	
	Again, it is possible that at least one of the endpoints of this alternating chain corresponds to the Type 1 vertex from $V_{-1}$, or that it has a neighbor in $V_1$ whose only neighbors from $V_{-1}$ are part of the chain. In this case, we add the appropriate vertex with label $1$ to the chain, ensuring  each  $-1$  is assigned to exactly one triplet. A more challenging case arises when the alternating chain $-1\,\,1-1\dotsc-1\,\,1-1$ is isolated and lacks an ``obvious'' adjacent $1$ to borrow. To determine all of the triplets containing  $-1$s from the chain, we need to identify an additional, chain-independent vertex labeled $\geqslant 1$, and incorporate it into the chain. A potential candidate for this extension can be found in the following source types in the graph:
	\begin{itemize}[label=\textbullet]
		\item Source 1: A label of a vertex  $w\in V_1\cup V_2$ with no neighbors from $V_{-1}$. 
		
		\item Source 2: A label of a vertex $w\in V_{-1}$ having at least two neighbors from $V_1$, each of which has no other neighbors from $V_{-1}$. In this case, label of any of these vertices from $V_1$ could be reassigned to the observed alternating chain.
		
		\item Source 3: A different (maximal) alternating chain, with both endpoints corresponding to vertices from $V_{-1}$ of type 1, which contains a ``spare'' label $1$ at one end.  This label can be  borrowed to extend the observed alternating chain. 
	\end{itemize}
	It should be noted that in the case of Source 1, if an isolated label 2 appears, a 'stronger' triplet of the form $(-1, 2, 2)$ will be created, which contributes an even greater value to the total sum.

	Let $\mathcal{C}$ denote the family of all alternating chains of maximal length, with endpoints labeled $-1$, such that both endpoints have no adjacent label $1$ outside of this chain. For each chain in the family $\mathcal C$, it is essential to ensure the existence of a label $1$ (or, in rare cases, a label $2$, although this is not a common scenario), selected from one of the previously specified sources, which can be added to the chain. Moreover, it is necessary that the label assignment be uniquely determined for each chain. This can be accomplished by selecting the next unassigned label $1$ from the nearest available source. In this manner, triplets corresponding to the $-1$s from the chains in $\mathcal C$ are determined; the remaining $-1$s in the graph labeling correspond one-to-one with their adjacent label $1$ and their triplets can be easily described.
	
	The precise procedure for appending the chain $C\in\mathcal C$ is highly dependent on its placement within the labeling of the grid graph $Grid_{3\times n}$. In particular, various shapes and sizes of $C$ can lead to numerous distinct cases that require consideration. 
	
	The detailed analysis of the construction for triplets associated with the alternating chains, as well as considerations for the remaining  -1's is provided in Appendix~\ref{app:B}. These cases involve identifying unique configurations and ensuring consistency across the labeling process.

	After all variants that can appear in the consideration of all chains of maximal length are exhausted, as it is shown in Appendix~\ref{app:B}, we can conclude that the whole set of vertices can be divided into triplets containing labels $-1$, $1$, and $2$, where the maximum number of such triplets is $n$. This concludes the proof of the lower bound, which is $2n$.
	
\end{proof}

The flower snarks graph $J_{2n+1}=\left(V(J_{2n+1}),E(J_{2n+1})\right),n\geqslant 2$, consists of the following sets of vertices and edges:
\begin{align*}
	V(J_{2n+1})&=\{a_i,b_i,c_i,d_i:i\in\{0,1,\dotsc,2n\}\},\\
	E(J_{2n+1})&=\{a_ia_{i+1},c_ic_{i+1},d_id_{i+1}:i\in\{0,1,\dotsc,2n-1\}\}\\
	&\cup\{a_{2n}a_0\}\cup\{c_{2n}d_0,d_{2n}c_0\}\cup\{a_ib_i,b_ic_i,b_id_i:i\in\{0,1,\dotsc,2n\}\}.
\end{align*}

Flower snarks are visually depicted as in Figures~\ref{fig:upper bound J_9}~and~\ref{fig:upper bound J_{11}}. The vertices $a_i$ are positioned in the innermost region of the graph, the vertices $b_i$ occupy the middle region, while the vertices $c_i$ and $d_i$ form two distinct layers in the outer region of the graph.

\begin{theorem} For the flower snarks graph $J_{2n+1}$, $n\geqslant 2$, it holds
	$$ \gamma_{plsR}(J_{2n+1})\leqslant
	\begin{cases}
		5n+3, &\mbox{for even number } n,\\
		5n+4  &\mbox{for odd number } n.
	\end{cases}$$
\end{theorem}
\begin{proof}

	We define function the $f:V(J_{2n+1})\rightarrow \{-1,1,2\}$, based on the parity of $n$, and therefore distinguish between the following two cases.
	
	$\bullet$ Case 1. $n\geqslant 2$ is an even number;
	
	We define 
	\begin{align*}
		&f(a_n)=2, f(a_{n-1})=f(a_{n+1})=1, f(a_{n-2})=f(a_{n+2})=-1, \\
		&f(b_n)=f(b_{n-1})=f(b_{n+1})=-1, f(b_{n-2})=f(b_{n+2})=2,\\
		&f(c_n)=1, f(c_{n-1})=-1, f(c_{n+1})=2, f(c_{n-2})=2, f(c_{n+2})=1,\\
		&f(d_n)=1, f(d_{n-1})=2, f(d_{n+1})=-1, f(d_{n-2})=1, f(d_{n+2})=2.
	\end{align*}
	Also, for $n\geqslant 4$ and $j\in\{1,2,\dotsc,n-2\}$, we define
	\begin{align*}
		&f(a_{n-2-j})=\begin{cases}
			-1, &\mbox{ for } j\equiv 0(\bmod 4) \mbox{ or } j\equiv 1(\bmod 4),\\
			2  &\mbox{ otherwise},
		\end{cases}\\
		&f(b_{n-2-j})=1,\\
		&f(c_{n-2-j})=f(d_{n-2-j})=\begin{cases}
			2, &\mbox{ for } j\equiv 0(\bmod 4) \mbox{ or } j\equiv 1(\bmod 4),\\
			-1  &\mbox{ otherwise}.
		\end{cases}
	\end{align*} 
	For $n\geqslant 4$, the remaining vertices are labeled in such a way that symmetry relative to the vertices $a_n$, $b_n$, $c_n$ and $d_n$ is preserved. More precisely, for $j\in\{1,2,\dotsc,n-2\}$, we set 
	\begin{align*}
		&f(a_{n+2+j})=f(a_{n-2-j}), f(b_{n+2+j})=f(b_{n-2-j}),\\
		&f(c_{n+2+j})=f(c_{n-2-j}), f(d_{n+2+j})=f(d_{n-2-j}).	
	\end{align*}
	
	An example of the labeling of the flower snarks graph, using the previously described pattern, is given in the Figure \ref{fig:upper bound J_9}. 
	
	\begin{figure}
		\centering
		\begin{tikzpicture}[scale=0.6, line cap=round,line join=round,>=triangle 45,x=1cm,y=1cm]
			\draw [line width=1pt] (-0.86,-6.53) node[above]{$1$} -- (0.84,-6.57);
			\draw [line width=1pt] (0.84,-6.57) node[above]{$-1$} -- (2.1679870576897233,-5.507902841257643) ;
			\draw [line width=1pt] (2.1679870576897233,-5.507902841257643) node[left]{$-1$}-- (2.502581269843993,-3.8406755882435677);
			\draw [line width=1pt] (2.502581269843993,-3.8406755882435677) node[left]{$2$}-- (1.6872222859953718,-2.3484324018100224);
			\draw [line width=1pt] (1.6872222859953718,-2.3484324018100224) node[below left]{$2$}-- (0.10342563639235468,-1.7294104533249488);
			\draw [line width=1pt] (0.10342563639235468,-1.7294104533249488) node[below]{$-1$}-- (-1.5077326246767149,-2.273256992147148);
			\draw [line width=1pt] (-1.5077326246767149,-2.273256992147148) node[below]{$-1$}-- (-2.3923736408280933,-3.725500178580693);
			\draw [line width=1pt] (-2.3923736408280933,-3.725500178580693) node[right]{$1$}-- (-2.136564048914801,-5.406619285808124);
			\draw [line width=1pt] (-2.136564048914801,-5.406619285808124) node[right]{$a_4(\textbf{2})$}-- (-0.86,-6.53);
			\draw [line width=1pt] (0.10342563639235468,-1.7294104533249488) -- (0.10017224365241753,-0.22941398151720493);
			\draw [line width=1pt] (1.6872222859953718,-2.3484324018100224)-- (2.6952636896379882,-1.2376432550579382);
			\draw [line width=1pt] (2.502581269843993,-3.8406755882435677)-- (3.9224196569869996,-3.356879190577002);
			\draw [line width=1pt] (2.1679870576897233,-5.507902841257643)-- (3.579512478375505,-6.015441989262826);
			\draw [line width=1pt] (0.84,-6.57)-- (1.6588207160472934,-7.826794587421428);
			\draw [line width=1pt] (-0.86,-6.53)-- (-1.493308308829422,-7.8897501924867);
			\draw [line width=1pt] (-2.136564048914801,-5.406619285808124) -- (-3.425541485037292,-6.17378104786818);
			\draw [line width=1pt] (-2.3923736408280933,-3.725500178580693)-- (-3.8430331937006743,-3.343941714306183);
			\draw [line width=1pt] (-1.5077326246767149,-2.273256992147148)-- (-2.603329719814657,-1.2487244132970154);
			\draw [line width=1pt] (0.10017224365241753,-0.22941398151720493) node[right]{$1$}-- (0.7199397440771068,0.5553714958059505) node[right]{$2$};
			\draw [line width=1pt] (0.10017224365241753,-0.22941398151720493)-- (-0.561876818362425,0.5200465174320507) node[left]{$2$};
			\draw [line width=1pt] (2.6952636896379882,-1.2376432550579382) node[left]{$1$}-- (3.6952636896379882,-1.2376432550579382) node[below]{$-1$};
			\draw [line width=1pt] (2.6952636896379882,-1.2376432550579382)-- (2.7928954996468978,-0.24242065198142837) node[left]{$-1$};
			\draw [line width=1pt] (3.9224196569869996,-3.356879190577002) node[below]{$1$}-- (4.908254342317568,-3.5245991154493264) node[below]{$-1$};
			\draw [line width=1pt] (3.9224196569869996,-3.356879190577002)-- (4.400779721658352,-2.478715396898613) node[left]{$-1$};
			\draw [line width=1pt] (3.579512478375505,-6.015441989262826) node[below]{$1$}-- (4.30258408754147,-6.706215069237823) node[below]{$2$};
			\draw [line width=1pt] (3.579512478375505,-6.015441989262826)-- (4.497924720264232,-5.619817225188429) node[above]{$2$};
			\draw [line width=1pt] (1.6588207160472934,-7.826794587421428) node[left]{$2$}-- (2.5927475348063744,-8.18425884745667) node[right]{$2$};
			\draw [line width=1pt] (1.6588207160472934,-7.826794587421428)-- (1.5124801628682145,-8.816028858206282) node[left]{$1$};
			\draw [line width=1pt] (-1.493308308829422,-7.8897501924867) node[right]{$-1$}-- (-1.110263518864173,-8.813479962971812) node[right]{$-1$};
			\draw [line width=1pt] (-1.493308308829422,-7.8897501924867)-- (-2.194024813288919,-8.603190015021204) node[left]{$2$};
			\draw [line width=1pt] (-3.425541485037292,-6.17378104786818) node[right]{$\ b_4(\textbf{-1})$}-- (-3.721169720956498,-7.129084112935647) node[right]{$d_4(\textbf{1})$};
			\draw [line width=1pt] (-3.425541485037292,-6.17378104786818)-- (-4.423336452189108,-6.240152755142491) node[above]{$c_4(\textbf{1})$};
			\draw [line width=1pt] (-3.8430331937006743,-3.343941714306183)-- (-4.797677572510462,-3.6416901148172083) node[below]{$2$};
			\draw [line width=1pt] (-3.8430331937006743,-3.343941714306183) node[below]{$-1$}-- (-4.496183938815856,-2.5867138066353927) node[right]{$-1$};
			\draw [line width=1pt] (-2.603329719814657,-1.2487244132970154) -- (-3.5778252714046817,-1.0243175953719486) node[below]{$1$};
			\draw [line width=1pt] (-2.603329719814657,-1.2487244132970154) node[right]{$2$}-- (-2.619073395763152,-0.24884835264368188) node[right]{$2$};
			\draw [shift={(4.605419214929519,-5.115407092343575)},line width=1pt]  plot[domain=-1.7589111188422972:1.382681534747496,variable=\t]({1*1.6193761557249466*cos(\t r)+0*1.6193761557249466*sin(\t r)},{0*1.6193761557249466*cos(\t r)+1*1.6193761557249466*sin(\t r)});
			\draw [shift={(2.907532125204842,-7.7611219637220525)},line width=1pt]  plot[domain=-2.4941496586791065:0.6474429949106866,variable=\t]({1*1.7489992949254258*cos(\t r)+0*1.7489992949254258*sin(\t r)},{0*1.7489992949254258*cos(\t r)+1*1.7489992949254258*sin(\t r)});
			\draw [shift={(-0.3407723252103523,-8.709609436613743)},line width=1pt]  plot[domain=3.084232588647188:6.225825242236981,variable=\t]({1*1.8563054376533754*cos(\t r)+0*1.8563054376533754*sin(\t r)},{0*1.8563054376533754*cos(\t r)+1*1.8563054376533754*sin(\t r)});
			\draw [shift={(-3.308680632739014,-7.421671385081847)},line width=1pt]  plot[domain=2.3270834956813458:5.468676149271139,variable=\t]({1*1.6243287440440541*cos(\t r)+0*1.6243287440440541*sin(\t r)},{0*1.6243287440440541*cos(\t r)+1*1.6243287440440541*sin(\t r)});
			\draw [shift={(-4.459760195502483,-4.413433280888942)},line width=1pt]  plot[domain=1.590733113942834:4.732325767532627,variable=\t]({1*1.8270825724893007*cos(\t r)+0*1.8270825724893007*sin(\t r)},{0*1.8270825724893007*cos(\t r)+1*1.8270825724893007*sin(\t r)});
			\draw [shift={(-3.557628667289504,-1.4177810796395374)},line width=1pt]  plot[domain=0.8942796124658572:4.03587226605565,variable=\t]({1*1.4990963004263174*cos(\t r)+0*1.4990963004263174*sin(\t r)},{0*1.4990963004263174*cos(\t r)+1*1.4990963004263174*sin(\t r)});
			\draw [shift={(-0.9495668258430225,0.1532615715811343)},line width=1pt]  plot[domain=0.23635377837499705:3.3779464319647903,variable=\t]({1*1.7172491310717157*cos(\t r)+0*1.7172491310717157*sin(\t r)},{0*1.7172491310717157*cos(\t r)+1*1.7172491310717157*sin(\t r)});
			\draw [shift={(2.2076017168575475,-0.34113587962599384)},line width=1pt]  plot[domain=-0.5423499168248593:2.5992427367649342,variable=\t]({1*1.7369120931874662*cos(\t r)+0*1.7369120931874662*sin(\t r)},{0*1.7369120931874662*cos(\t r)+1*1.7369120931874662*sin(\t r)});
			\draw [shift={(4.0480217056481695,-1.8581793259782755)},line width=1pt]  plot[domain=-1.0538810780443164:2.0877115755454767,variable=\t]({1*0.7137949517702474*cos(\t r)+0*0.7137949517702474*sin(\t r)},{0*0.7137949517702474*cos(\t r)+1*0.7137949517702474*sin(\t r)});
			\draw [shift={(4.449352220961291,-4.049266311043521)},line width=1pt]  plot[domain=-1.539879136206313:1.6017135173834798,variable=\t]({1*1.5713018365705365*cos(\t r)+0*1.5713018365705365*sin(\t r)},{0*1.5713018365705365*cos(\t r)+1*1.5713018365705365*sin(\t r)});
			\draw [shift={(3.545336127535303,-6.9020380363225495)},line width=1pt]  plot[domain=-2.209751372244157:0.931841281345636,variable=\t]({1*1.5973463104482766*cos(\t r)+0*1.5973463104482766*sin(\t r)},{0*1.5973463104482766*cos(\t r)+1*1.5973463104482766*sin(\t r)});
			\draw [shift={(0.7412420079711007,-8.49886940521424)},line width=1pt]  plot[domain=-2.9732788127777363:0.1683138408120568,variable=\t]({1*1.8780448660652638*cos(\t r)+0*1.8780448660652638*sin(\t r)},{0*1.8780448660652638*cos(\t r)+1*1.8780448660652638*sin(\t r)});
			\draw [shift={(-2.4157166199103353,-7.9712820379537295)},line width=1pt]  plot[domain=2.5686426736185357:5.710235327208329,variable=\t]({1*1.5535459902866748*cos(\t r)+0*1.5535459902866748*sin(\t r)},{0*1.5535459902866748*cos(\t r)+1*1.5535459902866748*sin(\t r)});
			\draw [shift={(-4.25942364673348,-5.385387113876428)},line width=1pt]  plot[domain=1.8702022711174957:5.011794924707289,variable=\t]({1*1.8248826025644669*cos(\t r)+0*1.8248826025644669*sin(\t r)},{0*1.8248826025644669*cos(\t r)+1*1.8248826025644669*sin(\t r)});
			\draw [shift={(-4.1877514219575716,-2.3330038550945784)},line width=1pt]  plot[domain=1.1346675632134142:4.276260216803207,variable=\t]({1*1.4438385074221678*cos(\t r)+0*1.4438385074221678*sin(\t r)},{0*1.4438385074221678*cos(\t r)+1*1.4438385074221678*sin(\t r)});
			\draw [shift={(-2.069851044883553,-0.25213553896994895)},line width=1pt]  plot[domain=0.4732536051223898:3.614846258712183,variable=\t]({1*1.6941816301923522*cos(\t r)+0*1.6941816301923522*sin(\t r)},{0*1.6941816301923522*cos(\t r)+1*1.6941816301923522*sin(\t r)});
			\draw [shift={(1.1155093406422365,0.13881293272531114)},line width=1pt]  plot[domain=-0.22348196253436914:2.9181106910554244,variable=\t]({1*1.7201637632878914*cos(\t r)+0*1.7201637632878914*sin(\t r)},{0*1.7201637632878914*cos(\t r)+1*1.7201637632878914*sin(\t r)});
			\draw [shift={(3.850574920982233,-1.8835098837153774)},line width=1pt]  plot[domain=-0.9982983750677557:2.1432942785220375,variable=\t]({1*1.9523984288124627*cos(\t r)+0*1.9523984288124627*sin(\t r)},{0*1.9523984288124627*cos(\t r)+1*1.9523984288124627*sin(\t r)});
			\begin{scriptsize}
				\draw [fill=black] (-0.86,-6.53) circle (2.5pt);
				\draw [fill=black] (0.84,-6.57) circle (2.5pt);
				\draw [fill=black] (2.1679870576897233,-5.507902841257643) circle (2.5pt);
				\draw [fill=black] (2.502581269843993,-3.8406755882435677) circle (2.5pt);
				\draw [fill=black] (1.6872222859953718,-2.3484324018100224) circle (2.5pt);
				\draw [fill=black] (0.10342563639235468,-1.7294104533249488) circle (2.5pt);
				\draw [fill=black] (-1.5077326246767149,-2.273256992147148) circle (2.5pt);
				\draw [fill=black] (-2.3923736408280933,-3.725500178580693) circle (2.5pt);
				\draw [fill=black] (-2.136564048914801,-5.406619285808124) circle (5.5pt);
				\draw [fill=black] (0.10017224365241753,-0.22941398151720493) circle (2.5pt);
				\draw [fill=black] (2.6952636896379882,-1.2376432550579382) circle (2.5pt);
				\draw [fill=black] (3.9224196569869996,-3.356879190577002) circle (2.5pt);
				\draw [fill=black] (3.579512478375505,-6.015441989262826) circle (2.5pt);
				\draw [fill=black] (1.6588207160472934,-7.826794587421428) circle (2.5pt);
				\draw [fill=black] (-1.493308308829422,-7.8897501924867) circle (2.5pt);
				\draw [fill=black] (-3.425541485037292,-6.17378104786818) circle (5.5pt);
				\draw [fill=black] (-3.8430331937006743,-3.343941714306183) circle (2.5pt);
				\draw [fill=black] (-2.603329719814657,-1.2487244132970154) circle (2.5pt);
				\draw [fill=black] (0.7199397440771068,0.5553714958059505) circle (2.5pt);
				\draw [fill=black] (-0.561876818362425,0.5200465174320507) circle (2.5pt);
				\draw [fill=black] (3.6952636896379882,-1.2376432550579382) circle (2.5pt);
				\draw [fill=black] (2.7928954996468978,-0.24242065198142837) circle (2.5pt);
				\draw [fill=black] (4.908254342317568,-3.5245991154493264) circle (2.5pt);
				\draw [fill=black] (4.400779721658352,-2.478715396898613) circle (2.5pt);
				\draw [fill=black] (4.30258408754147,-6.706215069237823) circle (2.5pt);
				\draw [fill=black] (4.497924720264232,-5.619817225188429) circle (2.5pt);
				\draw [fill=black] (2.5927475348063744,-8.18425884745667) circle (2.5pt);
				\draw [fill=black] (1.5124801628682145,-8.816028858206282) circle (2.5pt);
				\draw [fill=black] (-1.110263518864173,-8.813479962971812) circle (2.5pt);
				\draw [fill=black] (-2.194024813288919,-8.603190015021204) circle (2.5pt);
				\draw [fill=black] (-3.721169720956498,-7.129084112935647) circle (5.5pt);
				\draw [fill=black] (-4.423336452189108,-6.240152755142491) circle (5.5pt);
				\draw [fill=black] (-4.797677572510462,-3.6416901148172083) circle (2.5pt);
				\draw [fill=black] (-4.496183938815856,-2.5867138066353927) circle (2.5pt);
				\draw [fill=black] (-3.5778252714046817,-1.0243175953719486) circle (2.5pt);
				\draw [fill=black] (-2.619073395763152,-0.24884835264368188) circle (2.5pt);
			\end{scriptsize}
		\end{tikzpicture}
		\caption{The labeling of the flower snarks graph $J_9$ derived from the definition of the function $f$.  The ``central'' vertices $a_4,b_4,c_4$ and $d_4$ are  highlighted, with their respective labels shown in brackets.} \label{fig:upper bound J_9}
	\end{figure}
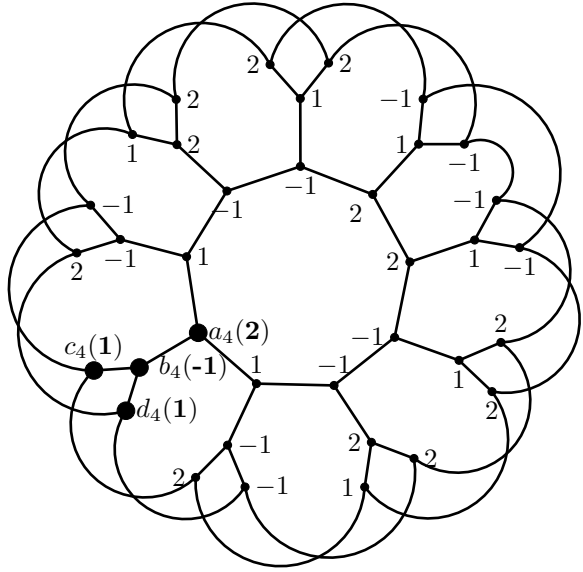
		\begin{figure}
		\centering
		\begin{tikzpicture}[scale=0.45, line cap=round,line join=round,>=triangle 45,x=1cm,y=1cm]
			\draw [line width=1pt] (-2.96,-3.44) node[above]{$1$}-- (-0.9,-3.46);
			\draw [line width=1pt] (-0.9,-3.46) node[above]{$-1$}-- (0.8437950939813437,-2.3631049866980938);
			\draw [line width=1pt] (0.8437950939813437,-2.3631049866980938) node[above]{$2$}-- (1.7177426606723198,-0.4975713765278249);
			\draw [line width=1pt] (1.7177426606723198,-0.4975713765278249) node[left]{$2$}-- (1.4443705226669716,1.5443070905123617);
			\draw [line width=1pt] (1.4443705226669716,1.5443070905123617) node[left]{$-1$}-- (0.11047240222677068,3.1142484283610394);
			\draw [line width=1pt] (0.11047240222677068,3.1142484283610394) node[below left]{$-1$}-- (-1.8604484722822634,3.713807354926675);
			\draw [line width=1pt] (-1.8604484722822634,3.713807354926675) node[below]{$2$}-- (-3.8426386490649564,3.15262814730562);
			\draw [line width=1pt] (-3.8426386490649564,3.15262814730562) node[below]{$2$}-- (-5.206766752479329,1.6088812388147546);
			\draw [line width=1pt] (-5.206766752479329,1.6088812388147546) node[right]{$-1$}-- (-5.5197317481599155,-0.4273046346944991);
			\draw [line width=1pt] (-5.5197317481599155,-0.4273046346944991) node[right]{$1$}-- (-4.68216946128312,-2.309454845384846);
			\draw [line width=1pt] (-4.68216946128312,-2.309454845384846) node[above right]{$a_6(\textbf{2})$}-- (-2.96,-3.44);
			\draw [line width=1pt] (0.11047240222677068,3.1142484283610394) -- (1.4066828078550726,4.637350061326126) node[below]{$1$};
			\draw [line width=1pt] (1.4443705226669716,1.5443070905123617)-- (3.365338663457841,2.1009770125186344) node[below]{$1$};
			\draw [line width=1pt] (1.7177426606723198,-0.4975713765278249)-- (3.6869738254455635,-0.8470404305023007) node[below]{$1$};
			\draw [line width=1pt] (0.8437950939813437,-2.3631049866980938)-- (2.6448570401658955,-3.232688718150054) node[below]{$1$};
			\draw [line width=1pt] (-0.9,-3.46)-- (0.15419012744811444,-5.159612654457164) node[left]{$1$};
			\draw [line width=1pt] (-2.96,-3.44)-- (-3.4689865985592876,-5.374149074525293) node[left]{$-1$};
			\draw [line width=1pt] (-4.68216946128312,-2.309454845384846)-- (-6.21525154164667,-3.5938454926834664) node[right]{$\ b_6(\textbf{-1})$};
			\draw [line width=1pt] (-5.5197317481599155,-0.4273046346944991)-- (-7.430773320893657,-1.0171521635384795) node[above]{$-1$};
			\draw [line width=1pt] (-5.206766752479329,1.6088812388147546)-- (-7.080965283688063,2.307007205682498) node[below]{$1$};
			\draw [line width=1pt] (-3.8426386490649564,3.15262814730562)-- (-5.1163766536628374,4.694572211026159) node[below]{$1$};
			\draw [line width=1pt] (-1.8604484722822634,3.713807354926675)-- (-1.9136549383894643,5.713099497652291) node[left]{$1$};
			\draw [line width=1pt] (1.4066828078550726,4.637350061326126)-- (2.7066828078550724,4.637350061326126) node[below]{$2$};
			\draw [line width=1pt] (1.4066828078550726,4.637350061326126)-- (1.4605980731975976,5.936231558677999) node[left]{$2$};
			\draw [line width=1pt] (3.365338663457841,2.1009770125186344)-- (4.657002939723892,1.9539960701550085) node[below]{$2$};
			\draw [line width=1pt] (3.365338663457841,2.1009770125186344)-- (3.9429444442253914,3.265610672526754) node[left]{$2$};
			\draw [line width=1pt] (3.6869738254455635,-0.8470404305023007)-- (4.752009121882505,-0.10158767053459106) node[above]{$-1$};
			\draw [line width=1pt] (3.6869738254455635,-0.8470404305023007)-- (4.738394846006681,-1.6115754771301312) node[below]{$-1$};
			\draw [line width=1pt] (2.6448570401658955,-3.232688718150054)-- (3.9266364514783123,-3.0157980988609885) node[above]{$-1$};
			\draw [line width=1pt] (2.6448570401658955,-3.232688718150054)-- (3.227399085091234,-4.39486115195418) node[left]{$-1$};
			\draw [line width=1pt] (0.15419012744811444,-5.159612654457164)-- (1.356145301913671,-5.654894148736401) node[above]{$2$};
			\draw [line width=1pt] (0.15419012744811444,-5.159612654457164)-- (-0.10376844672090513,-6.433762321495023) node[left]{$2$};
			\draw [line width=1pt] (-3.4689865985592876,-5.374149074525293)-- (-4.378684557303566,-6.302830733575633) node[left]{$1$};
			\draw [line width=1pt] (-3.4689865985592876,-5.374149074525293)-- (-2.8148305409168897,-6.497572348385042) node[right]{$2$};
			\draw [line width=1pt] (-6.21525154164667,-3.5938454926834664)-- (-6.484763386597155,-4.865601509154938) node[right]{$d_6(\textbf{1})$};
			\draw [line width=1pt] (-6.21525154164667,-3.5938454926834664)-- (-7.500027970810292,-3.792212644867425) node[above]{$c_6(\textbf{1})$};
			\draw [line width=1pt] (-7.430773320893657,-1.0171521635384795)-- (-8.446322667110964,-1.8287305795880933) node[below]{$2$};
			\draw [line width=1pt] (-7.430773320893657,-1.0171521635384795)-- (-8.576670631556278,-0.40321458249609654) node[above]{$1$};
			\draw [line width=1pt] (-7.080965283688063,2.307007205682498)-- (-8.140807274830701,1.5541891972562856) node[below]{$2$};
			\draw [line width=1pt] (-7.080965283688063,2.307007205682498)-- (-8.135191692556972,3.0676692060764728) node[right]{$2$};
			\draw [line width=1pt] (-5.1163766536628374,4.694572211026159)-- (-6.416361486297197,4.700851933998928) node[below]{$-1$};
			\draw [line width=1pt] (-5.1163766536628374,4.694572211026159)-- (-5.303061492088993,5.9810980642843825) node[right]{$-1$};
			\draw [line width=1pt] (-1.9136549383894643,5.713099497652291)-- (-1.128295138892324,6.749057983950758) node[right]{$-1$};
			\draw [line width=1pt] (-1.9136549383894643,5.713099497652291)-- (-2.751929371962511,6.706727181454321) node[left]{$-1$};
			\draw [shift={(4.347476783053948,1.5820115009960816)},line width=1pt]  plot[domain=-1.3349882242182263:1.8066044293715668,variable=\t]({1*1.7315174222446825*cos(\t r)+0*1.7315174222446825*sin(\t r)},{0*1.7315174222446825*cos(\t r)+1*1.7315174222446825*sin(\t r)});
			\draw [shift={(4.339322786680409,-1.5586928846977899)},line width=1pt]  plot[domain=-1.846791606174297:1.2948010474154958,variable=\t]({1*1.5144192340313558*cos(\t r)+0*1.5144192340313558*sin(\t r)},{0*1.5144192340313558*cos(\t r)+1*1.5144192340313558*sin(\t r)});
			\draw [shift={(2.6413908766959917,-4.335346123798694)},line width=1pt]  plot[domain=-2.3430262808437554:0.7985663727460377,variable=\t]({1*1.8420269209798588*cos(\t r)+0*1.8420269209798588*sin(\t r)},{0*1.8420269209798588*cos(\t r)+1*1.8420269209798588*sin(\t r)});
			\draw [shift={(-0.7293426195016093,-6.076233248560721)},line width=1pt]  plot[domain=-2.942242260733796:0.19935039285599698,variable=\t]({1*2.1276246161881556*cos(\t r)+0*2.1276246161881556*sin(\t r)},{0*2.1276246161881556*cos(\t r)+1*2.1276246161881556*sin(\t r)});
			\draw [shift={(-4.649796963757023,-5.68158692876999)},line width=1pt]  plot[domain=2.7231659296485295:5.864758583238323,variable=\t]({1*2.0082166162979194*cos(\t r)+0*2.0082166162979194*sin(\t r)},{0*2.0082166162979194*cos(\t r)+1*2.0082166162979194*sin(\t r)});
			\draw [shift={(-7.46554302685406,-3.3471660443715154)},line width=1pt]  plot[domain=2.14429420000467:5.285886853594463,variable=\t]({1*1.8076434835040094*cos(\t r)+0*1.8076434835040094*sin(\t r)},{0*1.8076434835040094*cos(\t r)+1*1.8076434835040094*sin(\t r)});
			\draw [shift={(-8.293564970970833,-0.13727069116590385)},line width=1pt]  plot[domain=1.4807295063046715:4.622322159894464,variable=\t]({1*1.6983437425536818*cos(\t r)+0*1.6983437425536818*sin(\t r)},{0*1.6983437425536818*cos(\t r)+1*1.6983437425536818*sin(\t r)});
			\draw [shift={(-7.278584380563949,3.1275205656276066)},line width=1pt]  plot[domain=1.0694716868096172:4.21106434039941,variable=\t]({1*1.7941014224672218*cos(\t r)+0*1.7941014224672218*sin(\t r)},{0*1.7941014224672218*cos(\t r)+1*1.7941014224672218*sin(\t r)});
			\draw [shift={(-4.584145429129855,5.703789557726624)},line width=1pt]  plot[domain=0.5008375319185491:3.642430185508342,variable=\t]({1*2.0887555044165422*cos(\t r)+0*2.0887555044165422*sin(\t r)},{0*2.0887555044165422*cos(\t r)+1*2.0887555044165422*sin(\t r)});
			\draw [shift={(-0.6456656493824567,6.32147937006616)},line width=1pt]  plot[domain=-0.1809060940267777:2.9606865595630154,variable=\t]({1*2.1412059091166724*cos(\t r)+0*2.1412059091166724*sin(\t r)},{0*2.1412059091166724*cos(\t r)+1*2.1412059091166724*sin(\t r)});
			\draw [shift={(3.0588005064607446,3.945113814416504)},line width=1pt]  plot[domain=-0.8944319126437428:2.2471607409460503,variable=\t]({1*2.553194252147931*cos(\t r)+0*2.553194252147931*sin(\t r)},{0*2.553194252147931*cos(\t r)+1*2.553194252147931*sin(\t r)});
			\draw [shift={(4.6976988928652865,0.1712102965124387)},line width=1pt]  plot[domain=-1.5479731209219079:1.5936195326678853,variable=\t]({1*1.7832501998610355*cos(\t r)+0*1.7832501998610355*sin(\t r)},{0*1.7832501998610355*cos(\t r)+1*1.7832501998610355*sin(\t r)});
			\draw [shift={(3.9828969655489574,-3.0032183145421554)},line width=1pt]  plot[domain=-2.068158260831911:1.073434392757882,variable=\t]({1*1.5834919748127243*cos(\t r)+0*1.5834919748127243*sin(\t r)},{0*1.5834919748127243*cos(\t r)+1*1.5834919748127243*sin(\t r)});
			\draw [shift={(1.5618153191851645,-5.414311736724601)},line width=1pt]  plot[domain=-2.592346820276365:0.549245833313428,variable=\t]({1*1.9528053605105142*cos(\t r)+0*1.9528053605105142*sin(\t r)},{0*1.9528053605105142*cos(\t r)+1*1.9528053605105142*sin(\t r)});
			\draw [shift={(-2.2412265020122355,-6.368296527535328)},line width=1pt]  plot[domain=3.110974352151567:6.25256700574136,variable=\t]({1*2.138460359302591*cos(\t r)+0*2.138460359302591*sin(\t r)},{0*2.138460359302591*cos(\t r)+1*2.138460359302591*sin(\t r)});
			\draw [shift={(-5.939356264056929,-5.047521689221529)},line width=1pt]  plot[domain=2.464211603288498:5.605804256878291,variable=\t]({1*2.00287218091861*cos(\t r)+0*2.00287218091861*sin(\t r)},{0*2.00287218091861*cos(\t r)+1*2.00287218091861*sin(\t r)});
			\draw [shift={(-8.038349301183285,-2.0977136136817607)},line width=1pt]  plot[domain=1.878400260241021:5.019992913830814,variable=\t]({1*1.7779529862804875*cos(\t r)+0*1.7779529862804875*sin(\t r)},{0*1.7779529862804875*cos(\t r)+1*1.7779529862804875*sin(\t r)});
			\draw [shift={(-8.355931162056624,1.3322273117901882)},line width=1pt]  plot[domain=1.4442807238436541:4.5858733774334475,variable=\t]({1*1.7494240428892465*cos(\t r)+0*1.7494240428892465*sin(\t r)},{0*1.7494240428892465*cos(\t r)+1*1.7494240428892465*sin(\t r)});
			\draw [shift={(-6.719126592322983,4.524383635180428)},line width=1pt]  plot[domain=0.7995470446038112:3.9411396981936044,variable=\t]({1*2.0315652330310146*cos(\t r)+0*2.0315652330310146*sin(\t r)},{0*2.0315652330310146*cos(\t r)+1*2.0315652330310146*sin(\t r)});
			\draw [shift={(-3.2156783154906585,6.365078024117571)},line width=1pt]  plot[domain=0.18191902189975312:3.3235116754895464,variable=\t]({1*2.12240644917491*cos(\t r)+0*2.12240644917491*sin(\t r)},{0*2.12240644917491*cos(\t r)+1*2.12240644917491*sin(\t r)});
			\draw [shift={(0.7891938344813743,5.693204022638442)},line width=1pt]  plot[domain=-0.5033375755019058:2.638255078087887,variable=\t]({1*2.188970431647862*cos(\t r)+0*2.188970431647862*sin(\t r)},{0*2.188970431647862*cos(\t r)+1*2.188970431647862*sin(\t r)});
			\draw [shift={(3.324813626040232,3.95148036692644)},line width=1pt]  plot[domain=-0.8372984833603416:2.3042941702294515,variable=\t]({1*0.9233108610247003*cos(\t r)+0*0.9233108610247003*sin(\t r)},{0*0.9233108610247003*cos(\t r)+1*0.9233108610247003*sin(\t r)});
			\begin{scriptsize}
				\draw [fill=black] (-2.96,-3.44) circle (2.5pt);
				\draw [fill=black] (-0.9,-3.46) circle (2.5pt);
				\draw [fill=black] (0.8437950939813437,-2.3631049866980938) circle (2.5pt);
				\draw [fill=black] (1.7177426606723198,-0.4975713765278249) circle (2.5pt);
				\draw [fill=black] (1.4443705226669716,1.5443070905123617) circle (2.5pt);
				\draw [fill=black] (0.11047240222677068,3.1142484283610394) circle (2.5pt);
				\draw [fill=black] (-1.8604484722822634,3.713807354926675) circle (2.5pt);
				\draw [fill=black] (-3.8426386490649564,3.15262814730562) circle (2.5pt);
				\draw [fill=black] (-5.206766752479329,1.6088812388147546) circle (2.5pt);
				\draw [fill=black] (-5.5197317481599155,-0.4273046346944991) circle (2.5pt);
				\draw [fill=black] (-4.68216946128312,-2.309454845384846) circle (5.5pt);
				\draw [fill=black] (1.4066828078550726,4.637350061326126) circle (2.5pt);
				\draw [fill=black] (3.365338663457841,2.1009770125186344) circle (2.5pt);
				\draw [fill=black] (3.6869738254455635,-0.8470404305023007) circle (2.5pt);
				\draw [fill=black] (2.6448570401658955,-3.232688718150054) circle (2.5pt);
				\draw [fill=black] (0.15419012744811444,-5.159612654457164) circle (2.5pt);
				\draw [fill=black] (-3.4689865985592876,-5.374149074525293) circle (2.5pt);
				\draw [fill=black] (-6.21525154164667,-3.5938454926834664) circle (5.5pt);
				\draw [fill=black] (-7.430773320893657,-1.0171521635384795) circle (2.5pt);
				\draw [fill=black] (-7.080965283688063,2.307007205682498) circle (2.5pt);
				\draw [fill=black] (-5.1163766536628374,4.694572211026159) circle (2.5pt);
				\draw [fill=black] (-1.9136549383894643,5.713099497652291) circle (2.5pt);
				\draw [fill=black] (2.7066828078550724,4.637350061326126) circle (2.5pt);
				\draw [fill=black] (1.4605980731975976,5.936231558677999) circle (2.5pt);
				\draw [fill=black] (4.657002939723892,1.9539960701550085) circle (2.5pt);
				\draw [fill=black] (3.9429444442253914,3.265610672526754) circle (2.5pt);
				\draw [fill=black] (4.752009121882505,-0.10158767053459106) circle (2.5pt);
				\draw [fill=black] (4.738394846006681,-1.6115754771301312) circle (2.5pt);
				\draw [fill=black] (3.9266364514783123,-3.0157980988609885) circle (2.5pt);
				\draw [fill=black] (3.227399085091234,-4.39486115195418) circle (2.5pt);
				\draw [fill=black] (1.356145301913671,-5.654894148736401) circle (2.5pt);
				\draw [fill=black] (-0.10376844672090513,-6.433762321495023) circle (2.5pt);
				\draw [fill=black] (-4.378684557303566,-6.302830733575633) circle (2.5pt);
				\draw [fill=black] (-2.8148305409168897,-6.497572348385042) circle (2.5pt);
				\draw [fill=black] (-6.484763386597155,-4.865601509154938) circle (5.5pt);
				\draw [fill=black] (-7.500027970810292,-3.792212644867425) circle (5.5pt);
				\draw [fill=black] (-8.446322667110964,-1.8287305795880933) circle (2.5pt);
				\draw [fill=black] (-8.576670631556278,-0.40321458249609654) circle (2.5pt);
				\draw [fill=black] (-8.140807274830701,1.5541891972562856) circle (2.5pt);
				\draw [fill=black] (-8.135191692556972,3.0676692060764728) circle (2.5pt);
				\draw [fill=black] (-6.416361486297197,4.700851933998928) circle (2.5pt);
				\draw [fill=black] (-5.303061492088993,5.9810980642843825) circle (2.5pt);
				\draw [fill=black] (-1.128295138892324,6.749057983950758) circle (2.5pt);
				\draw [fill=black] (-2.751929371962511,6.706727181454321) circle (2.5pt);
			\end{scriptsize}
		\end{tikzpicture}
		\caption{The labeling of the flower snarks graph $J_{11}$ derived from the definition of the function $f$. The ``central'' vertices $a_6,b_6,c_6$ and $d_6$  are  highlighted, with their respective labels shown in brackets.} \label{fig:upper bound J_{11}}
	\end{figure}
			
	\begin{table}[h]
		
		\caption{The neighbor structure of vertices in the flower snarks graph $J_{2n+1}$, with respect to the labeling given by the function $f$, for even $n\geqslant 2$. 
			Vertices are grouped based on their assignments and  cardinalities of their neighbor sets intersected with $V_{-1}$, $V_1$ and $V_2$. respectively.} \label{opt-J_{2n+1},n even}
		\begin{tabular}{l|c|c|c|c|c } \hline
			Vertex $v$ &$f(v)$& $\left|N(v) \cap V_{-1}\right|$  & $\left|N(v) \cap V_{1}\right|$  & $\left|N(v) \cap V_{2}\right|$ & $\sum_{u\in N[v]} f(u)$ \\ \hline
			\multirow{2}{70pt}{$a_n$, $c_{n+1}$, $c_{n+3}$, $d_{n-1}$, $d_{n-3}$}   &   &  &   & &    \\
			&  $2$ & $1$  &  $2$  & $0$ & $3$    \\
			\hline
			\multirow{3}{90pt}{$a_{n-1}$, $a_{n+1}$, $c_n$, $d_n$, $b_{n-2-j}$, $b_{n+2+j}$, {\tiny $j\equiv 2(\bmod 4)$ or $j\equiv 3(\bmod 4)$}} &  &&  &  &        \\
			&   &   &   &   &        \\ 
			& $1$ &  $2$  & $0$  & $1$   &  $1$\\
			\hline 
			\multirow{7}{90pt}{$b_{n-1}$, $b_{n+1}$, $c_{n-1}$, $d_{n+1}$, $a_{n-2-i}$, $a_{n+2+i}$, {\tiny $i\equiv 0(\bmod 4)$ or $i\equiv 1(\bmod 4)$}, $c_{n-2-j}$, $c_{n+2+j}$, $d_{n-2-j}$, $d_{n+2+j}$, {\tiny $j\equiv 2(\bmod 4)$ or $j\equiv 3(\bmod 4)$} }& & & & & \\
			&   & & & & \\
			&  & & & & \\
			& $-1$ & $1$ & $1$ & $1$ & $1$ \\
			&  & & & & \\
			&  & & & & \\
			&  & & & & \\ 
			\hline
			\multirow{6}{91pt}{$b_{n-2}$, $b_{n+2}$, $a_{n-2-i}$, $a_{n+2+i}$, {\tiny $i\equiv 2(\bmod 4)$ or $i\equiv 3(\bmod 4)$}, $c_{n-2-j}$, $c_{n+2+j}$, $d_{n-2-j}$, $d_{n+2+j}$, \tiny{$j\geqslant 1$ and} {\tiny $j\equiv 0(\bmod 4)$ or $j\equiv 1(\bmod 4)$}} && & & & \\
			&  & & & & \\
			&  & & & & \\
			& $2$ & $1$ & $1$ & $1$ & $4$\\
			&  & & & & \\
			&  & & & & \\   
			\hline 
			$b_n$&$-1$ & $0$ & $2$ & $1$ & $3$\\
			\hline
			\multirow{2}{95pt}{$b_{n-2-j}$, $b_{n+2+j}$, \tiny{$j\geqslant 1$ and}  {\tiny $j\equiv 0(\bmod 4)$ or $j\equiv 1(\bmod 4)$} }& & & & \\
			&$1$  & $1$ & $0$ & $2$ & $4$ \\
			\hline
			$c_{n-2}$, $d_{n+2}$&$2$ & $1$ & $0$ & $2$ & $5$\\
			\hline
			$c_{n+2}$, $d_{n-2}$ &$1$& $0$ & $0$ & $3$ & $7$\\
			\hline
		\end{tabular}

	\end{table}
	
	In Table \ref{opt-J_{2n+1},n even}, we provide an overview of the neighbor structure for each vertex of flower snarks graph $J_{2n+1}$, according to the labeling given by the function $f$, for even $n\geqslant 2$. This table demonstrates that the conditions (C1)--(C3) from Definition~\ref{def:1} are satisfied. Specifically, each vertex assigned a value of $-1$ has exactly one (exclusive) neighbor assigned a value of $2$, and the total sum of labels in the closed neighborhood of each vertex is positive.
	
	
	
	This means that $f$ is a valid $PLSRD$ function on a graph $J_{2n+1}$, for even $n\geqslant 2$. Furthermore, in this case, we have
	\begin{align*}
		\sum_{v\in V(J_{2n+1})} f(v)&=\sum_{i=0}^{2n} f(a_i)+ \sum_{i=0}^{2n} f(b_i)+\sum_{i=0}^{2n} f(c_i)+ \sum_{i=0}^{2n} f(d_i)\\
		&=((n-2)+2)+(2(n-2)+1)+((n-2)+5)+((n-2)+5)\\
		&=5n+3.
	\end{align*} 
	
	$\bullet$ Case 2. $n\geqslant 2$ is an odd number;
	
	We define 
	\begin{align*}
		&f(a_n)=2, f(a_{n-1})=f(a_{n+1})=1,  \\
		&f(b_n)=f(b_{n-1})=f(b_{n+1})=-1,\\
		&f(c_n)=f(c_{n-1})=f(c_{n+1})=1,\\
		&f(d_n)=1, f(d_{n-1})=f(d_{n+1})=2.
	\end{align*}
	Also, for  $j\in\{1,2,\dotsc,n-1\}$, we put
	\begin{align*}
		&f(a_{n-1-j})=\begin{cases}
			-1, &\mbox{ for } j\equiv 0(\bmod\ 4) \mbox{ or } j\equiv 1(\bmod\ 4),\\
			2  &\mbox{ otherwise},
		\end{cases}\\
		&f(b_{n-1-j})=1,\\
		&f(c_{n-1-j})=f(d_{n-1-j})=\begin{cases}
			2, &\mbox{ for } j\equiv 0(\bmod\ 4) \mbox{ or } j\equiv 1(\bmod\ 4),\\
			-1  &\mbox{ otherwise}.
		\end{cases}
	\end{align*} 
	The remaining vertices are labeled in such a way that symmetry relative to the vertices $a_n$, $b_n$, $c_n$ and $d_n$ is preserved. More precisely, for $j\in\{1,2,\dotsc,n-1\}$, we set 
	\begin{align*}
		&f(a_{n+1+j})=f(a_{n-1-j}), f(b_{n+1+j})=f(b_{n-1-j}),\\
		&f(c_{n+1+j})=f(c_{n-1-j}), f(d_{n+1+j})=f(d_{n-1-j}).	
	\end{align*}
	An example of the labeling of the flower snarks graph, using the previously described pattern, is given in the Figure \ref{fig:upper bound J_{11}}.

	In Table \ref{opt-J_{2n+1},n odd},  an overview of the neighbor structure for each vertex of the flower snarks graph $J_{2n+1}$, based on the labeling given by the function $f$, for odd $n\geqslant 2$ is presented. From the table, it is evident that conditions (C1)--(C3) from Definition~\ref{def:1} are satisfied and  $f$ is a valid $PLSRD$ function on a graph $J_{2n+1}$, for odd $n\geqslant 2$. Specifically, each vertex assigned a value of -1 has exactly one exclusive neighbor assigned a value of 2, and the total sum of labels in the closed neighborhood of each vertex remains positive.
	
	
	
	Furthermore, in this case, we have
	\begin{align*}
		f(V(J_{2n+1}))& = \sum_{v\in V(J_{2n+1})} f(v)=\sum_{i=0}^{2n} f(a_i)+ \sum_{i=0}^{2n} f(b_i)+\sum_{i=0}^{2n} f(c_i)+ \sum_{i=0}^{2n} f(d_i)\\
		&=((n-1)+4)+(2(n-1)-3)+((n-1)+3)+((n-1)+5)\\
		&=5n+4.
	\end{align*}
	
	Thus, 
	$$ \gamma_{plsR}(J_{2n+1})\leqslant
	\begin{cases}
		5n+3, &\mbox{for even number } n,\\
		5n+4  &\mbox{for odd number } n.
	\end{cases}$$

\begin{table}[h]

	\caption{The neighbor structure of vertices in the flower snarks graph $J_{2n+1}$, with respect to the labeling given by the function $f$, for odd $n\geqslant 2$. 
		Vertices are grouped based on their assignments and  cardinalities of their neighbor sets intersected with $V_{-1}$, $V_1$ and $V_2$. respectively.}  \label{opt-J_{2n+1},n odd}
	\begin{tabular}{l|c|c|c|c|c } \hline
		Vertex $v$ &$f(v)$& $\left|N(v) \cap V_{-1}\right|$  & $\left|N(v) \cap V_{1}\right|$  & $\left|N(v) \cap V_{2}\right|$ & $\sum_{u\in N[v]} f(u)$ \\ \hline
		$a_n$, $c_{n-2}$, $c_{n+2}$ &$2$& $1$  &  $2$  & $0$ & $3$    \\
		\hline
		\multirow{3}{90pt}{$a_{n-1}$, $a_{n+1}$, $b_{n-1-j}$, $b_{n+1+j}$, {\tiny $j\equiv 2(\bmod 4)$ or $j\equiv 3(\bmod 4)$}}& &  &  &  &        \\
		&  & &   &   &        \\ 
		&$1$&  $2$  & $0$  & $1$   &  $1$\\
		\hline 
		\multirow{6}{90pt}{$a_{n-1-i}$, $a_{n+1+i}$, {\tiny $i\geqslant 2$ and $i\equiv 0(\bmod 4)$ or $i\equiv 1(\bmod 4)$}, $c_{n-1-j}$, $c_{n+1+j}$, $d_{n-1-j}$, $d_{n+1+j}$, {\tiny $j\equiv 2(\bmod 4)$ or $j\equiv 3(\bmod 4)$} }& & & & & \\
		&	& & & & \\
		&	& & & & \\
		&$-1$	& $1$ & $1$ & $1$ & $1$ \\
		&	& & & & \\
		&	& & & & \\
		\hline
		\multirow{7}{91pt}{ $a_{n-1-i}$, $a_{n+1+i}$, {\tiny $i\equiv 2(\bmod 4)$ or $i\equiv 3(\bmod 4)$}, $c_{n-1-j}$, $c_{n+1+j}$, {\tiny$j\geqslant 2$ and} {\tiny $j\equiv 0(\bmod 4)$ or $j\equiv 1(\bmod 4)$}$, d_{n-1-l}$, $d_{n+1+l}$,  {\tiny $l\equiv 0(\bmod 4)$ or $l\equiv 1(\bmod 4)$}} & && & & \\
		&	& & & & \\
		&	& & & & \\
		&$2$	& $1$ & $1$ & $1$ & $4$\\
		&	& & & & \\
		&	& & & & \\
		&	& & & & \\
		\hline 
		$c_{n-1}$, $c_{n+1}$	& $1$	&$1$	&$1$ &$1$ &$3$\\ \hline
		\multirow{2}{80pt}{$a_{n-2}$, $a_{n+2}$, $b_n$, $b_{n-1}$, $b_{n+1}$}& & & & & \\
		&	$-1$ & $0$ & $2$ & $1$ & $3$\\
		\hline
		\multirow{3}{95pt}{$d_n$, $b_{n-1-j}$, $b_{n+1+j}$, \tiny{$j\geqslant 1$ and}  {\tiny $j\equiv 0(\bmod 4)$ or $j\equiv 1(\bmod 4)$} }& & & & \\
		&$1$	& $1$ & $0$ & $2$ & $4$ \\
		&	& & & &\\
		\hline
		$c_{n}$&$1$ & $1$ & $2$ & $0$ & $2$\\
		\hline
	\end{tabular}

\end{table}

\end{proof}

\section{Conclusions}
This paper presents a study of the Perfect Location Signed Roman Domination (PLSRD) problem, a novel variant in the domain of Roman domination problems. The considered problem combines elements of Perfect Roman, Locating Roman, and Signed Roman Domination. Through detailed analysis, we provided exact values for the PLSRD number in well-known graph classes such as complete graphs, complete bipartite graphs, wheels, paths, cycles, ladders, prism graphs, and $3 \times n$ grids. In addition to these results,  a lower bound for 3 regular graphs and an upper bound for flower snark graphs  have been determined.

The results contribute to a deeper understanding of how the interplay between weak and strong vertices can affect the optimization of the PLSRD function, setting the groundwork for further exploration. In future research, a natural direction would be to extend these findings to other graph classes. In particular, exploring additional classes of 3-regular graphs could yield valuable insights, as the proposed lower bound was shown to be attainable in some instances of prism graphs, while remaining sharp for flower snarks.  Specifically, generalized Petersen graphs $GP(n,k), k\geqslant 2$ and certain classes of 3 regular convex polytopes could be analyzed in greater depth. Another promising avenue for further investigation is the development of efficient algorithms for computing the PLSRD number, especially for large and complex graphs where exact values are difficult to compute. Additionally, exploring the relationship between PLSRD and other types of domination problems could provide further insights into their applicability in diverse theoretical and practical contexts.

\backmatter

\bmhead{Acknowledgements}

 This research is partially supported by the grant of Ministry of Scientific and Technological Development and Higher Education, Republika Srpska, Bosnia and Herzegovina, under the project ``Analysis of biological networks by using machine learning methods''.

\begin{appendices}

\section {Detailed Proof for the cases of Theorem~\ref{thm:ladder}}\label{app:A}

In this appendix, we provide a detailed analysis of the cases where a column in the ladder graph $L_n$ requires charging. Specifically, we examine two main types of columns: boundary columns (denoted as type (a)) and interior columns (denoted as type (b)). For each of these types, we describe the specific scenarios that arise and outline how the $\mathrm{CBDP}$ procedure is applied. Using arrows above the scheme, we graphically illustrate the redistribution of charges, leading to proper discharging.

We begin with the boundary column case (type (a)), followed by an exploration of the interior column case (type (b)). Each case is examined thoroughly to ensure that all configurations and behaviors are accounted for.


$\bullet$ \textbf{Type (a)}; 

Due to the symmetry of the structure of the ladder graph $L_n$, it is sufficient to consider only leftmost columns. Column $\begin{tikzcd}
	-1\\[-20pt]
	-1
\end{tikzcd}$ cannot occupy the leftmost position in any feasible solution. Hence, it is sufficient to examine the columns 
$\begin{tikzcd}
	-1\\[-20pt]
	1
\end{tikzcd}$ and
$\begin{tikzcd}
	-1\\[-20pt]
	2
\end{tikzcd}$, as the same CBDP is applied equally to their symmetric counterparts $\begin{tikzcd}
	1\\[-20pt]
	-1
\end{tikzcd}$ and $\begin{tikzcd}
	2\\[-20pt]
	-1
\end{tikzcd}$, respectively. 

\textbf{Case 1.} CBDP applied for the leftmost column $\begin{tikzcd}
	-1\\[-20pt]
	1
\end{tikzcd}$. This column determines the following configuration: 
$\begin{tikzcd}[column sep=tiny]
	-1 & 2 & \geqslant 1\\[-20pt]
	1 & \geqslant 1 & a_1
\end{tikzcd}$. Depending on the value of $a_1$, we distinguish  two cases.

\textbf{Case 1.1.} $a_1\geqslant 1$. CBDP is performed as 
$\begin{tikzcd}[column sep=tiny]
	-1 & 2  \arrow[l,bend left=300,swap,"\frac{6}{5}"]& \geqslant 1\\[-20pt]
	1 & \geqslant 1  & \boxed{\geqslant 1} \\[-20pt]
\end{tikzcd}$,

\textbf{Case 1.2.} $a_1=-1$. We have three subcases, which are directly resolved as follows.:

\textbf{Case 1.2.1.} $2$ is located left from $a_1$: 
$\begin{tikzcd}[column sep=tiny]
	-1 & 2 \arrow[l,bend left=300,swap,"\frac{6}{5}"] \arrow[r,bend left=60,"\frac{6}{5}"]& 1 \\[-20pt]
	1 & 2 & \boxed{-1}
\end{tikzcd}$  

\textbf{Case 1.2.2.} $2$ is located above $a_1$: 
$\begin{tikzcd}[column sep=tiny]
	-1 & 2 \arrow[l,bend left=300,swap,"\frac{6}{5}"] \arrow[r,bend left=60,"\frac{1}{5}"] & 2\\[-20pt]
	1 & 1 & \boxed{-1}
\end{tikzcd}$ 

\textbf{Case 1.2.3.} $2$ is located right from $a_1$: 
$\begin{tikzcd}[column sep=tiny]
	-1 & 2 \arrow[l,bend left=300,swap,"\frac{6}{5}"] \arrow[r,bend left=60,"\frac{3}{5}"]& 1 & \geqslant 1 \arrow[l,bend left=300,swap,"\frac{3}{5}"] \arrow[r,bend left=60, dashrightarrow, "\frac{6}{5}"]& \phantom{-1}\\[-20pt]
	1 & 1 & \boxed{-1} & 2 & \geqslant 1
\end{tikzcd}$

\textbf{Case 2.} CBDP for the leftmost column $\begin{tikzcd}[column sep=tiny]
	-1\\[-20pt]
	2
\end{tikzcd}$. This configuration forces the assignment

$\begin{tikzcd}[column sep=tiny]
	-1 & 1 & a_2\\[-20pt]
	2 & \geqslant 1 & a_3
\end{tikzcd}$. Depending on the value of $a_2 + a_3$, we distinguish  4 cases.

\textbf{Case 2.1.} $a_2+a_3\geqslant 2$. CBDP is performed as 
$\begin{tikzcd}[column sep=tiny]
	-1 & 1 \arrow[l,bend left=300,swap,"\frac{1}{5}"] & \boxed{\geqslant 1}\\[-20pt]
	2 & \geqslant 1 & \boxed{\geqslant 1}
\end{tikzcd}$,

\textbf{Case 2.2.} $a_2+a_3= 1$. CBDP is performed as 
$\begin{tikzcd}[column sep=tiny]
	-1 & 1 \arrow[l,bend left=300,swap,"\frac{1}{5}"] \arrow[r,bend left=60,"\frac{1}{5}"]& \boxed{2}\\[-20pt]
	2 & 1 & \boxed{-1}
\end{tikzcd}$ or
$\begin{tikzcd}[column sep=tiny]
	-1 & 1 \arrow[l,bend left=300,swap,"\frac{1}{5}"] \arrow[r,bend left=60,"\frac{1}{5}"]& \boxed{-1}\\[-20pt]
	2 & 2 & \boxed{2}
\end{tikzcd}$,

\textbf{Case 2.3.} $a_2+a_3= 0$. We distinguish two subcases.

\textbf{Case 2.3.1.} $a_2=-1$ and $a_3=1$. In this case, $2$ must be located left from $a_3$, so we immediately have:  $\begin{tikzcd}[column sep=tiny]
	-1 & 1 \arrow[l,bend left=300,swap,"\frac{1}{5}"] \arrow[r,bend left=60,"\frac{6}{5}"]& \boxed{-1}\\[-20pt]
	2 & 2 & \boxed{1}
\end{tikzcd}$ 

\textbf{Case 2.3.2.} $a_3=-1$ and $a_2 = 1$. We get two subcases, which are immediately resolved as follows.

\textbf{Case 2.3.2.1.} $2$ is left from $a_3$: 
$\begin{tikzcd}[column sep=tiny]
	-1 & 1 \arrow[l,bend left=300,swap,"\frac{1}{5}"] \arrow[r,bend left=60,"\frac{6}{5}"]& \boxed{1}\\[-20pt]
	2 & 2 & \boxed{-1}
\end{tikzcd}$;

\textbf{Case 2.3.2.2.} $2$ is located right from $a_3$: 	 
$\begin{tikzcd}[column sep=tiny]
	-1 & 1 \arrow[l,bend left=300,swap,"\frac{1}{5}"] \arrow[r,bend left=60,"\frac{3}{5}"]& \boxed{1} & \geqslant 1 \arrow[l,bend left=300,swap,"\frac{3}{5}"] \arrow[r,bend left=60,dashrightarrow,"\frac{6}{5}"]& \phantom{-1}\\[-20pt]
	2 & 1 & \boxed{-1} & 2 & \geqslant 1
\end{tikzcd}$

\textbf{Case 2.4.} $a_2+a_3= -2$. We get the configuration 

$\begin{tikzcd}[column sep=tiny]
	-1 & 1 & \boxed{-1} & 2 &\geqslant 1\\[-20pt]
	2 & 2 & \boxed{-1} & 1 & a_4
\end{tikzcd}$. Depending on the value of $a_4$, we have two subcases.

\textbf{Case 2.4.1.} If $a_4\geqslant 1$ we have only one case:
$\begin{tikzcd}[column sep=tiny]
	-1 & 1 \arrow[l,bend left=300,swap,"\frac{1}{5}"] \arrow[r,bend left=60,"\frac{8}{5}"] & -1 & 2 \arrow[l,bend left=300,swap,"\frac{8}{5}"] &\geqslant 1\\[-20pt]
	2 & 2 & -1 & 1 & \boxed{\geqslant 1}
\end{tikzcd}$,

\textbf{Case 2.4.2.} If $a_4=-1$ we distinguish two sub-variants:

\textbf{Case 2.4.2.1.}  $2$ is located above $a_4$:
$\begin{tikzcd}[column sep=tiny]
	-1 & 1 \arrow[l,bend left=300,swap,"\frac{1}{5}"] \arrow[r,bend left=60,"\frac{8}{5}"] & -1 & 2 \arrow[l,bend left=300,swap,"\frac{8}{5}"] \arrow[r,bend left=60,"\frac{1}{5}"] & 2\\[-20pt]
	2 & 2 & -1 & 1 & \boxed{-1}
\end{tikzcd}$ 

\textbf{Case 2.4.2.2.}  $2$ is located right from $a_4$: 
$\begin{tikzcd}[column sep=tiny]
	-1 & 1 \arrow[l,bend left=300,swap,"\frac{1}{5}"] \arrow[r,bend left=60,"\frac{8}{5}"] & -1 & 2 \arrow[l,bend left=300,swap,"\frac{8}{5}"] \arrow[r,bend left=60,"\frac{1}{5}"] & 1 & \geqslant 1 & a_5\\[-20pt]
	2 & 2 & -1 & 1 & \boxed{-1} & 2 & \geqslant 1
\end{tikzcd}$. Depending on the value of $a_5$, new subcases appear: 

\textbf{Case 2.4.2.2.1.} $a_5\geqslant 1$. We have  the case 
$\begin{tikzcd}[column sep=tiny]
	-1 & 1 \arrow[l,bend left=300,swap,"\frac{1}{5}"] \arrow[r,bend left=60,"\frac{8}{5}"] & -1 & 2 \arrow[l,bend left=300,swap,"\frac{8}{5}"] \arrow[r,bend left=60,"\frac{1}{5}"] & 1 & \geqslant 1 \arrow[l,bend left=300,swap,"1"]& \boxed{\geqslant 1}\\[-20pt]
	2 & 2 & -1 & 1 & -1 & 2 & \geqslant 1
\end{tikzcd}$,

\textbf{Case 2.4.2.2.2.} $a_5=-1$. Here we have three new cases and we establish CBDP for each. 

If $2$ is located left from $a_5$:
$\begin{tikzcd}[column sep=tiny]
	-1 & 1 \arrow[l,bend left=300,swap,"\frac{1}{5}"] \arrow[r,bend left=60,"\frac{8}{5}"] & -1 & 2 \arrow[l,bend left=300,swap,"\frac{8}{5}"] \arrow[r,bend left=60,"\frac{1}{5}"] & 1 & 2\arrow[l,bend left=300,swap,"1"] \arrow[r,bend left=60,"\frac{6}{5}"]& \boxed{-1}\\[-20pt]
	2 & 2 & -1 & 1 & -1 & 2 & 1
\end{tikzcd}$  

If $2$ is located below $a_5$:
$\begin{tikzcd}[column sep=tiny]
	-1 & 1 \arrow[l,bend left=300,swap,"\frac{1}{5}"] \arrow[r,bend left=60,"\frac{8}{5}"] & -1 & 2 \arrow[l,bend left=300,swap,"\frac{8}{5}"] \arrow[r,bend left=60,"\frac{1}{5}"] & 1 & 1\arrow[l,bend left=300,swap,"1"] \arrow[r,bend left=60,"\frac{1}{5}"]& \boxed{-1}\\[-20pt]
	2 & 2 & -1 & 1 & -1 & 2 & 2
\end{tikzcd}$ 

If $2$ is located right from $a_5$: 
$\begin{tikzcd}[column sep=tiny]
	-1 & 1 \arrow[l,bend left=300,swap,"\frac{1}{5}"] \arrow[r,bend left=60,"\frac{8}{5}"] & -1 & 2 \arrow[l,bend left=300,swap,"\frac{8}{5}"] \arrow[r,bend left=60,"\frac{1}{5}"] & 1 & 1\arrow[l,bend left=300,swap,"1"] \arrow[r,bend left=60,"\frac{4}{5}"]& \boxed{-1} & 2 \arrow[l,bend left=300,swap,"\frac{2}{5}"] \arrow[r,bend left=60,dashrightarrow,"\frac{6}{5}"] & \geqslant 1\\[-20pt]
	2 & 2 & -1 & 1 & -1 & 2 & 1  & \geqslant 1
\end{tikzcd}$.

$\bullet$ \textbf{Type (b)}; 

Due to the inherent symmetry of the problem, it suffices to consider only the following interior column configurations:  
$\begin{tikzcd}[column sep=tiny]
	-1\\[-20pt]
	-1 
\end{tikzcd}$, 
$\begin{tikzcd}[column sep=tiny]
	-1\\[-20pt]
	1 
\end{tikzcd}$ and 
$\begin{tikzcd}[column sep=tiny]
	-1\\[-20pt]
	2 
\end{tikzcd}$. 
Moreover, for each possible graph labeling, it is sufficient to analyze only one representative non-symmetrical instance, as symmetric configurations exhibit identical CBDP.

\textbf{Case 1.} CBDP for the interior column $\begin{tikzcd}[column sep=tiny]
	-1\\[-20pt]
	-1
\end{tikzcd}$. Depending on the position of $2$s, which protect these $-1$s, we distinguish two cases.

\textbf{Case 1.1.} If two $2$'s are located on the same side of two $-1$'s, we immediately can perform the following CBDP:  
$\begin{tikzcd}[column sep=tiny]
	\geqslant 1 &	1 \arrow[r,bend left=60,"\frac{2}{5}"] & -1 & 2 \arrow[l,bend left=300,swap,"\frac{14}{5}"] & \geqslant 1\\[-20pt]
	\geqslant 1 &	1& -1 & 2 & \geqslant 1
\end{tikzcd}$

\textbf{Case 1.2.} If two $2$'s are located on opposite sides of two $-1$'s, we get the configuration
$\begin{tikzcd}[column sep=tiny]
	b_1 &	1 & -1 & 2 & \geqslant 1\\[-20pt]
	\geqslant 1 &	2& -1 & 1 & b_2
\end{tikzcd}$.

Depending on values of $b_1$ and $b_2$, we have four subcases.

\textbf{Case 1.2.1.} $b_1,b_2\geqslant 1$. CBDP is defined as: 
$\begin{tikzcd}[column sep=tiny]
	\boxed{\geqslant 1} &	1 \arrow[r,bend left=60,"\frac{8}{5}"] & -1 & 2 \arrow[l,bend left=300,swap,"\frac{8}{5}"] & \geqslant 1\\[-20pt]
	\geqslant 1 &	2& -1 & 1 & 	\boxed{\geqslant 1} 
\end{tikzcd}$,

\textbf{Case 1.2.2.} $b_1=-1$, $b_2\geqslant 1.$ 

\textbf{Case 1.2.2.1.}  If $2$ is located below $b_1$: 
$\begin{tikzcd}[column sep=tiny]
	& 	\boxed{\geqslant -1} &	1 \arrow[l,bend left=300,swap,"\frac{1}{5}"]\arrow[r,bend left=60,"\frac{7}{5}"] & -1 & 2 \arrow[l,bend left=300,swap,"\frac{9}{5}"] & \geqslant 1\\[-20pt]
	\geqslant 1 & 2&	2& -1 & 1 & 	\boxed{\geqslant 1}
\end{tikzcd}$ 

\textbf{Case 1.2.2.2.}  If $2$ is located left of $b_1$, we get the configuration
$\begin{tikzcd}[column sep=tiny]
	\geqslant 1	&2 & \boxed{\geqslant -1} &	1 & -1 & 2  & \geqslant 1\\[-20pt]
	b_3 &\geqslant 1 & 1&	2& -1 & 1 & 	\boxed{\geqslant 1}
\end{tikzcd}$. Now depending on $b_3$, we consider:

\textbf{Case 1.2.2.2.1.} $b_3\geqslant 1$. CBDP is defined as: 
$\begin{tikzcd}[column sep=tiny]
	\geqslant 1 &2 \arrow[r,bend left=60,"\frac{4}{5}"] & -1 &	1 \arrow[l,bend left=300,swap,"\frac{2}{5}"]\arrow[r,bend left=60,"\frac{7}{5}"] & -1 & 2 \arrow[l,bend left=300,swap,"\frac{9}{5}"] & \geqslant 1\\[-20pt]
	\boxed{\geqslant 1}& \geqslant 1 & 1&	2& -1 & 1 & \geqslant 1
\end{tikzcd}$,

\textbf{Case 1.2.2.2.2.}  $b_3=-1$. If $2$ is located right of $b_3$: 
$\begin{tikzcd}[column sep=tiny]
	1 &2 \arrow[r,bend left=60,"\frac{4}{5}"] \arrow[l,bend left=300,swap,"\frac{6}{5}"]& -1 &	1 \arrow[l,bend left=300,swap,"\frac{2}{5}"]\arrow[r,bend left=60,"\frac{7}{5}"] & -1 & 2 \arrow[l,bend left=300,swap,"\frac{9}{5}"] & \geqslant 1\\[-20pt]
	\boxed{-1}& 2 & 1&	2& -1 & 1 & \geqslant 1
\end{tikzcd}$ 

If $2$ is located above $b_3$: 
$\begin{tikzcd}[column sep=tiny]
	2 &2 \arrow[r,bend left=60,"\frac{4}{5}"] \arrow[l,bend left=300,swap,"\frac{1}{5}"]& -1 &	1 \arrow[l,bend left=300,swap,"\frac{2}{5}"]\arrow[r,bend left=60,"\frac{7}{5}"] & -1 & 2 \arrow[l,bend left=300,swap,"\frac{9}{5}"] & \geqslant 1\\[-20pt]
	\boxed{-1}& 1 & 1&	2& -1 & 1 & \geqslant 1
\end{tikzcd}$

If $2$ is located left of $b_3$:
$\begin{tikzcd}[column sep=tiny]
	\phantom{-1}&\geqslant 1 \arrow[r,bend left=60,"\frac{1}{5}"] \arrow[l,bend left=300,dashrightarrow,swap,"\frac{6}{5}"]&1 &2 \arrow[r,bend left=60,"\frac{4}{5}"] \arrow[l,bend left=300,swap,"1"]& -1 &	1 \arrow[l,bend left=300,swap,"\frac{2}{5}"]\arrow[r,bend left=60,"\frac{7}{5}"] & -1 & 2 \arrow[l,bend left=300,swap,"\frac{9}{5}"] & \geqslant 1\\[-20pt]
	\geqslant 1&2 &\boxed{-1}& 1 & 1&	2& -1 & 1 & \geqslant 1
\end{tikzcd}$

\textbf{Case 1.2.3.}  $b_1\geqslant 1$, $b_2=-1$. This case is symmetrical to the case Case 1.2.2.

\textbf{Case 1.2.4.}  $b_1=b_2=-1$. Depending on position of $2$s, which protect $b_1$ and $b_2$, we distinguish three different non-symmetrical subcases:

\textbf{Case 1.2.4.1.}  If $2$'s are located below $b_1$ and above $b_2$, we immediately have:

$\begin{tikzcd}[column sep=tiny]
	\boxed{-1} &	1 \arrow[r,bend left=60,"\frac{8}{5}"] \arrow[l,bend left=300,swap,"\frac{1}{5}"] & -1 & 2 \arrow[l,bend left=300,swap,"\frac{8}{5}"] \arrow[r,bend left=60,"\frac{1}{5}"] & 2\\[-20pt]
	2 &	2& -1 & 1 & \boxed{-1}
\end{tikzcd}$

\textbf{Case 1.2.4.2.}  $2$'s are located left of $b_1$ and above $b_2$:
$\begin{tikzcd}[column sep=tiny]
	\geqslant 1 &2  & \boxed{-1} &	1  & -1 & 2  & 2\\[-20pt]
	b_3 & \geqslant 1&1 &	2& -1 & 1 & \boxed{-1}
\end{tikzcd}$

\textbf{Case 1.2.4.2.1.}
If $b_3\geqslant 1$, CBDP is defined as: 
$\begin{tikzcd}[column sep=tiny]
	\geqslant 1 &2 \arrow[r,bend left=60,"1"] & -1 &	1 \arrow[r,bend left=60,"\frac{8}{5}"] \arrow[l,bend left=300,swap,"\frac{1}{5}"]  & -1 & 2 \arrow[l,bend left=300,swap,"\frac{8}{5}"] \arrow[r,bend left=60,"\frac{1}{5}"] & 2\\[-20pt]
	\boxed{\geqslant 1} & \geqslant 1&1 &	2& -1 & 1 & -1
\end{tikzcd}$.

\textbf{Case 1.2.4.2.2.}
If $b_3=-1$, we have three sub-variants, depending on the position of 2, which protects $b_3$.

If $2$ is located right of $b_3$:
$\begin{tikzcd}[column sep=tiny]
	1 &2 \arrow[r,bend left=60,"1"] \arrow[l,bend left=300,swap,"\frac{6}{5}"]& -1 &	1 \arrow[r,bend left=60,"\frac{8}{5}"] \arrow[l,bend left=300,swap,"\frac{1}{5}"]  & -1 & 2 \arrow[l,bend left=300,swap,"\frac{8}{5}"] \arrow[r,bend left=60,"\frac{1}{5}"] & 2\\[-20pt]
	\boxed{-1} & 2 &1 &	2& -1 & 1 & -1
\end{tikzcd}$

If $2$ is located above $b_3$:
$\begin{tikzcd}[column sep=tiny]
	2 &2 \arrow[r,bend left=60,"1"] \arrow[l,bend left=300,swap,"\frac{1}{5}"]& -1 &	1 \arrow[r,bend left=60,"\frac{8}{5}"] \arrow[l,bend left=300,swap,"\frac{1}{5}"]  & -1 & 2 \arrow[l,bend left=300,swap,"\frac{8}{5}"] \arrow[r,bend left=60,"\frac{1}{5}"] & 2\\[-20pt]
	\boxed{-1} & 1 &1 &	2& -1 & 1 & -1
\end{tikzcd}$

If $2$ is located left of $b_3$:
$\begin{tikzcd}[column sep=tiny]
	\phantom{-1}&\geqslant 1 \arrow[r,bend left=60,"\frac{2}{5}"] \arrow[l,bend left=300,swap,dashrightarrow,"\frac{6}{5}"]&1 &2 \arrow[r,bend left=60,"1"] \arrow[l,bend left=300,swap,"\frac{4}{5}"]& -1 &	1 \arrow[r,bend left=60,"\frac{8}{5}"] \arrow[l,bend left=300,swap,"\frac{1}{5}"]  & -1 & 2 \arrow[l,bend left=300,swap,"\frac{8}{5}"] \arrow[r,bend left=60,"\frac{1}{5}"] & 2\\[-20pt]
	\geqslant 1&2&\boxed{-1}& 1 &1 &	2& -1 & 1 & -1
\end{tikzcd}$ .

\textbf{Case 1.2.4.3.} $2$'s are located left of $b_1$ and right of $b_2$:

$\begin{tikzcd}[column sep=tiny]
	\geqslant 1 &2 & \boxed{-1} & 1	 & -1 & 2  & 1 &\geqslant 1 & b_4\\[-20pt]
	b_3 & \geqslant 1&1 &	2& -1 & 1 & \boxed{-1} & 2 & \geqslant 1
\end{tikzcd}$.

Note that for this case,  there is a symmetry with respect to the column $\begin{tikzcd}[column sep=tiny]
	-1\\[-20pt]
	-1
\end{tikzcd}$. This implies that, in this case,  it suffices to describe the CBDP only for the columns to the left of this one. Consequently, we can apply the same types of CBDPs as described in the previous case.

\textbf{Case 2.} CBDP for the interior column $\begin{tikzcd}[column sep=tiny]
	-1\\[-20pt]
	1
\end{tikzcd}$. The following configuration is formed: 
$\begin{tikzcd}[column sep=tiny]
	-1 & 2 & \geqslant 1\\[-20pt]
	1 & \geqslant 1 & b_5
\end{tikzcd}$.

\textbf{Case 2.1.} $b_5\geqslant 1$. CBDP is defined as:
$\begin{tikzcd}[column sep=tiny]
	-1 & 2 \arrow[l,bend left=300,swap,"\frac{6}{5}"]& \geqslant 1\\[-20pt]
	1 & \geqslant 1 & \boxed{\geqslant 1}
\end{tikzcd}$,

\textbf{Case 2.2.}  $b_5=-1$. We have three subcases, which are immediately solved as follows.

\textbf{Case 2.2.1.} $2$ is located left of $b_5$:  $\begin{tikzcd}[column sep=tiny]
	-1 & 2 \arrow[l,bend left=300,swap,"\frac{6}{5}"] \arrow[r,bend left=60,"\frac{6}{5}"] & 1\\[-20pt]
	1 & 2 & \boxed{-1}
\end{tikzcd}$ 

\textbf{Case 2.2.2.} $2$ is located above $b_5$:   $\begin{tikzcd}[column sep=tiny]
	-1 & 2 \arrow[l,bend left=300,swap,"\frac{6}{5}"] \arrow[r,bend left=60,"\frac{1}{5}"] & 2\\[-20pt]
	1 & 1 &\boxed{-1}
\end{tikzcd}$

\textbf{Case 2.2.3.} $2$ is located right of $b_5$:   $\begin{tikzcd}[column sep=tiny]
	-1 & 2 \arrow[l,bend left=300,swap,"\frac{6}{5}"] \arrow[r,bend left=60,"\frac{3}{5}"] & 1 & \geqslant 1 \arrow[l,bend left=300,swap,"\frac{3}{5}"] \arrow[r,bend left=60,dashrightarrow,"\frac{6}{5}"]&\phantom{-1}\\[-20pt]
	1 & 1 & \boxed{-1} & 2 & \geqslant 1
\end{tikzcd}$

\textbf{Case 3.} CBDP for the interior column $\begin{tikzcd}[column sep=tiny]
	-1\\[-20pt]
	2
\end{tikzcd}$. In this case the configuration is as follows: 
$\begin{tikzcd}[column sep=tiny]
	b_6 & -1 & b_7 \\[-20pt]
	\geqslant 1 &2 & \geqslant 1 
\end{tikzcd}$; 

The values $b_6 = -1$ and $b_7 = -1$ cannot occur simultaneously, as this violates condition C3 from Definition~\ref{def:1} for the vertex between them. Furthermore, the CBDP corresponding to the case where $b_6=-1$ and $b_7\geqslant 1$ is reduced to the CBDP for the type (a) column $\begin{tikzcd}[column sep=tiny]
	-1\\[-20pt]
	2
\end{tikzcd}$, because this column can be only charged by the neighboring columns from its right-hand side.  Due to symmetry, the same argument applies to the CBDP for the case $b_6\geqslant 1$ and $b_7=-1$.
Therefore, the only remaining course of action is to outline the CBDP for the case where $b_6=b_7=1$:

\textbf{Case 3.1.} $b_6=b_7=1$. This situation gives the configuration
$\begin{tikzcd}[column sep=tiny]
	b_8& 1 & -1 & 1 & b_{10}\\[-20pt]
	b_9 &\geqslant 1 &2 & \geqslant 1 & b_{11}
\end{tikzcd}$, that  generates 9 subcases, depending on the values $b_8 + b_9$ and $b_{10}+ b_{11}$. Let us analyze one by one.

\textbf{Case 3.1.1.}  $b_8+b_9\geqslant 1$ and $b_{10}+b_{11}\geqslant 1$. CBDP is defined as: 
$\begin{tikzcd}[column sep=tiny]
	b_8& 1 \arrow[l,bend left=300,swap,dashrightarrow,"\frac{1}{5}"] \arrow[r,bend left=60,"\frac{1}{5}"] & -1 & 1 \arrow[r,bend left=60,,dashrightarrow,"\frac{1}{5}"]& b_{10}\\[-20pt]
	b_9 &\geqslant 1 &2 & \geqslant 1 & b_{11}
\end{tikzcd}$,

\textbf{Case 3.1.2.}  $b_8+b_9\geqslant 1$ and $b_{10}+b_{11}=0$. We have three variants, which are resolved as follows.

\textbf{Case 3.1.2.1.}  $b_{10}=-1$, $b_{11}=1$:
$\begin{tikzcd}[column sep=tiny]
	b_8& 1 \arrow[l,bend left=300,swap,dashrightarrow,"\frac{1}{5}"] \arrow[r,bend left=60,"\frac{1}{5}"] & -1 & 1 \arrow[r,bend left=60,"\frac{6}{5}"]& \boxed{-1}\\[-20pt]
	b_9 &\geqslant 1 &2 & 2 ^*& \boxed{1}
\end{tikzcd}$  Note that the label $2$, here and after marked with the asterisk sign, is forced because of the condition (C3) from Definition~\ref{def:1} for the vertex above it.

\textbf{Case 3.1.2.2.} $b_{10}=1,b_{11}=-1$, $2$ is located left of $b_{11}$:
$\begin{tikzcd}[column sep=tiny]
	b_8& 1 \arrow[l,bend left=300,swap,dashrightarrow,"\frac{1}{5}"] \arrow[r,bend left=60,"\frac{1}{5}"] & -1 & 1 \arrow[r,bend left=60,"\frac{6}{5}"]& \boxed{1}\\[-20pt]
	b_9 &\geqslant 1 &2 & 2 & \boxed{-1} 
\end{tikzcd}$ 

\textbf{Case 3.1.2.3.} $b_{10}=1,b_{11}=-1$, $2$ is located right of $b_{11}$:
$\begin{tikzcd}[column sep=tiny]
	b_8& 1 \arrow[l,bend left=300,swap,dashrightarrow,"\frac{1}{5}"] \arrow[r,bend left=60,"\frac{1}{5}"] & -1 & 1 \arrow[r,bend left=60,"\frac{4}{5}"]& \boxed{1} & \geqslant 1 \arrow[l,bend left=300,swap,"\frac{2}{5}"] \arrow[r,bend left=60,dashrightarrow,"\frac{6}{5}"] &\phantom{-1}\\[-20pt]
	b_9 &\geqslant 1 &2 & 1 & \boxed{-1} & 2 &\geqslant 1
\end{tikzcd}$

\textbf{Case 3.1.3.} $b_8+b_9\geqslant 1$ and $b_{10}+b_{11}=-2.$ Now, we have a configuration: 

$\begin{tikzcd}[column sep=tiny]
	b_8& 1 \arrow[l,bend left=300,swap,dashrightarrow,"\frac{1}{5}"] \arrow[r,bend left=60,"\frac{1}{5}"] & -1 & 1 & \boxed{-1} & 2 & \geqslant 1\\[-20pt]
	b_9 &\geqslant 1 &2 & 2^* & \boxed{-1} & 1 & b_{12}
\end{tikzcd}$;

Depending on the value of $b_{12}$, we have the following sub-variants:

\textbf{Case 3.1.3.1.} 
$b_{12}\geqslant 1$. We immediately have:  
$\begin{tikzcd}[column sep=tiny]
	b_8& 1 \arrow[l,bend left=300,swap,dashrightarrow,"\frac{1}{5}"] \arrow[r,bend left=60,"\frac{1}{5}"] & -1 & 1 \arrow[r,bend left=60,"\frac{8}{5}"]  & -1 & 2 \arrow[l,bend left=300,swap,"\frac{8}{5}"]& \geqslant 1\\[-20pt]
	b_9 &\geqslant 1 &2 & 2^* & -1 & 1 & \boxed{\geqslant 1}
\end{tikzcd}$,

\textbf{Case 3.1.3.2.} $b_{12}=-1$. Depending on the position of 2, which protects $b_{12}$, we have

\textbf{Case 3.1.3.2.1.} $2$ is located above $b_{12}$: $\begin{tikzcd}[column sep=tiny]
	b_8& 1 \arrow[l,bend left=300,dashrightarrow,swap,"\frac{1}{5}"] \arrow[r,bend left=60,"\frac{1}{5}"] & -1 & 1 \arrow[r,bend left=60,"\frac{8}{5}"]  & -1 & 2 \arrow[l,bend left=300,swap,"\frac{8}{5}"] \arrow[r,bend left=60,"\frac{1}{5}"]& 2\\[-20pt]
	b_9 &\geqslant 1 &2 & 2^* & -1 & 1 &\boxed{-1}
\end{tikzcd}$ 

\textbf{Case 3.1.3.2.2.} $2$ is located right of $b_{12}$: $\begin{tikzcd}[column sep=tiny]
	b_8& 1 \arrow[l,bend left=300,swap,dashrightarrow,"\frac{1}{5}"] \arrow[r,bend left=60,"\frac{1}{5}"] & -1 & 1 \arrow[r,bend left=60,"\frac{8}{5}"]  & -1 & 2 \arrow[l,bend left=300,swap,"\frac{8}{5}"] \arrow[r,bend left=60,"\frac{1}{5}"]& 1 & \geqslant 1 & b_{13}\\[-20pt]
	b_9 &\geqslant 1 &2 & 2^* & -1 & 1 & \boxed{-1} & 2 & \geqslant 1
\end{tikzcd}$ . Now, we have to analyze the configuration, depending on the value of $b_{13}$.

\textbf{Case 3.1.3.2.2.1.} $b_{13}\geqslant 1$. We immediately have
$\begin{tikzcd}[column sep=tiny]
	b_8& 1 \arrow[l,bend left=300,swap, dashrightarrow,"\frac{1}{5}"] \arrow[r,bend left=60,"\frac{1}{5}"] & -1 & 1 \arrow[r,bend left=60,"\frac{8}{5}"]  & -1 & 2 \arrow[l,bend left=300,swap,"\frac{8}{5}"] \arrow[r,bend left=60,"\frac{1}{5}"]& 1 & \geqslant 1 \arrow[l,bend left=300,swap,"1"]& \boxed{\geqslant 1}\\[-20pt]
	b_9 &\geqslant 1 &2 & 2^* & -1 & 1 & -1 & 2 & \geqslant 1
\end{tikzcd}$

\textbf{Case 3.1.3.2.2.2.} $b_{13}=-1$.  We have three sub-variants:

If $2$ is located left of $b_{13}$: 
$\begin{tikzcd}[column sep=tiny]
	b_8& 1 \arrow[l,bend left=300,swap,dashrightarrow,"\frac{1}{5}"] \arrow[r,bend left=60,"\frac{1}{5}"] & -1 & 1 \arrow[r,bend left=60,"\frac{8}{5}"]  & -1 & 2 \arrow[l,bend left=300,swap,"\frac{8}{5}"] \arrow[r,bend left=60,"\frac{1}{5}"]& 1 & 2 \arrow[l,bend left=300,swap,"1"] \arrow[r,bend left=60,"\frac{6}{5}"]& \boxed{-1}\\[-20pt]
	b_9 &\geqslant 1 &2 & 2^* & -1 & 1 & -1 & 2 & 1
\end{tikzcd}$ 

If $2$ is located below $b_{13}$:  $\begin{tikzcd}[column sep=tiny]
	b_8& 1 \arrow[l,bend left=300,swap,dashrightarrow,"\frac{1}{5}"] \arrow[r,bend left=60,"\frac{1}{5}"] & -1 & 1 \arrow[r,bend left=60,"\frac{8}{5}"]  & -1 & 2 \arrow[l,bend left=300,swap,"\frac{8}{5}"] \arrow[r,bend left=60,"\frac{1}{5}"]& 1 & 1 \arrow[l,bend left=300,swap,"1"] \arrow[r,bend left=60,"\frac{1}{5}"]& \boxed{-1}\\[-20pt]
	b_9 &\geqslant 1 &2 & 2^* & -1 & 1 & -1 & 2 & 2
\end{tikzcd}$

If $2$ is located right of $b_{13}$:
$\begin{tikzcd}[column sep=tiny]
	b_8& 1 \arrow[l,bend left=300,swap,dashrightarrow,"\frac{1}{5}"] \arrow[r,bend left=60,"\frac{1}{5}"] & -1 & 1 \arrow[r,bend left=60,"\frac{8}{5}"]  & -1 & 2 \arrow[l,bend left=300,swap,"\frac{8}{5}"] \arrow[r,bend left=60,"\frac{1}{5}"]& 1 & 1 \arrow[l,bend left=300,swap,"1"] \arrow[r,bend left=60,"\frac{4}{5}"]& \boxed{-1} & 2 \arrow[l,bend left=300,swap,"\frac{2}{5}"] \arrow[r,bend left=60,dashrightarrow,"\frac{6}{5}"]&\geqslant 1\\[-20pt]
	b_9 &\geqslant 1 &2 & 2^* & -1 & 1 & -1 & 2 & 1 &\geqslant 1 &\phantom{-1}
\end{tikzcd}$

\textbf{Case 3.1.4.} $b_8+b_9=0$, $b_{10}+b_{11}\geqslant 1$;

This case is symmetrical to Case 3.1.2. 

\textbf{Case 3.1.5.} $b_8+b_9=0=b_{10}+b_{11}$. Now we have 4 subcases.

\textbf{Case 3.1.5.1.} $b_8=b_{10}=-1$, $b_9=b_{11}=1$: $\begin{tikzcd}[column sep=tiny]
	\boxed{-1}& 1 \arrow[r,bend left=60,"\frac{1}{5}"] \arrow[l,bend left=300,swap,"\frac{6}{5}"]& -1 & 1 \arrow[r,bend left=60,"\frac{6}{5}"] & \boxed{-1}\\[-20pt]
	\boxed{1} & 2^* &2 & 2^* & \boxed{1}
\end{tikzcd}$

\textbf{Case 3.1.5.2.} $b_8=b_{11}=-1$, $b_9=b_{10}=1$. We distinguish two sub-variants.

\textbf{Case 3.1.5.2.1.} $2$ is located left of $b_{11}$: 
$\begin{tikzcd}[column sep=tiny]
	\boxed{-1}& 1 \arrow[r,bend left=60,"\frac{1}{5}"] \arrow[l,bend left=300,swap,"\frac{6}{5}"]& -1 & 1 \arrow[r,bend left=60,"\frac{6}{5}"] & \boxed{1}\\[-20pt]
	\boxed{1} & 2^* &2 & 2 & \boxed{-1}
\end{tikzcd}$

\textbf{Case 3.1.5.2.2.} $2$ is located right of $b_{11}$:

$\begin{tikzcd}[column sep=tiny]
	\boxed{-1}& 1 \arrow[r,bend left=60,"\frac{1}{5}"] \arrow[l,bend left=300,swap,"\frac{6}{5}"]& -1 & 1 \arrow[r,bend left=60,"\frac{4}{5}"] & \boxed{1} & \geqslant 1 \arrow[l,bend left=300,swap,"\frac{2}{5}"]\arrow[r,bend left=60,dashrightarrow,"\frac{6}{5}"]& \phantom{-1} \\[-20pt]
	\boxed{1} & 2^* &2 & 1 & \boxed{-1} & 2 & \geqslant 1
\end{tikzcd}$

\textbf{Case 3.1.5.3.} $b_8=b_{11}=1$, $b_9=b_{10}=-1$. Again, we distinguish two sub-variants.

\textbf{Case 3.1.5.3.1}  $2$ is located right of $b_9$:
$\begin{tikzcd}[column sep=tiny]
	\boxed{1}& 1 \arrow[r,bend left=60,"\frac{1}{5}"] \arrow[l,bend left=300,swap,"\frac{6}{5}"]& -1 & 1 \arrow[r,bend left=60,"\frac{6}{5}"] & \boxed{-1}  \\[-20pt]
	\boxed{-1} & 2 &2 & 2^* & \boxed{1} 
\end{tikzcd}$

\textbf{Case 3.1.5.3.2}  $2$ is located left of $b_9$: 
$\begin{tikzcd}[column sep=tiny]
	\phantom{-1}&\geqslant 1 \arrow[r,bend left=60,"\frac{3}{5}"] \arrow[l,bend left=300,swap,dashrightarrow,"\frac{6}{5}"]&\boxed{1}& 1 \arrow[r,bend left=60,"\frac{1}{5}"] \arrow[l,bend left=300,swap,"\frac{3}{5}"]& -1 & 1 \arrow[r,bend left=60,"\frac{6}{5}"] & \boxed{-1}  \\[-20pt]
	\geqslant 1 &2& \boxed{-1} & 1 &2 & 2^* & \boxed{1} 
\end{tikzcd}$

\textbf{Case 3.1.5.4.} $b_8=b_{10}=1$, $b_9=b_{11}=-1$. We distinguish four sub-variants.

\textbf{Case 3.1.5.4.1.}  Twos are located right of $b_9$ and left of $b_{11}$:

$\begin{tikzcd}[column sep=tiny]
	\boxed{1}& 1 \arrow[r,bend left=60,"\frac{1}{5}"] \arrow[l,bend left=300,swap,"\frac{6}{5}"]& -1 & 1 \arrow[r,bend left=60,"\frac{6}{5}"] & \boxed{1} \\[-20pt]
	\boxed{-1} & 2 &2 & 2 & \boxed{-1} 
\end{tikzcd}$

\textbf{Case 3.1.5.4.2.}  Twos are located left of $b_9$ and left of $b_{11}$:

$\begin{tikzcd}[column sep=tiny]
	\phantom{-1}&\geqslant 1 \arrow[r,bend left=60,"\frac{3}{5}"] \arrow[l,bend left=300,swap,dashrightarrow,"\frac{6}{5}"] &\boxed{1}& 1 \arrow[r,bend left=60,"\frac{1}{5}"] \arrow[l,bend left=300,swap,"\frac{3}{5}"]& -1 & 1 \arrow[r,bend left=60,"\frac{6}{5}"] & \boxed{1} \\[-20pt]
	\geqslant 1 &2&\boxed{-1} & 1 &2 & 2 & \boxed{-1} 
\end{tikzcd}$

\textbf{Case 3.1.5.4.3.} Twos are located left of $b_9$ and right of $b_{11}$:

$\begin{tikzcd}[column sep=tiny]
	\phantom{-1}&\geqslant 1 \arrow[r,bend left=60,"\frac{3}{5}"] \arrow[l,bend left=300,swap,dashrightarrow,"\frac{6}{5}"] &\boxed{1}& 1 \arrow[r,bend left=60,"\frac{1}{5}"] \arrow[l,bend left=300,swap,"\frac{3}{5}"]& -1 & 1 \arrow[r,bend left=60,"\frac{4}{5}"] & \boxed{1} &\geqslant 1  \arrow[r,bend left=60,dashrightarrow,"\frac{6}{5}"] \arrow[l,bend left=300,swap,"\frac{2}{5}"]&\phantom{-1}\\[-20pt]
	\geqslant 1 &2&\boxed{-1} & 1 &2 & 1 & \boxed{-1} & 2 & \geqslant 1
\end{tikzcd}$

\textbf{Case 3.1.5.4.4.} Twos are located right of $b_9$ and right of $b_{11}$:

$\begin{tikzcd}[column sep=tiny]
	\boxed{1}& 1 \arrow[r,bend left=60,"\frac{1}{5}"] \arrow[l,bend left=300,swap,"\frac{6}{5}"]& -1 & 1 \arrow[r,bend left=60,"\frac{4}{5}"] & \boxed{1} &\geqslant 1  \arrow[r,bend left=60,dashrightarrow,"\frac{6}{5}"] \arrow[l,bend left=300,swap,"\frac{2}{5}"]&\phantom{-1}\\[-20pt]
	\boxed{-1} & 2 &2 & 1 & \boxed{-1} & 2 & \geqslant 1
\end{tikzcd}$

\textbf{Case 3.1.6.}  $b_8+b_9=0, b_{10}+b_{11}=-2$;

In this case, for the columns to the right of column 
$\begin{tikzcd}[column sep=tiny]
	-1\\[-20pt]
	2
\end{tikzcd}$ we apply the same types of CBDP as described in Case 3.1.3. Therefore, it suffices to provide
the CBDP for the columns to the left of column 
$\begin{tikzcd}[column sep=tiny]
	-1\\[-20pt]
	2
\end{tikzcd}$, and we get three sub-variants, which can be immediately resolved as follows.


\textbf{Case 3.1.6.1.} 	 $b_8=b_{10}=b_{11}=-1$, $b_9=1$.

$
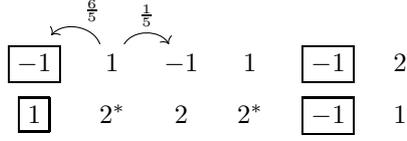
\begin{tikzcd}[column sep=tiny]
	\boxed{-1}& 1 \arrow[r,bend left=60,"\frac{1}{5}"] \arrow[l,bend left=300,swap,"\frac{6}{5}"]& -1 & 1  & \boxed{-1} & 2 \\[-20pt]
	\boxed{1} & 2^* &2 & 2^* & \boxed{-1} & 1
\end{tikzcd}$ 

\textbf{Case 3.1.6.2.} $b_8=1,b_9=b_{10}=b_{11}=-1$, $2$ is located right of $b_9$.

$
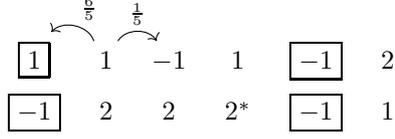
\begin{tikzcd}[column sep=tiny]
	\boxed{1}& 1 \arrow[r,bend left=60,"\frac{1}{5}"] \arrow[l,bend left=300,swap,"\frac{6}{5}"]& -1 & 1  & \boxed{-1} & 2 \\[-20pt]
	\boxed{-1} &  2 &2 & 2^* & \boxed{-1} & 1
\end{tikzcd}$

\textbf{Case 3.1.6.3.} $b_8=1,b_9=b_{10}=b_{11}=-1$, $2$ is located left of $b_9$.

$
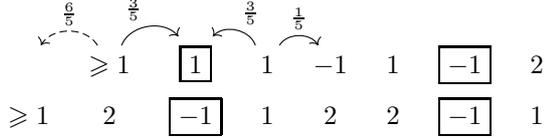
\begin{tikzcd}[column sep=tiny]
	\phantom{-1}&\geqslant 1 \arrow[r,bend left=60,"\frac{3}{5}"] \arrow[l,bend left=300,swap,dashrightarrow,"\frac{6}{5}"]& \boxed{1}& 1 \arrow[r,bend left=60,"\frac{1}{5}"] \arrow[l,bend left=300,swap,"\frac{3}{5}"]& -1 & 1  & \boxed{-1} & 2 \\[-20pt]
	\geqslant 1 &2&\boxed{-1} &  1 &2 & 2 & \boxed{-1} & 1
\end{tikzcd}$

\textbf{Case 3.1.7.} $b_8+b_9=-2, b_{10}+b_{11}\geqslant 1$. 
This case is symetrical to   Case 3.1.3.

\textbf{Case 3.1.8.} $b_8+b_9=-2, b_{10}+b_{11}=0$. 
This case is symetrical to   Case 3.1.6.

\textbf{Case 3.1.9.} $b_8+b_9=b_{10}+b_{11}=-2$. 

$
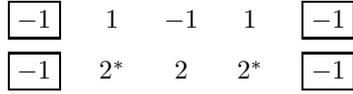
\begin{tikzcd}[column sep=tiny]
	\boxed{-1}& 1 & -1 & 1  & \boxed{-1}  \\[-20pt]
	\boxed{-1} &  2^* &2 & 2^* & \boxed{-1} 
\end{tikzcd}$

Once again, for the columns to the right of column 
$\begin{tikzcd}[column sep=tiny]
	-1\\[-20pt]
	2
\end{tikzcd}$, we apply the same types of CBDP as outlined in Case 3.1.3. Note that in this case,
column $\begin{tikzcd}[column sep=tiny]
	1\\[-20pt]
	2
\end{tikzcd}$, as a right-adjacent column to column $\begin{tikzcd}[column sep=tiny]
	-1\\[-20pt]
	2
\end{tikzcd}$, has an excess charge of $\frac{1}{5}$, which allows us to apply the same CBDP schemes described in Case 3.1.3, to define CBDP for the ``left'' portion of the graph.

Finally, after examining all possible cases, we conclude that a proper CBDP can be achieved, resulting in a charge of at least $\frac{6}{5}$ on each column. As stated in the main part of the paper at the conclusion of the proof of Theorem~\ref{thm:ladder}, this leads us to the inequality  $f(V(L_n)) \geqslant \left\lceil \frac{6n}{5} \right\rceil$.

\section{Detailed Analysis of Alternating Chains and Triplet Construction of Theorem~\ref{thm:grid3}}\label{app:B}

In the main body of the paper (proof of Theorem~\ref{thm:grid3}), we introduced the family $\mathcal C$ consisting of all alternating chains of maximal length with endpoints labeled -1, where neither endpoint has an adjacent vertex labeled with 1 outside of the chain. The process of ensuring the existence of a vertex labeled with 1 for each chain was briefly outlined. Here, we provide a comprehensive examination of the specific configurations and methods used to resolve these chains.

For each alternating chain $-1\,\,1-1\ldots-1\,\,1-1 \in \mathcal{C}$, the goal is to assign a label $1$ from a source external to the chain to construct valid triplets while maintaining the integrity of the labeling.

Before we dive into details of the proper constructions, three general observations can be made: 
\begin{enumerate}
	\item[O1.]  The endpoints of chain from $\mathcal C$ do not correspond to labels in the middle row of the grid graph, since vertices in $V_{-1}$ that belong to the middle row must have at least 2 neighbors from $V_1$.
	
	\item[O2.]  The chain from $\mathcal C$ cannot be placed in such a way that two consecutive $-1$s in this chain appear as diagonally placed labels in the graph. Indeed, the only non-symmetrical feasible labeling of the graph with this configuration is
	$\begin{tikzcd}[column sep=tiny]
		&2 & -1 & \boxed{1} & 2 \\[-20pt]
		\dotsc &\geqslant 1& \boxed{1} & -1 & &\dotsc\\[-20pt]
		&\phantom{-1} &\geqslant 1& \phantom{-1} & \phantom{-1}
	\end{tikzcd}$
		However,  since the endpoint of the given alternating chain is adjacent to two labels $1$ (boxed in the previous scheme), this chain does not belong to the family $\mathcal C$.
		
		\item[O3.] If the alternating chain, ending with a triplet which forms a column $\begin{tikzcd}[column sep=tiny]
			-1 & \\[-20pt]
			1\\[-20pt]
			-1
		\end{tikzcd}$, is adjacent to a column $\begin{tikzcd}[column sep=tiny]
			2 & \\[-20pt]
			\geqslant 1\\[-20pt]
			2
		\end{tikzcd}$, then label $1$ from one of the three mentioned source types can be found in the immediate neighborhood. Indeed, the list of all non-symmetrical feasible labelings containing
		$\begin{tikzcd}[column sep=tiny]
			-1 & 2 \\[-20pt]
			1 & \geqslant 1\\[-20pt]
			-1 & 2
		\end{tikzcd}$ is given by:
		
		$\begin{tikzcd}[column sep=tiny]
			-1 & 2 & \geqslant 1 \\[-20pt]
			1 & \boxed{\geqslant 1} & \geqslant 1\\[-20pt]
			-1 & 2 & \geqslant 1
		\end{tikzcd}$	(the boxed label is Source 1 type of label),
		
		$\begin{tikzcd}[column sep=tiny]
			-1 & 2 & 1 & \geqslant 1 \\[-20pt]
			1 & 2 & -1 &\phantom{-1}\\[-20pt]
			-1 & 2 & \boxed{1} & \geqslant 1
		\end{tikzcd}$	(the boxed label is Source 2 type of label),
		
		$\begin{tikzcd}[column sep=tiny]
			-1 & 2 & 1 & -1 & 2\\[-20pt]
			1 & 2 & -1 &1 &\geqslant 1\\[-20pt]
			-1 & 2 & \boxed{1} & \geqslant 1 &\phantom{-1}
		\end{tikzcd}$	(the boxed label is Source 3 type of label),
		
		$\begin{tikzcd}[column sep=tiny]
			-1 & 2 & 2 & \phantom{2} \\[-20pt]
			1 & \boxed{1} & -1 &\phantom{-1}\\[-20pt]
			-1 & 2 & 1 & \geqslant 1
		\end{tikzcd}$	(the boxed label is Source 2 type of label),
		
		$\begin{tikzcd}[column sep=tiny]
			-1 & 2 & 2 & \geqslant 1 &\phantom{-1}\\[-20pt]
			1 & \boxed{1} & -1 &1 &\geqslant 1\\[-20pt]
			-1 & 2 & 1 & -1 & 2
		\end{tikzcd}$	(the boxed label is Source 3 type of label),
		
		$\begin{tikzcd}[column sep=tiny]
			-1 & 2 & 1 & \geqslant 1 \\[-20pt]
			1 & \boxed{1} & -1 &2 \\[-20pt]
			-1 & 2 & 1 & \geqslant 1 
		\end{tikzcd}$	(the boxed label is Source 2 type of label).
		
		The similar procedure can be applied for any alternating chain with endpoint which is adjacent to a column $\begin{tikzcd}[column sep=tiny]
			2 & \\[-20pt]
			\geqslant 1\\[-20pt]
			2
		\end{tikzcd}$. 
	\end{enumerate}
	Now, we describe the appending procedure by examining five cases representing characteristic configurations of chains in $\mathcal C$. Let $l(C)$
	denote the length of the chain $C$. In the following schemes, and for the remainder of the proof, the considered chain is highlighted in bold.
	
	\textbf{Case 1.} The entire chain $C\in\mathcal C$ is located either in the first row or the third row of the graph labeling. Due to the symmetry, it is sufficient to assume that $C$ is located in the first row. We distinguish three possible cases.
	
	\textbf{Case 1.1.} Endpoints of $C$ are in the first and last column of the graph labeling (which is possible only for an odd number of columns).
	
	$\begin{tikzcd}[column sep=tiny]
		\textbf{-1} & \textbf{1} & \textbf{-1} & \textbf{1} &\dotsc & \textbf{-1} & \textbf{1} & \textbf{-1} \\[-20pt]
		2 & 2 & 2  &2 & \dotsc & 2 & 2 & 2\\[-20pt]
		\geqslant 1 & a_1 & \geqslant 1 & a_2 & \dotsc & \geqslant 1 & a_{\lfloor\frac{n}{2}\rfloor} & \geqslant 1
	\end{tikzcd}$;
	
	\textbf{Case 1.1.1.} $a_i\geqslant 1$ for at least one $i\in\{1,\dotsc,\lfloor\frac{n}{2}\rfloor\}$. The label $2$ directly above label $a_i$ is a Source 1 type label that can append the chain $C$. 
	
	\textbf{Case 1.1.2.}  $a_i = - 1$ for all $i\in\{1,\dotsc,\lfloor\frac{n}{2}\rfloor\}$.
	
	$\begin{tikzcd}[column sep=tiny]
		\textbf{-1} & \textbf{1} & \textbf{-1} & \textbf{1} &\dotsc & \textbf{-1} & \textbf{1} & \textbf{-1} \\[-20pt]
		2 & 2 & 2  &2 & \dotsc & 2 & 2 & 2\\[-20pt]
		\boxed{1} & -1 & 1 & -1 & \dotsc & 1 & -1 & 1
	\end{tikzcd}$, and the boxed label is a Source 3 type of label that can be added to the chain $C$.
	Notice that this configuration corresponds to the construction of a valid PLSRD function, which was initially presented at the beginning of the proof for the upper bound.
	
	\textbf{Case 1.2} Exactly one endpoint of $C$ is  in the leftmost (or rightmost) column of the graph. Depending on the position of label two which protects the last -1 in the chain, we distinguish two subcases.
	
	\textbf{Case 1.2.1} 2 is positioned to the right of the last -1.
	
	$\begin{tikzcd}[column sep=tiny]
		\textbf{-1} & \textbf{1} & \textbf{-1} & \textbf{1} &\dotsc & \textbf{1} & \textbf{-1} & \textbf{1} & \textbf{-1} & 2\\[-20pt]
		2 & 2 & 2  &2 & \dotsc & 2 & 2 & 2 & -1  & 1 &\dotsc\\[-20pt]
		\geqslant 1 & a_1 & \geqslant 1 & a_2 & \dotsc & a_{\lfloor\frac{l(C)}{2}\rfloor-1}& \geqslant 1 & \boxed{\geqslant 1 } &  1 
	\end{tikzcd}$ 
	
	The boxed label is a Source 1 type of label that can be appended the chain $C$.
	
	\textbf{Case 1.2.2}  2 is located bellow the last -1.
	
	$\begin{tikzcd}[column sep=tiny]
		\textbf{-1} & \textbf{1} & \textbf{-1} & \textbf{1} &\dotsc & \textbf{-1} & \textbf{1} & \textbf{-1} & -1\\[-20pt]
		2 & 2 & 2  &2 & \dotsc & 2 & 2 & 2  &\dotsc\\[-20pt]
		\geqslant 1 & a_1 & \geqslant 1 & a_2 & \dotsc & \geqslant 1 & a_{\lfloor\frac{l(C)}{2}\rfloor} & \geqslant 1  
	\end{tikzcd}$
	
	\textbf{Case 1.2.2.1} 	If $a_i\geqslant 1$ for at least one $i\in\{1,\dotsc,\lfloor\frac{l(C)}{2}\rfloor\}$, then the label $2$ directly above label $a_i$ is a Source 1 type label that can be appended to the chain $C$.

	\textbf{Case 1.2.2.2} Each $a_i=-1$, $i\in\{1,\dotsc,\lfloor\frac{l(C)}{2}\rfloor\}$. This results in the following configuration.
	
	$\begin{tikzcd}[column sep=tiny]
		\textbf{-1} & \textbf{1} & \textbf{-1} & \textbf{1} &\dotsc & \textbf{-1} & \textbf{1} & \textbf{-1} & -1\\[-20pt]
		2 & 2 & 2  &2 & \dotsc & 2 & 2 & 2  & \phantom{-1}&\dotsc\\[-20pt]
		1 & -1 & 1 & -1 & \dotsc & 1 & -1 & 1 & b 
	\end{tikzcd}$ 
	
	Depending on the value $b$ we have two new subcases.
	
	\textbf{Case 1.2.2.2.1} $b\geqslant 1$.
	
	$\begin{tikzcd}[column sep=tiny]
		\textbf{-1} & \textbf{1} & \textbf{-1} & \textbf{1} &\dotsc & \textbf{-1} & \textbf{1} & \textbf{-1} & -1\\[-20pt]
		2 & 2 & 2  &2 & \dotsc & 2 & 2 & 2  & \phantom{-1}&\dotsc\\[-20pt]
		\boxed{1} & -1 & 1 & -1 & \dotsc & 1 & -1 & 1 & \geqslant 1  
	\end{tikzcd}$
	
	The boxed label is Source 3 type label that can be appended to the chain $C$.
	
	\textbf{Case 1.2.2.2.1} $b=-1$. From the scheme below, it can be easily observed that another alternating ``L-shaped'' chain appears, containing equal number of -1s and 1s. The endpoints of this chain are marked with an asterisk ('*').  Consequently, in this case, we can apply  the general observation O3 and conclude that the label 1 can be identified immediately to the right of the last column on the right hand side.
	
	$\begin{tikzcd}[column sep=tiny]
		\textbf{-1} & \textbf{1} & \textbf{-1} & \textbf{1} &\dotsc & \textbf{-1} & \textbf{1} & \textbf{-1} & -1^* & 2 &\\[-20pt]
		2 & 2 & 2  &2 & \dotsc & 2 & 2 & 2  & 1&  \geqslant 1 &\dotsc\\[-20pt]
		1^* & -1 & 1 & -1 & \dotsc & 1 & -1 & 1 & -1  & 2 &
	\end{tikzcd}$ 
	
	\textbf{Case 1.3.} Neither endpoint of $C$ is in the leftmost (or rightmost) column of the graph. Depending on the location of 2s which protect the endpoints of $C$, we distinguish 2 cases.
	
	\textbf{Case 1.3.1} At least one endpoint is protected by a label $2$ located in the first row. WLOG, assume that the left $-1$ is protected by a label $2$ positioned to its left.
	
	$\begin{tikzcd}[column sep=tiny]
		&2	&\textbf{-1} & \textbf{1} & \textbf{-1} & \textbf{1} &\dotsc & \textbf{-1} & \textbf{1} & \textbf{-1} & \phantom{1}\\[-20pt]
		\dotsc& &-1& 2 & 2  &2 & \dotsc & 2 & 2 & b  &\dotsc\\[-20pt]
		&&\geqslant 1 & \boxed{\geqslant 1} & \geqslant 1 & a_2 & \dotsc & \geqslant 1 & a_{\lfloor\frac{l(C)}{2}\rfloor} &  
	\end{tikzcd}$
	
	and the boxed label is a Source 1 type label that can be appended to the chain $C$.	 
	
	\textbf{Case 1.3.2} Both endpoints are protected with 2 located  below each of them.
	
	$\begin{tikzcd}[column sep=tiny]
		&-1	&\textbf{-1} & \textbf{1} & \textbf{-1} & \textbf{1} &\dotsc & \textbf{-1} & \textbf{1} & \textbf{-1} & -1&\\[-20pt]
		\dotsc &&2& 2 & 2  &2 & \dotsc & 2 & 2 & 2  &&\dotsc\\[-20pt]
		&	&\geqslant 1 & a_1 & \geqslant 1 & a_2 & \dotsc & \geqslant 1 & a_{\lfloor\frac{l(C)}{2}\rfloor} &  &
	\end{tikzcd}$
	
	\textbf{Case 1.3.2.1} 
	$a_i\geqslant 1$ for at least one $i\in\{1,\dotsc,\lfloor\frac{l(C)}{2}\rfloor\}$.

	The label $2$ directly above label $a_i$ is a Source 1 type label that can be appended the chain $C$.

	\textbf{Case 1.3.2.2}  	$a_i= -1$ for all $i\in\{1,\dotsc,\lfloor\frac{l(C)}{2}\rfloor\}$.
	
	$\begin{tikzcd}[column sep=tiny]
		&&-1	&\textbf{-1} & \textbf{1} & \textbf{-1} & \textbf{1} &\dotsc & \textbf{-1} & \textbf{1} & \textbf{-1}  &-1\\[-20pt]
		\dotsc&\phantom{-1}&\phantom{-1} &2& 2 & 2  &2 & \dotsc & 2 & 2 & 2 &\phantom{-1} &\dotsc\\[-20pt]
		&&b_1&1 & -1 & 1 & -1 & \dotsc & -1 & -1 & 1 & b_2 
	\end{tikzcd}$
	
	Depending on the values of $b_1$ and $b_2$, we have three non symmetric sub-variants.
	
	\textbf{Case 1.3.2.2.1}   $b_1,b_2\geqslant 1$.
	
	$\begin{tikzcd}[column sep=tiny]
		&&-1	&\textbf{-1} & \textbf{1} & \textbf{-1} & \textbf{1} &\dotsc & \textbf{-1} & \textbf{1} & \textbf{-1}  &-1\\[-20pt]
		\dotsc&\phantom{-1}&\phantom{-1} &2& 2 & 2  &2 & \dotsc & 2 & 2 & 2 &\phantom{-1} &\dotsc\\[-20pt]
		&&\geqslant 1&\boxed{1} & -1 & 1 & -1 & \dotsc & -1 & -1 & 1 & \geqslant 1 
	\end{tikzcd}$ 
	
	The boxed label is Source 3 type label that can be  appended the chain $C$.
	
	\textbf{Case 1.3.2.2.2}  $b_1\geqslant 1$, $b_2=-1$, we get
	
	$\begin{tikzcd}[column sep=tiny]
		&&-1	&\textbf{-1} & \textbf{1} & \textbf{-1} & \textbf{1} &\dotsc & \textbf{-1} & \textbf{1} & \textbf{-1}  &-1 & 2\\[-20pt]
		\dotsc&\phantom{-1}&\phantom{-1} &2& 2 & 2  &2 & \dotsc & 2 & 2 & 2 &1 & \geqslant 1&\dotsc\\[-20pt]
		&&\geqslant 1&1 & -1 & 1 & -1 & \dotsc & -1 & -1 & 1 & -1 & 2
	\end{tikzcd}$, where general observation O3 can be applied. The same conclusion applies for $b_1=-1$, $b_2\geqslant 1$. 
	
	\textbf{Case 1.3.2.2.3} $b_1=b_2=-1$.
	
	$\begin{tikzcd}[column sep=tiny]
		&2&-1	&\textbf{-1} & \textbf{1} & \textbf{-1} & \textbf{1} &\dotsc & \textbf{-1} & \textbf{1} & \textbf{-1}  &-1 & 2\\[-20pt]
		\dotsc&\geqslant 1&1 &2& 2 & 2  &2 & \dotsc & 2 & 2 & 2 &1 & \geqslant 1&\dotsc\\[-20pt]
		&2&-1&1 & -1 & 1 & -1 & \dotsc & -1 & -1 & 1 & -1 & 2
	\end{tikzcd}$ and general observation O3 can be applied on both sides of the chain $C$.
	
	\textbf{Case 2.} The chain $C\in\mathcal C$ contains a $\bigsqcup$ shaped subchain $C'$. Let $l(C')$ denote the length of subchain $C'$. The subchain $C'$ of the chain $C$ is as follows:
	
	$\begin{tikzcd}[column sep=tiny]
		&\textbf{-1} & a_1 &  & a_2 & &\dotsc &  & a_{\lfloor\frac{l(C')-4}{2}\rfloor} & \textbf{-1} & \phantom{1}\\[-20pt]
		\dotsc &\textbf{1}& 2 & 2  &2 & 2& \dotsc & 2 & 2 & \textbf{1}  &\dotsc\\[-20pt]
		&\textbf{-1} & \textbf{1} & \textbf{-1} & \textbf{1} & \textbf{-1} & \dotsc & \textbf{-1} & \textbf{1} &  \textbf{-1}
	\end{tikzcd}$
	
	\textbf{Case 2.1.}
	$a_i\geqslant 1$ for at least one $i\in\{1,\dotsc,\lfloor\frac{l(C')-4}{2}\rfloor\}$. The label $2$ directly below label $a_i$ is a Source 1 type label that can be  appended to the chain $C$.

	\textbf{Case 2.2.} $a_i = -1$ for all $i\in\{1,\dotsc,\lfloor\frac{l(C')-4}{2}\rfloor\}$
	
	$\begin{tikzcd}[column sep=tiny]
		&2&\textbf{-1} & -1^* & 1 & -1 & 1&\dotsc & 1 & -1^*& \textbf{-1} & 2 &\\[-20pt]
		\dotsc &\geqslant 1&\textbf{1}& 2 & 2  &2 & 2& \dotsc & 2 & 2 & \textbf{1} &\geqslant 1 &\dotsc\\[-20pt]
		&2&\textbf{-1} & \textbf{1} & \textbf{-1} & \textbf{1} & \textbf{-1} & \dotsc & \textbf{-1} & \textbf{1} &  \textbf{-1} & 2
	\end{tikzcd}$
	
	and general observation O3 can be applied to both sides of the chain $C$ (one application for chain $C$ and the other for the
	$\mathcal C$-chain located in the first row of the graph, with the endpoints marked with the asterisk sign '*'.
	
	
	\textbf{Case 3.} The chain $C\in\mathcal C$ is L-shaped

	We distinguish two different configurations.
	
	\textbf{Case 3.1.} The chain spans to the end of graph. 
	
	$\begin{tikzcd}[column sep=tiny]
		& \textbf{-1} & a_1 & \geqslant 1 & a_2& & \geqslant 1&a_q & \geqslant 1 & &\\[-20pt]
		\dotsc & \textbf{1} & 2  &2 & 2& \dotsc & 2 & 2 & 2 &  \\[-20pt]
		& \textbf{-1} & \textbf{1} & \textbf{-1} & \textbf{1} &  & \textbf{-1} & \textbf{1}& \textbf{-1} & 
	\end{tikzcd}$ 
	
	\textbf{Case 3.1.1.}  $a_i\geqslant 1$ for at least one $i\in\{1,\dotsc,q\}$. Then,  the label $2$ directly below label $a_i$ is a Source 1 type label that can be  appended to the chain $C$.
	
	\textbf{Case 3.1.2.}  $a_i=-1$ for all $i\in\{1,\dotsc,q\}$.
	
	$\begin{tikzcd}[column sep=tiny]
		&2& \textbf{-1} & -1 & 1 & -1& & 1&-1 & 1 & &\\[-20pt]
		\dotsc & \geqslant 1 &\textbf{1} & 2  &2 & 2& \dotsc & 2 & 2 & 2 &  \\[-20pt]
		&2& \textbf{-1} & \textbf{1} & \textbf{-1} & \textbf{1} &  & \textbf{-1} & \textbf{1} & \textbf{-1} & 
	\end{tikzcd}$ and general observation O3 can be applied on the left  side of the chain $C$.
	
	\textbf{Case 3.2.} There are columns right to the end of the chain.
	
	$\begin{tikzcd}[column sep=tiny]
		& \textbf{-1} & a_1 & \geqslant 1 & a_2& & \geqslant 1&a_q & \geqslant 1 & a_{q+1}&\\[-20pt]
		\dotsc & \textbf{1} & 2  &2 & 2& \dotsc & 2 & 2 & 2 &  & \dotsc \\[-20pt]
		& \textbf{-1} & \textbf{1} & \textbf{-1} & \textbf{1} &  & \textbf{-1} & \textbf{1} & \textbf{-1} & -1
	\end{tikzcd}$
	
	\textbf{Case 3.2.1.} If $a_i\geqslant 1$ for at least one $i\in\{1,\dotsc,q\}$. Then, the label $2$ directly below label $a_i$ is a Source 1 type label that can be appended to the chain $C$. 
	
	\textbf{Case 3.2.2.}   $a_i=-1$ for all $i\in\{1,\dotsc,q\}$.
	
	$\begin{tikzcd}[column sep=tiny]
		&2& \textbf{-1} & -1 & 1 & -1& & 1&-1 & 1 & a_{q+1}&\\[-20pt]
		\dotsc &\geqslant 1& \textbf{1} & 2  &2 & 2& \dotsc & 2 & 2 & 2 &  & \dotsc \\[-20pt]
		&2& \textbf{-1} & \textbf{1} & \textbf{-1} & \textbf{1} &  & \textbf{-1} & \textbf{1} & \textbf{-1} & -1
	\end{tikzcd}$. 
	
	\textbf{Case 3.2.2.1.} $a_{q+1}\geqslant 1$.  We again apply general observation O3 on the left  side of the chain $C$.

	\textbf{Case 3.2.2.2.}  $a_{q+1}=-1$.
	
	$\begin{tikzcd}[column sep=tiny]
		&2& \textbf{-1} & -1^* & 1 & -1& & 1&-1 & 1 & -1& 2&\\[-20pt]
		\dotsc &\geqslant 1& \textbf{1} & 2  &2 & 2& \dotsc & 2 & 2 & 2 & 1 &\geqslant 1& \dotsc \\[-20pt]
		&2& \textbf{-1} & \textbf{1} & \textbf{-1} & \textbf{1} &  & \textbf{-1} & \textbf{1} & \textbf{-1} & -1^* & 2
	\end{tikzcd}$, and general observation O3 can be applied to both sides of the chain $C$ (one application for chain $C$ and 	
	the other for the other ``L-shaped'' chain from 	$\mathcal C$,  with the endpoints marked with the asterisk sign '*'

	\textbf{Case 4.} 
	The chain $C\in\mathcal C$ is of the shape \begin{tikzpicture}
		\draw[thick] (0,0) -- (0.15,0) -- (0.15,-0.3) -- (0.3,-0.3) ;  
	\end{tikzpicture} . This configuration is easy to resolve, as shown below.
	
	$\begin{tikzcd}[column sep=tiny]
		&\textbf{-1}& \textbf{1} & \textbf{-1} & 2& \geqslant 1 & \phantom{-1}&&\\[-20pt]
		\dotsc& & \phantom{-1}  &\textbf{1}  &\boxed{2} & 2& &&	\dotsc  \\[-20pt]
		&& 2 & \textbf{-1} & \textbf{1} & \textbf{-1} & \textbf{1} &\textbf{-1}&
	\end{tikzcd}$

	The boxed label 2  is a Source 1 type label that can be appended to the chain $C$.

	\textbf{Case 5.} The chain $C\in\mathcal C$ is positioned within a single column of the graph labeling,  i.e., the chain $C$ is of length 3 and is represented by a column $\begin{tikzcd}[column sep=tiny]
		-1 & \\[-20pt]
		1\\[-20pt]
		-1
	\end{tikzcd}$. 
	
	\textbf{Case 5.1.} $2$s are on the same side of $-1$s. In this case, the following configuration is forced:

	$\begin{tikzcd}[column sep=tiny]
		2& -1 & \textbf{-1} & 2\\[-20pt]
		\geqslant 1& 1 &\textbf{1} &\geqslant 1\\[-20pt]
		2& -1 & \textbf{-1} &2
	\end{tikzcd}$;

	From the scheme, it can be easily observed that general observation O3 is applied to both sides -- one application for chain $C$ and the other for the adjacent $\begin{tikzcd}[column sep=tiny]
		-1 & \\[-20pt]
		1\\[-20pt]
		-1
	\end{tikzcd}$$\mathcal C$-chain.

	\textbf{Case 5.2.} $2$s are on the opposite sides of $-1$s. The following configuration is formed:
	
	$\begin{tikzcd}[column sep=tiny]
		-1&\textbf{-1} & 2 &\geqslant 1 & a_1\\[-20pt]
		\geqslant 1& \textbf{1} &\geqslant 1&\geqslant 1 & b_1\\[-20pt]
		2&\textbf{-1} & -1 & \geqslant 1
	\end{tikzcd}$
	
	Depending on values $a_1$ and $b_1$, we distinguish two subcases.
	
	\textbf{Case 5.2.1.}  $a_1\geqslant 1$ or $b_1\geqslant 1$. The label $\geqslant 1$, positioned to the left of $a_1$ (respectively to the left of $b_1$), represents a Source 1 type label that can be appended to the chain $C$.

	\textbf{Case 5.2.2.} $a_1 = b_1=-1$. We have two subcases.
	
	\textbf{Case 5.2.2.1.}  2 is right from $b_1=-1$, ($b_1$ is marked with  asterisk sign '*').
	
	$\begin{tikzcd}[column sep=tiny]
		-1&\textbf{-1} & 2 &\geqslant 1 & -1 & \geqslant 1\\[-20pt]
		\geqslant 1& \textbf{1} &\geqslant 1&\boxed{1} & -1^* & 2\\[-20pt]
		2&\textbf{-1} & -1 & \geqslant 1 & 1 & \geqslant 1
	\end{tikzcd}$ 
	
	The boxed label is a Source 2 type label that can be appended to the chain $C$.

	\textbf{Case 5.2.2.2.} 2 is left or below $b_1 = -1$, ($b_1$ is marked with  asterisk sign '*').
	
	$\begin{tikzcd}[column sep=tiny]
		-1&\textbf{-1} & 2 &\geqslant 1 & -1 & \geqslant 1\\[-20pt]
		\geqslant 1& \textbf{1} &\geqslant 1&\geqslant 1 & -1^* & 1\\[-20pt]
		2&\textbf{-1} & -1 & \geqslant 1 & \geqslant 1 & c_1
	\end{tikzcd}$. 
	We distinguish two new subcases.
	
	\textbf{Case 5.2.2.2.1.}  $c_1\geqslant 1$, ($c_1$ is marked with  asterisk sign '*').
	
	$\begin{tikzcd}[column sep=tiny]
		-1&\textbf{-1} & 2 &\geqslant 1 & -1 & \geqslant 1 & \\[-20pt]
		\geqslant 1& \textbf{1} &\geqslant 1&\geqslant 1 & -1 & 1 & b_2\\[-20pt]
		2&\textbf{-1} & -1 & \geqslant 1 & \geqslant 1 & \geqslant 1^* & a_2
	\end{tikzcd}$;
	
	Depending on values $a_2$ and $b_2$, we have three sub-varaints.
	
	\textbf{Case 5.2.2.2.1.1.}
	$a_2\geqslant 1$ or 	$b_2\geqslant 1$.
	
	If $a_2\geqslant 1$:	
	$\begin{tikzcd}[column sep=tiny]
		-1&\textbf{-1} & 2 &\geqslant 1 & -1 & \geqslant 1 & \\[-20pt]
		\geqslant 1& \textbf{1} &\geqslant 1&\geqslant 1 & -1 & 1 & b_2\\[-20pt]
		2&\textbf{-1} & -1 & \geqslant 1 & \geqslant 1 & \boxed{\geqslant 1} & \geqslant 1
	\end{tikzcd}$ 
	
	The boxed label is a Source 1 type label that can be appended to the chain $C$.
	
	If $b_2\geqslant 1$:		
	$\begin{tikzcd}[column sep=tiny]
		-1&\textbf{-1} & 2 &\geqslant 1 & -1 & \geqslant 1 & \\[-20pt]
		\geqslant 1& \textbf{1} &\geqslant 1&\geqslant 1 & -1 & \boxed{1} & \geqslant 1\\[-20pt]
		2&\textbf{-1} & -1 & \geqslant 1 & \geqslant 1 & \geqslant 1 & a_2
	\end{tikzcd}$ 
	
	The boxed label is a Source 2 type label that can be appended  the chain $C$.
	
	\textbf{Case 5.2.2.2.1.2.} 
	$a_2=b_2=-1$.
	
	If 2 is located right from $b_2 = -1$ ($b_2$ is marked with asterisk sign '*'), we get the configuration 
	
	$\begin{tikzcd}[column sep=tiny]
		-1&\textbf{-1} & 2 &\geqslant 1 & -1 & \geqslant 1 & \boxed{1} &\geqslant 1\\[-20pt]
		\geqslant 1& \textbf{1} &\geqslant 1&\geqslant 1 & -1 & 1 & -1^* & 2\\[-20pt]
		2&\textbf{-1} & -1 & \geqslant 1 & \geqslant 1 & \geqslant 1 & -1
	\end{tikzcd}$ 
	
	and the boxed label is a Source 3 type label that can be appended the chain $C$.
	
	If 2 is located above $b_2 = -1 $ , we get
	
	$\begin{tikzcd}[column sep=tiny]
		-1&\textbf{-1} & 2 &\geqslant 1 & -1 & \geqslant 1 & 2 &\geqslant 1\\[-20pt]
		\geqslant 1& \textbf{1} &\geqslant 1&\geqslant 1 & -1 & 1 & -1^* & 1\\[-20pt]
		2&\textbf{-1} & -1 & \geqslant 1 & \geqslant 1 & \geqslant 1 & -1 &\geqslant 1
	\end{tikzcd}$. 
	
	Continuing with the described procedure, we obtain the alternating chain $C'$, placed in the central row of the graph labeling, starting with $b_1$. Denote by $b_m$ the other endpoint of chain $C'$. The label $b_m=-1$ cannot occupy the last column, as  neither of the columns $\begin{tikzcd}[column sep=tiny]
		-1 & \\[-20pt]
		-1\\[-20pt]
		*
	\end{tikzcd}$ nor $\begin{tikzcd}[column sep=tiny]
		* & \\[-20pt]
		-1\\[-20pt]
		-1
	\end{tikzcd}$ can occupy the last column in a valid graph labeling. 
	
	Depending on parity of $m$, we consider:
	
	$\bullet$ $m$ is even. We have the following configuration:
	
	$\begin{tikzcd}[column sep=tiny]
		&  &\geqslant 1 & -1 & \geqslant 1 & 2 &\geqslant 1 & & \geqslant 1 & \geqslant 1\\[-20pt]
		&\dotsc&\geqslant 1 & -1 & 1 & -1 & 1 &\dotsc & 1 & -1 & d\\[-20pt]
		&  & \geqslant 1 & \geqslant 1 & \geqslant 1 & -1 &\geqslant 1 & &\geqslant 1 &-1
	\end{tikzcd}$
	
	If $d=2$: 
	$\begin{tikzcd}[column sep=tiny]
		&  &\geqslant 1 & -1 & \geqslant 1 & 2 &\geqslant 1 & & \geqslant 1 & \boxed{1} & \geqslant 1\\[-20pt]
		&\dotsc&\geqslant 1 & -1 & 1 & -1 & 1 &\dotsc & 1 & -1 & 2\\[-20pt]
		&  & \geqslant 1 & \geqslant 1 & \geqslant 1 & -1 &\geqslant 1 & &\geqslant 1 &-1
	\end{tikzcd}$ 
	
	The boxed label is a Source 3 type label that can be appended to the chain $C$.
	
	If	$d=1$: 
	$\begin{tikzcd}[column sep=tiny]
		&  &\geqslant 1 & -1 & \geqslant 1 & 2 &\geqslant 1 & & \geqslant 1 & 2 & \geqslant 1\\[-20pt]
		&\dotsc&\geqslant 1 & -1 & 1 & -1 & 1 &\dotsc & 1 & -1 & \boxed{1} &\geqslant 1\\[-20pt]
		&  & \geqslant 1 & \geqslant 1 & \geqslant 1 & -1 &\geqslant 1 & &\geqslant 1 &-1 &\geqslant 1
	\end{tikzcd}$ 
	
	The boxed label is a Source 3 type label that can be  appended to the chain $C$.

	$\bullet$ 	$m$ is odd. We have the following configuration:
	
	$\begin{tikzcd}[column sep=tiny]
		&  &\geqslant 1 & -1 & \geqslant 1 & 2 &\geqslant 1 & & \geqslant 1 & -1\\[-20pt]
		&\dotsc&\geqslant 1 & -1 & 1 & -1 & 1 &\dotsc & 1 & -1 & e\\[-20pt]
		&  & \geqslant 1 & \geqslant 1 & \geqslant 1 & -1 &\geqslant 1 & &\geqslant 1 &\geqslant 1
	\end{tikzcd}$
	
	If $e=2$:
	$\begin{tikzcd}[column sep=tiny]
		&  &\geqslant 1 & -1 & \geqslant 1 & 2 &\geqslant 1 & & \geqslant 1 & -1 & \\[-20pt]
		&\dotsc&\geqslant 1 & -1 & 1 & -1 & 1 &\dotsc & 1 & -1 & 2\\[-20pt]
		&  & \geqslant 1 & \geqslant 1 & \geqslant 1 & -1 &\geqslant 1 & &\geqslant 1 &\boxed{1} & \geqslant 1
	\end{tikzcd}$ 
	
	The boxed label is a Source 3 type label that can be appended to the chain $C$.
	
	If $e=1$: 
	$\begin{tikzcd}[column sep=tiny]
		&  &\geqslant 1 & -1 & \geqslant 1 & 2 &\geqslant 1 & & \geqslant 1 & -1 & \geqslant 1\\[-20pt]
		&\dotsc&\geqslant 1 & -1 & 1 & -1 & 1 &\dotsc & 1 & -1 & \boxed{1} &\geqslant 1\\[-20pt]
		&  & \geqslant 1 & \geqslant 1 & \geqslant 1 & -1 &\geqslant 1 & &\geqslant 1 &2 &\geqslant 1
	\end{tikzcd}$ 
	
	The boxed label is a Source 3 type label that can be appended to the chain $C$.
	
	\textbf{Case 5.2.2.2.2.}  $c_1 = -1$ ($c1$ is marked with asterisk sign *).

	$\begin{tikzcd}[column sep=tiny]
		-1&\textbf{-1} & 2 &\geqslant 1 & -1 & \geqslant 1 &\geqslant 1 & a_1'\\[-20pt]
		\geqslant 1& \textbf{1} &1&2 & -1 & 1 & \geqslant 1 & b_1'\\[-20pt] 
		2&\textbf{-1} & -1 & 2 & 1 & -1^* & 2
	\end{tikzcd}$. 
	
	\textbf{Case 5.2.2.2.2.1.}  $a_1'\geqslant 1$ or $b_1'\geqslant 1$.
	Then, at least one of the  $\geqslant 1$ placed to the left of $a_1'$ and $b_1'$ is a Source 1 type label that can be appended to the chain $C$. 
	
	\textbf{Case 5.2.2.2.2.2.}  $a_1'=b_1'=-1$.

	If 2 is located right to $b_1'$ (marked with asterisk), we get the  configuration 
	
	$\begin{tikzcd}[column sep=tiny]
		-1&\textbf{-1} & 2 &\geqslant 1 & -1 & \geqslant 1 &\geqslant 1 & -1& \geqslant 1\\[-20pt]
		\geqslant 1& \textbf{1} &1&2 & -1 & 1 &1 & -1^*& 2\\[-20pt] 
		2&\textbf{-1} & -1 & 2 & 1 & -1 & 2 &\boxed{1} & \geqslant 1
	\end{tikzcd}$  and the boxed label is a Source 2 type label that can be appended to the chain $C$.
	
	If 2 is located left or below $b_1'$ (marked with asterisk), we get the  configuration 
	
	$\begin{tikzcd}[column sep=tiny]
		-1&-1 & 2 &\geqslant 1 & -1 & \geqslant 1 &\geqslant 1 & -1& \geqslant 1\\[-20pt]
		\geqslant 1& 1 &1&2 & -1 & 1 &\geqslant 1 & -1^*& 1\\[-20pt] 
		2&-1 & -1 & 2 & 1 & -1 & 2 &\geqslant 1 & c_2
	\end{tikzcd}$.

	If $c_2\geqslant 1$, we essentially switched to the setup of Case 4.2.2.2.1.
	
	If $c_2=-1$, we have 
	$\begin{tikzcd}[column sep=tiny]
		-1&-1 & 2 &\geqslant 1 & -1 & \geqslant 1 &\geqslant 1 & -1& \geqslant 1 &\geqslant 1\\[-20pt]
		\geqslant 1& 1 &1&2 & -1 & 1 &2& -1& 1 &\geqslant 1\\[-20pt] 
		2&-1 & -1 & 2 & 1 & -1 & 2 &1 & -1 & 2
	\end{tikzcd}$. Continuing with the described procedure, either it holds $c_l\geqslant 1$, for some $l$ in which case
	we continue by applying procedure described in Case 4.2.2.2.1., or the pattern $\begin{tikzcd}[column sep=tiny]
		-1 & \geqslant 1& \geqslant 1\\[-20pt]
		-1& 1 & \geqslant 1\\[-20pt]
		1 & -1 & 2 
	\end{tikzcd}$ repeats until it reaches the rightmost column of the graph labeling. In the last case scenario, pattern can end only with the column 
	$\begin{tikzcd}[column sep=tiny]
		\geqslant 1 & \\[-20pt]
		\geqslant 1\\[-20pt]
		2
	\end{tikzcd}$, so we obtain the configuration
	
	$\begin{tikzcd}[column sep=tiny]
		& &  &\geqslant 1 & -1 & \geqslant 1 &\geqslant 1 & -1& \geqslant 1 &\geqslant 1 & &-1 &\geqslant 1 & \boxed{\geqslant 1}\\[-20pt]
		& \dotsc & &2 & -1 & 1 &2& -1& 1 &2 & \dotsc & -1 & 1 & \geqslant 1\\[-20pt] 
		& &  & 2 & 1 & -1 & 2 &1 & -1 & 2 & & 1 & -1 & 2
	\end{tikzcd}$ 
	
	and the boxed label is a Source 1 type label that can be appended to the chain $C$.

	Having thoroughly analyzed the five basic configurations and presented a systematic method for identifying the label $1$ in each of them, we can now extend this reasoning to more complex configurations, i.e. more complex chains from $\mathcal{C}$. These complex configurations, which may arise as combinations of the basic ones, can be resolved by utilizing a label $1$ identified within one of their constituent basic configurations. In other words, the structure of such configurations ensures that at least one label $1$ corresponding to the original five configurations remains accessible and can be used to address the requirements of the extended structure.
	
	With all possible scenarios for the maximal alternating chains and the corresponding triplet constructions thoroughly examined, we can now conclude that the labeling process satisfies all necessary conditions thus completing the analysis in Appendix~\ref{app:B}.

\end{appendices}


\bibliography{bib}

\end{document}